\documentclass[reqno,11pt]{amsart}
\usepackage[top=2.0cm,bottom=2.0cm,left=3cm,right=3cm]{geometry}
\usepackage{amsthm,amsmath,amssymb,dsfont}
\usepackage{mathrsfs,amsfonts,functan,extarrows,mathtools}

\usepackage{xcolor}
\usepackage{cite}
\usepackage{indentfirst, latexsym, amssymb, enumerate,amsmath,graphicx}
\usepackage{float}
\usepackage{relsize}
\usepackage{marginnote}
\usepackage{stmaryrd}
\usepackage{esint}
\usepackage{graphicx}
\usepackage{bm}

\allowdisplaybreaks
\theoremstyle{plain}
\newtheorem{thm}{Theorem}[section]
\newtheorem{cor}[thm]{Corollary}
\newtheorem{lem}[thm]{Lemma}

\newtheorem{rem}{Remark}[section]

\numberwithin{equation}{section}

\begin{document}

\title[a fluid-particle coupled model with energy exchanges]
{Global existence and time decay of strong solutions to a fluid-particle coupled model with 
energy exchanges}

\author[F.-C. Li]{Fucai Li}
\address{School of Mathematics, Nanjing University, Nanjing
 210093, P. R. China}
\email{fli@nju.edu.cn}

\author[J.-K. Ni]{Jinkai Ni$^*$}\thanks{$^*$\! Corresponding author}
\address{School  of Mathematics, Nanjing University, Nanjing
 210093, P. R. China}
\email{602023210006@smail.nju.edu.cn}

\author[M. Wu]{Man Wu}
\address{Department of Mathematics, Nanjing Audit University, Nanjing
 211815, P. R. China}
\email{manwu@nau.edu.cn}

\begin{abstract}
In this paper, we investigate a three-dimensional fluid-particle coupled model. 
This model combines the full compressible Navier-Stokes equations with the Vlasov-Fokker-Planck equation
via the momentum and energy exchanges.   We  obtain the global existence  and  optimal time decay rates of strong solutions to the model in   whole space $\mathbb{R}^3$ when the initial data are a small perturbation of the given equilibrium in $H^2$. We show that the $L^2$-norms of the solutions and their gradients decay as $(1+t)^{-3/4}$ and $(1+t)^{-5/4}$ respectively. Moreover, we also obtain the decay rates of solutions in $L^p$-norms for $p\in [2,\infty]$,   and the optimal time decay rates of the highest-order derivatives of strong solutions  which reads as  $(1+t)^{-{7}/{4}}$ in $L^2$-norm.
When   the model is considered in a  periodic domain, besides the global existence results, we  show 
the strong solution decay exponentially.  Our proofs rely on the energy method,
  Fourier analysis techniques, and the method of frequency decomposition. And some new ideas are introduced to achieve the desired convergence rates.


\end{abstract}



\maketitle

\section{Introduction and main result}\label{Sec:intro-resul}
\setcounter{equation}{0}
 \indent \allowdisplaybreaks

Fluid-particle models have been  widely used in various fields such as 
combustion theory \cite{Wfa-1958,Wfa-1985},
chemical engineering \cite{CP-siam-1983},
dynamics of sprays \cite{BBBDLLT-irma-2005}, pharmaceutical industry \cite{Tu}, 
pollution treatment \cite{BWC-zamm-2000}, 
biomedical sprays \cite{BBJM-esaim-2005},  rain formation \cite{FFS-02}, and
sedimentation process \cite{BBKT.siam-2003}.
Due to their importance in applications, there are substantial mathematical results 
on fluid-particle models in the   past four decades. Among others, we mention 
\cite{BGM-jde-2017,CKL-jde-2011,
CJ-ccm-2023,FLR-jmfm-2021,GHMZ-siam-2010,Mn-jmpa-2011,SY-jde-2020,
SWYZ-jde-2023,WY-jde-2015,Yc-jmpa-2013}  on   incompressible fluid-particle models,
and \cite{MV-mmmas-2007,MV-cmp-2008, DL-krm-2013, CKL-jhde-2013,LMW-siam-2017,GHM-arma-2018}
on compressible fluid-particle models.

It should be pointed out that in the aforementioned mathematical results on incompressible and 
compressible models, the effect of temperature is ignored (the temperature is taking as a positive constant for simplicity).  However, according to the seminal work of Einstein \cite{Ei-1905}, the variety of temperature   play a crucial role in the modeling of dispersive two-phase  flows. 
In \cite{BBBGLLM-esaim-2009},  Boudin et al. introduced a fluid-particle model including a temperature 
equation of fluid. In their model the fluid is describing by non-isentropic Euler equations 
and  the particle by a  Vlasov-Fokker-Planck equation which  couple  together via 
the momentum and energy exchanges. 
They formally derived the model and its hydrodynamic limit, and prosed a numerical scheme 
to the limit system, see also \cite{GJY-cms-2012} for different numerical schemes on this model. 
By introducing the effects of viscosity  and the heat conductivity of the fluid, Mu and Wang \cite{MW-cvpde-2020} obtained the global well-posedness and large time behavior 
of classical solutions to this model  when the initial data are a small perturbation near the equilibrium state
in $H^4$. As pointed out in \cite{MW-cvpde-2020}, they  only considered the shear
viscosity of the fluid and skipped the bulk viscosity   for the sake of simplicity of presentation.
In this paper, we revisit  the model studied in \cite{MW-cvpde-2020} by considering the 
whole shear and bulk viscosities of the fluid since it is natural to include them  in
the momentum and  energy equations. Compared to the results in \cite{MW-cvpde-2020}, there are at least two 
advantages in our paper:  we lower the regularity of the initial data 
from $H^4$ to $H^2$; besides the decay rates of the solution themselves, we also 
obtain the  decay rates $(1+t)^{-5/4}$ of   their gradients,  and  $(1+t)^{-7/4}$
of their second-order derivatives in $x$ in the whole space $\mathbb{R}^3$ case,
 which are  consistent with  those of non-isentropic compressible Navier-Stokes equations. 
Furthermore, we obtain the decay of gradients of solutions in $L^p$ norm with $p\in [2,6]$, which 
are new even in the non-isentropic compressible Navier-Stokes equations.
  
%


More precisely, inspired by the models studied \cite{BBBGLLM-esaim-2009, GJY-cms-2012,MW-cvpde-2020}, in this paper we  study the following  non-isentropic fluid-particle coupled  model in the domain $\Omega$ ( $\Omega=\mathbb{R}^3$ or   $\Omega=\mathbb{T} ^3$, a periodic domain in  $\mathbb{R} ^3$):
\begin{equation}\label{I1}\left\{
\begin{aligned}
&\partial_t F+v\cdot\nabla_x F=\tilde{\mathbf{L}}f ,\\
&\partial_t n+{\rm div}_x (nu)=0,\\
&\partial_t (nu)+\nabla_x (nu\otimes u)-{\rm div}_x(\mathds{T})+\nabla_x P=\mathcal{M},\\
&\partial_t (nE)+{\rm div}_x \left((nE+P)u\right)
-{\rm div}_x(\kappa \nabla_x \tilde{\theta}) -{\rm div}_x(\mathds{T}u)=\mathcal{F} ,
\end{aligned}\right.\end{equation}
where the unknowns  $n=n(t,x)\geq0$,
$u=u(t,x)=(u_1(t,x),u_2(t,x),u_3(t,x))\in\mathbb{R}^3$,
$P=P(t,x)\geq0$,
$E=e+\frac{1}{2}|u|^2\geq0$, $e=e(t,x)\geq0$, and 
$\tilde{\theta}=\tilde{\theta} ( t, x) \geq 0$ 
for $( t, x) \in \mathbb{R} ^{+ }\times \Omega$ 
denote the density, velocity, pressure, total energy, internal energy,
and temperature of the fluid, respectively.
$F= F( t, x, v) \geq 0$ for $( t, x, v) 
\in \mathbb{R} ^+ \times \Omega\times \mathbb{R} ^3$
denotes the density distribution function of particles
in the phase space $\Omega\times \mathbb{R}^3$. 
$\mathds{T} $ is a stress tensor given by 
\begin{align}\label{ttt}
\mathds{T}=2\mu_1 D(u)+\mu_2   {\rm div}_x u\mathds{I}_3,
\end{align} 
where $D(u)=\frac{1}{2}(\nabla_x u+(\nabla_x u)^\top)$ denotes the deformation tensor,
$(\nabla_x u)^\top$ the transpose of $\nabla_x u$,
$\mathds{I}_3$   the $3\times 3$ identity matrix. 
The parameter $\mu_1 > 0$  denotes  the shear viscosity coefficient, 
and $\mu_2+\frac{2}{3}\mu_1$    the bulk viscosity coefficient. 
And $\kappa>0$ denotes the coefficient of thermal conductivity.
Generally speaking, the coefficients $\mu, \lambda$ and $\kappa$ may dependent on 
the density $\rho$ and the temperature $\tilde \theta$. In this paper, for simplicity, we assume that
they are constants independent of $\rho$ and $\tilde \theta$.  
For the   perfect gas case,  the  pressure $P$, temperature $\tilde{\theta}$,
and internal energy $e$ satisfy
the following relations:
\begin{align}
&  
 P=Rn\tilde{\theta} ,\quad 
 e=\frac{R}{ \gamma-1 }\tilde{\theta}=:C_\textsl{v} \tilde{\theta},
\end{align}
where $R>0$ is   the gas constant, and  $\gamma> 1$ is the adiabatic constant.
And $C_\textsl{v}>0$ is named as the heat conductive constant.

In the system \eqref{I1}, the fluid and the particles interacts mutually   via the 
Fokker-Planck type operator
\begin{align}\label{fpf}
{\tilde{\mathbf{L}}}f={\rm div}_x\big((v-u)F+\tilde{\theta}\nabla_v F\big),
\end{align} 
and the   coupling terms $\mathcal{M}$ and $\mathcal{F}$ in \eqref{I1}$_2$ and \eqref{I1}$_3$, defined as 
 \begin{align}\label{mf}
 \mathcal{M}=\int_{\mathbb{R}^3} (v-u)F\mathrm{d}v , \quad \
&\mathcal{F}=\int_{\mathbb{R}^3} [v\cdot(v-u)-3\tilde{\theta}]F\mathrm{d}v ,
\end{align}
indicate the momentum exchanges and 
 energy exchanges, respectively.
Thus, putting \eqref{ttt}--\eqref{mf} into \eqref{I1}, we can  rewrite it as follows:
\begin{equation}\label{I2}\left\{
\begin{aligned}
&\partial_t F+v\cdot\nabla_x F=\tilde{\mathbf{L}}f,\\
&\partial_t n+{\rm div}_x(nu)=0,\\
&\partial_t (nu)+\nabla_x (nu\otimes u)-\mu_1 \Delta_x u-(\mu_1 +\mu_2)\nabla_x({\rm div}_x  u)+\nabla_x P=\mathcal{M},\\
&C_\textsl{v}[\partial_t (n\tilde{\theta}+{\rm div}_x (nu\tilde{\theta})]+
{P {\rm  div} u} -
\kappa\Delta_x \tilde{\theta}-2\mu_1 |D(u)|^2-\mu_2|{\rm div}_x u|^2=\mathcal{N} ,
\end{aligned}\right.
\end{equation}
where 
$\mathcal{N}=\int_{\mathbb{R}^3} [|v-u|^2-3\tilde{\theta}]F\mathrm{d}v$.
We supplement the system \eqref{I2} with the initial data 
\begin{equation}\label{I--2}
(F,n,u,\tilde{\theta})|_{t=0}=\big(F_0(x,v),n_0(x),u_0(x),\tilde{\theta}_0(x)\big). 
\end{equation}

It is readily to verify that $(F, n, u, \tilde{\theta})\equiv
(M, 1, 0, 1)$ is an equilibrium state of the system \eqref{I2}, where 
$M$ is the normalized global Maxwellian:
\begin{align*}
M=M(v)=\frac{1}{(2\pi)^{3/2}}\exp\Big\{-\frac{|v|^{2}}{2}\Big\}.
\end{align*}
The purpose of this paper is to establish the global existence and large time behavior of   
strong solutions to the system \eqref{I2} when the initial data $(F_0,n_0,u_0,\tilde \theta_0)$ is 
a small perturbation of the given equilibrium state $(M, 1, 0, 1)$ in the whole space $\mathbb{R}^3$ 
and the torus $\mathbb{T}^3$, respectively.

 To do so, we first perturb  the system \eqref{I2} near the equilibrium  $(M, 1, 0, 1)$. We set  
$F=M+\sqrt{M}f$, $n=1+\rho$, $\tilde{\theta}=1+\theta$, $C_\textsl{v}=R=1$, and pull them into \eqref{I2} to 
get 
\begin{align}\label{G1}
&\partial_t f +v \cdot \nabla_x f+u \cdot \nabla_v f-\frac{1}{2} u \cdot v f-
u \cdot v \sqrt{M}-\left(|v|^2-3\right) \sqrt{M} \theta \nonumber\\
&\quad =\mathbf{L} f+\frac{\theta}{\sqrt{M}} \Delta_v(\sqrt{M} f), \\
\label{G2}
&\partial_t \rho +u \cdot \nabla_x \rho+(1+\rho) \operatorname{div} u=0, \\
\label{G3}
&\partial_t u +u \cdot \nabla_x u+\frac{1+\theta}{1+\rho} \nabla_x \rho+\nabla_x \theta=\frac{1}{1+\rho}\big(\mu_1\Delta_x u+(\mu_1+\mu_2)\nabla_x({\rm div}_x  u)\big)\nonumber\\
&\qquad+\frac{1}{1+\rho}\big(-u(1+a)+b\big), \\
\label{G4}&\partial_t \theta +u \cdot \nabla_x \theta+\theta \operatorname{div} u+\operatorname{div} u-\sqrt{6} \omega+3 \theta\nonumber \\
&\quad =\frac{1}{1+\rho}\left(\kappa\Delta_x \theta+|u|^2-2 u \cdot b+a|u|^2-3 a \theta+2\mu_1 |D(u)|^2+\mu_2|{\rm div}_x  u|^2\right)\nonumber \\
&\qquad  -\frac{\rho}{1+\rho}(\sqrt{6} \omega-3 \theta) .
\end{align}
Accordingly, the initial data to \eqref{G1}--\eqref{G4} reads
\begin{align}\label{G5}
\left.(f, \rho, u, \theta)\right|_{t=0} &\, =\left(\frac{F_0(x,v)-M(v)}{\sqrt{M(v)}}, n_0(x)-1, u_0(x), \tilde{\theta}_0(x)-1\right)\nonumber\\
&\,=:\left(f_0(x, v), \rho_0(x), u_0(x), \theta_0(x)\right).
\end{align}
In \eqref{G1}, we denote the linearized Fokker-Planck operator $\mathbf{L}$ by
$$
\mathbf{L}f=\frac{1}{\sqrt{M}} \nabla_v \cdot\left[M \nabla_v\left(\frac{f}{\sqrt{M}}\right)\right],
$$
and $a=a^f, b=b^f$ and $\omega=\omega^f$ are defined in the following way:
\begin{align*}
a^f(t, x) & \,=\int_{\mathbb{R}^3} \sqrt{M(v)} f(t, x, v) \mathrm{d} v, \\
b^f(t, x) & \,=\int_{\mathbb{R}^3} v \sqrt{M(v)} f(t, x, v) \mathrm{d} v, \\
\omega^f(t, x) & \,=\int_{\mathbb{R}^3} \frac{|v|^2-3}{\sqrt{6}} \sqrt{M(v)} f(t, x, v) \mathrm{d} v .
\end{align*}

Before stating our results, we need to give some nations.
We define the weight function $\nu(v)=1+|v|^2$ and the norm $|\cdot|_\nu$  as
\begin{align*}
|g|_\nu^2=\int_{\mathbb{R}^3}\big(\left|\nabla_v g(v)\right|^2+\nu(v)|g(v)|^2\big)\mathrm{d} v, \quad g=g(v).
\end{align*}
We use $\langle\cdot, \cdot\rangle$ to denote the inner 
product of the space $L_v^2:=L^2(\mathbb{R}^3_v)$, i.e.
\begin{align*}
\langle g, h\rangle=\int_{\mathbb{R}^3} g(v) h(v) \mathrm{d} v, \quad \forall g, h \in L_v^2 .
\end{align*}
For brevity, we use $\Omega$ to denote $\mathbb{R}^3_x$ or $\mathbb{T}^3_x$ and  define $\|\cdot\|$ the norm of $L_x^2=L^2(\Omega)$ or $L_{x, v}^2=L^2(\Omega\times \mathbb{R}^3_v)$:
\begin{align*}
\|g\|_\nu^2=\iint_{\Omega \times \mathbb{R}^3}\big(\left|\nabla_v g(x, v)\right|^2+\nu(v)|g(x, v)|^2\big) \mathrm{d} x \mathrm{d} v, \quad g=g(x, v) .
\end{align*}
For   $q \geq 1$, we also define 
\begin{align*}
Z_q=L_v^2\left(L_x^q\right)=L^2\left(\mathbb{R}_v^3 ; L^q\left(\mathbb{R}_x^3\right)\right), \quad\|g\|_{Z_q}^2=\int_{\mathbb{R}^3}\left(\int_{\mathbb{R}^3}|g(x, v)|^q \mathrm{d} x\right)^{\frac{2}{q}} \mathrm{d} v,
\end{align*}
and the norm $\|(\cdot, \cdot, \cdot, \cdot)\|_{\mathcal{Z}_q}$ reads 
\begin{align*}
\|(f, \rho, u, \theta)\|_{\mathcal{Z}_q}=\|f\|_{Z_q}+\|(\rho, u, \theta)\|_{L^1_x}.
\end{align*}
For an integrable function 
$g: \mathbb{R}^3 \rightarrow \mathbb{R}$, 
its Fourier transform can be defined by
\begin{align*}
\hat{g}(k)=\mathcal{F} g(k)=\int_{\mathbb{R}^3} e^{-i x \cdot k} g(x) \mathrm{d} x, \quad x \cdot k=\sum_{j=1}^3 x_j k_j,
\end{align*}
for $k \in \mathbb{R}^3$. {$(f|h)$ denotes the dot product of $f$ with complex conjugate of $ h $.}
For the multi-indices $\alpha=\left(\alpha_1, \alpha_2, \alpha_3\right)$ and $\beta=\left(\beta_1, \beta_2, \beta_3\right)$, we denote 
\begin{align*}
\partial_\beta^\alpha = \partial_{x_1}^{\alpha_1} \partial_{x_2}^{\alpha_2} \partial_{x_3}^{\alpha_3} \partial_{v_1}^{\beta_1} \partial_{v_2}^{\beta_2} \partial_{v_3}^{\beta_3},
\end{align*}
the partial derivatives with regard to $x=\left(x_1, x_2, x_3\right)$ and $v=\left(v_1, v_2, v_3\right)$. The length of $\alpha$ and $\beta$ are defined as $|\alpha|=\alpha_1+\alpha_2+\alpha_3$ and $|\beta|=\beta_1+\beta_2+\beta_3$. Define
\begin{align*}
\|g\|_{H^s}=\sum_{|\alpha| \leq s}\|\partial^\alpha g\|, \quad\|g\|_{H_{x, v}^s}=\sum_{|\alpha|+|\beta| \leq s}\|\partial_\beta^\alpha g\| .
\end{align*}

In this paper, we use $C$ to denote a positive (generally large) constant.
Additionally, the symbol $A \sim B$ denotes $ C A \leq B \leq\frac{1}{C} A$ for a generic constant $C>0$.

Motivated by \cite{DFT-cmp-2010, GY-iumj-2004,MW-cvpde-2020}, we introduce the  macro-micro decomposition to $f(t,x,v)$, which plays a crucial role in the  following estimates.  Consider  the  velocity orthogonal projection $\mathbf{P}$ as
\begin{equation*}
\mathbf{P}:L_{v}^{2}\rightarrow \mathrm{Span}\{\sqrt{M}, v_{1}\sqrt{M}, v_{2}\sqrt{M}, v_{3}\sqrt{M}, |v|^2\sqrt{M} \},
\end{equation*}
and
\begin{equation*}
\mathbf{P}=\mathbf{P}_{0}\oplus \mathbf{P}_{1}\oplus\mathbf{P}_{2},\quad \mathbf{P}_{0}f=a\sqrt M,\quad \mathbf{P}_{1}f=b\cdot v\sqrt M,\quad
\mathbf{P}_{2}f=\omega \frac{|v|^2-3}{\sqrt{6}} \sqrt{M}.
\end{equation*}
Next, we   decompose $f(t,x,v)$ as 
\begin{equation}\label{a}
f(t,x,v)=\mathbf{P}f+\{\mathbf{I}-\mathbf{P}\}f,
\end{equation}
where $\mathbf{P}f$ is the macro part of $f$ and $\{\mathbf{I}-\mathbf{P}\}f$   the micro part of it. 
Additionally, we   decompose $\mathbf{L}f$ as
\begin{equation*}
\mathbf{L}f=\mathbf{L}\{\mathbf{I}-\mathbf{P}\}f+\mathbf{L}\mathbf{P}f
=\mathbf{L}\{\mathbf{I}-\mathbf{P}\}f-\mathbf{P}_{1}f-\mathbf{P}_{2}f.
\end{equation*}
Using the property of $\mathbf{L}$, it is direct to verify that  there exist a positive constant $\lambda_0>0$ such that
\begin{align}
-\langle f,\mathbf{L}\{\mathbf{I}-\mathbf{P}\}f\rangle &\geq\lambda_{0}|\{\mathbf{I}-\mathbf{P}\}f|_{\nu}^{2},\label{G1.121}\\
-\langle f, \mathbf{L}f\rangle
&\geq\lambda_{0}|\{\mathbf{I}-\mathbf{P}\}f|_{\nu}^{2}+|b|^{2}+2|\omega|^2.\label{G1.12}
\end{align} 

Now we    present our main results as follows. First, for the case   $\Omega=\mathbb{R}^3$, we have 
\begin{thm}\label{T1.1}
Assume that  the initial data  $\left(f_0, \rho_0, u_0, \theta_0\right)$
satisfy   that $F_0=M+$ $\sqrt{M} f_0 \geq 0$, 
and  $\left\|f_0\right\|_{H_{x, v}^2(\mathbb{R}^3\times \mathbb{R}^3)}+\left\|\left(\rho_0, u_0, \theta_0\right)\right\|_{H^2(\mathbb{R}^3)}<\epsilon_0$
for some  $\epsilon_0>0$. Then the Cauchy problem \eqref{G1}--\eqref{G5} admits a unique global solution 
  $(f, \rho, u, \theta)$ satisfying $F=M+\sqrt{M} f \geq$ 0 and
\begin{gather*}
   f \in C\left([0, \infty) ; H_{x,v}^2(\mathbb{R}^3\times  \mathbb{R}^3 )\right), \quad\rho, u, \theta \in C\left([0, \infty) ; H^2(\mathbb{R}^3)\right), \\
  \sup _{t \geq 0}\big\{\|f(t)\|_{H_{x, v}^2}+\|(\rho, u, \theta)(t)\|_{H^2}\big\} \leq C\big(\left\|f_0\right\|_{H_{x, v}^2}+\left\|(\rho_0, u_0, \theta_0)\right\|_{H^2}\big),
\end{gather*}
for some constant $C>0$. Moreover, if we further assume that
\begin{equation*}
\left\|\left(f_0, \rho_0, u_0, \theta_0\right)\right\|_{\mathcal{Z}_1} \leq \epsilon_1, 
\end{equation*}
for some  small $\epsilon_1>0$, then
\begin{align}
\|f(t)\|_{{L_v^2(H_x^{2})}}+\|(\rho, u, \theta)(t)\|_{H^2(\mathbb{R}^3)} 
&\leq C(1+t)^{-\frac{3}{4}},\label{G1.13}\\
\|\nabla_x f(t)\|_{{L_v^2(H_x^{1})}}+\|\nabla_x(\rho, u, \theta)(t)\|_{H^1(\mathbb{R}^3)} &\leq C(1+t)^{-\frac{5}{4}},\label{G1.14}\\
 {\|\nabla_x^2 f(t)\|_{L_v^2(L_{x}^{2})}+\|\nabla_x^2(\rho, u, \theta)(t)\|_{L^2(\mathbb{R}^3)} }&\leq C(1+t)^{-\frac{7}{4}},\label{NJK1.17}
\end{align} 
for some constant $C>0$ and all $t \geq 0$.
\end{thm}

\begin{cor}\label{Cor1}
  Under the assumptions in Theorem \ref{T1.1},  the following decay estimates
    \begin{align}
\|(\rho,u,\theta)(t)\|_{L^p}&\leq C(1+t)^{-\frac{3}{2}(1-\frac{1}{p})} , \\
\|f(t)\|_{L_v^2(L_x^{p})}&\leq C(1+t)^{-\frac{3}{2}(1-\frac{1}{p})}, 
\end{align} 
  hold for $p\in [2,\infty]$. {
  Furthermore, 
\begin{align}
\|\nabla_x(\rho,u,\theta)\|_{L^p}
&\,\leq C(1+t)^{-\frac{3}{2}(\frac{4}{3}-\frac{1}{p})} , \label{glpa}\\
\|\nabla_x f\|_{L_v^2(L_x^{p})}
&\,\leq C(1+t)^{-\frac{3}{2}(\frac{4}{3}-\frac{1}{p})},\label{glpb}
\end{align}
hold for $p\in [2,6]$.
}
\end{cor}

\begin{rem}
  In Theorem \ref{T1.1} and Corollary \ref{Cor1}, we have obtained the optimal decay rates up to the highest-order  derivatives of the strong solutions
 to the non-isentropic fluid-particle model \eqref{G1}--\eqref{G4}. 
To the best of our knowledge, our  decay results in general $L^p$ norm with $p\in [6, +\infty]$ of the solution 
and  in $L^q$ norm with $q\in [2, 6]$ of the gradient of the solution 
 are  new even for the compressible Navier-Stokes equations.  It is very interesting to give the lower bounds of decay to the solutions to the system  \eqref{G1}--\eqref{G4} where more delicate spectral analysis are involved,  and we shall report it in a forthcoming paper. 
\end{rem}


\medskip 
For the case   $\Omega=\mathbb{T}^3$, we have 
\begin{thm}\label{T1.2}
Assume that  the initial data  $\left(f_0, \rho_0, u_0, \theta_0\right)$
satisfy   that $F_0=M+$ $\sqrt{M} f_0 \geq 0$, 
  $\left\|f_0\right\|_{H_{x, v}^2(\mathbb{T}^3\times \mathbb{R}^3)}+\left\|\left(\rho_0, u_0, \theta_0\right)\right\|_{H^2(\mathbb{T}^3 )}<\epsilon_2$
for some  $\epsilon_2>0$, and
\begin{gather*}
  \int_{\mathbb{T}^3} a_0 \mathrm{d} x=0, \quad \int_{\mathbb{T}^3} \rho_0 \mathrm{d} x=0, \\
  \int_{\mathbb{T}^3}\left(b_0+\left(1+\rho_0\right) u_0\right) \mathrm{d} x=0, \\
  \int_{\mathbb{T}^3}\left(1+\rho_0\right)\Big(\theta_0+\frac{1}{2}\left|u_0\right|^2\Big)+\frac{\sqrt{6}}{2} \omega_0 \mathrm{d} x=0,
\end{gather*}
where
\begin{align*}
a_0=\int_{\mathbb{T}^3} \sqrt{M} f_0(x, v) \mathrm{d} v,\ \  b_0=\int_{\mathbb{T}^3} v \sqrt{M} f_0(x, v) \mathrm{d} v,\ \  \omega_0=\int_{\mathbb{T}^3} \frac{|v|^2-3}{\sqrt{6}} \sqrt{M} f_0(x, v) \mathrm{d} v .
\end{align*}
Then, the Cauchy problem \eqref{G1}--\eqref{G5} admits a unique global solution $(f, \rho, u,\theta)$ satisfying $F=M+\sqrt{M} f \geq 0$, and
\begin{gather*}
  f \in C\left([0, \infty) ; H_{x,v}^2(\mathbb{T}^3\times \mathbb{R}^3)\right), \quad\rho, u, \theta \in C\left([0, \infty) ; H^2 (\mathbb{T}^3 )\right), \\
  \|f(t)\|_{H_{x, v}^2}+\|(\rho, u, \theta)(t)\|_{H^2} \leq C\big(\left\|f_0\right\|_{H_{x, v}^2}+\left\|(\rho_0, u_0, \theta_0)\right\|_{H^2}\big) e^{-\zeta_0 t}, 
\end{gather*}
with some constant $\zeta_0>0$, for all $t \geq 0$.
\end{thm}

Now we describe our strategies on the proofs of main results. For the existence
part of our results, we mainly follow the classical wildly-used energy method \cite{MN-jmku-1980, LMW-siam-2017, MW-cvpde-2020} to compressible fluid and fluid-particle models. This method  
usually include three steps: (a) constructing local solutions to a  system, (b) deriving uniform-in-time a priori estimates of solutions, and  (c) applying a continuity argument to extend the solutions from locally  to globally. 
Compared to the results in \cite{MN-jmku-1980, LMW-siam-2017, MW-cvpde-2020} where the  $H^4$ regularity of 
  initial data is needed, here we  only assume that the initial data have  $H^2$ regularity.
Thus our solutions are only $H^2$ regular.  Due to the lower regularity of solutions, we need to choose the Sobolev space techniques very carefully and
introduce some new ideas in  the arguments of deriving 
the uniform-in-time a priori estimates. As mentioned in \cite{MW-cvpde-2020}, we also need to 
introduce subtle micro-macro decomposition of $f$ to extract the damping effect in Fokker-Planck operator 
to balance the effect caused by the exchanges of momentum and energy changes between the fluid 
and particles. Finally, we  construct delicate Lyapunov-type energy   functional inequality to obtain the uniform-in-time estimates. 
To obtain the decay rates of solutions,  we employ  our global existence results and adapt the 
the energy-spectrum method developed in \cite{duyz-2008}.
The main impediments  come from the nonlinear 
terms of the system \eqref{G1}--\eqref{G5}, 
which  cause some difficulties in our analysis.
We ingeniously choose the decomposition of the  nonlinear items,
and then combine  them with the Duhamel's principle and energy estimates to obtain the decay rates of solutions.
To get further   decay results of the gradients of solutions, we introduce high-order energy functional and 
higher-order dissipation rate  (see \eqref{hh} and \eqref{mm} below) and derive 
a Lyapunov-type   inequality of them. 
 However, the optimal time-decay rates of the highest-order derivatives of strong solutions  
can not  be achieved by the standard energy method. To obtain the desired decay rates, we shall introduce the method with regard to frequency decomposition. Combining the refined inequality \eqref{NJK4.16}   with the estimate of $\nabla^2(f^L,\rho^L,u^L,\theta^L)$, we  
 obtain the estimate of the highest-order derivative with $(f,\rho,u,\theta)$. 
 Furthermore, we  shall use the 
interpolation techniques to archive the decay rates of solutions in general $L^p$ norms.

In the periodic case, we can obtain  the global existence of solutions by applying the arguments developed in the whole space case  with slight modifications.
As for the decay, thanks to the conservations of mass, momentum  and 
energy, and the Poincar\'{e}'s inequality, we obtain the desired exponential decay rates of solutions 
by constructing suitable functional inequalities. 

 The rest of this article is arranged as follows.
In Section 2, we establish the global existence 
of strong solutions in the whole space $\mathbb{R}^3$, and
prove the existence part of Theorem \ref{T1.1}.
In Section 3, we obtain the time-decay rates of
$(\rho, u, \theta,f)$ and their gradients.
 In Section 4, we establish time-decay rates of the highest-order derivatives of $(f,\rho,u,\theta)$,
 and the decay in $L^p$ norms of the solutions, 
  and then complete  the proofs of Theorem \ref{T1.1}  and Corollary \ref{Cor1}. 
Finally, in Section 5, we present the results on global existence and decay of 
  solutions in the periodic case.

\section{Global existence of strong solutions to the problem \eqref{G1}--\eqref{G5} in $\mathbb{R}^3$}

In this section, we shall establish the global existence of strong solutions to the problem 
\eqref{G1}--\eqref{G5} in the whole space $\mathbb{R}^3$.  We mainly employ the  classical energy method (see, e.g.  \cite{MN-jmku-1980, LMW-siam-2017, MW-cvpde-2020}) to archive it.  
We first derive \emph{a priori} estimates of solutions to \eqref{G1}--\eqref{G5}, and then   construct  local solutions. Finally, we apply  the continuity argument  to extend the local solution    to be global. 
The whole process is similar to that in \cite{MW-cvpde-2020} in some sense. However, due to the lower regularity assumptions
on the initial data, the solutions also have lower regularity. Thus, 
we need to control some terms in the  \eqref{G1}--\eqref{G4}  very carefully. 
Below we shall focus on the differences between our arguments and those in \cite{MW-cvpde-2020}.

\smallskip
We first recall some  basic results which will be used  frequently.

\begin{lem}[see\cite{CDM-krm-2011}] \label{L2.1}
There exists a constant $C>0$ such that for any $f, g \in H^2\left(\mathbb{R}^3\right)$ and any multi-index $\gamma$ with $1 \leq|\gamma| \leq 2$, we have 
\begin{align}
\|f\|_{L^{\infty}\left(\mathbb{R}^3\right)} & \leq C\left\|\nabla_x f\right\|_{L^2\left(\mathbb{R}^3\right)}^{1 / 2}\left\|\nabla_x^2 f\right\|_{L^2\left(\mathbb{R}^3\right)}^{1 / 2}, \\
\|f g\|_{H^1\left(\mathbb{R}^3\right)} & \leq C\|f\|_{H^1\left(\mathbb{R}^3\right)}\left\|\nabla_x g\right\|_{H^1\left(\mathbb{R}^3\right)}, \\
\left\|\partial_x^\gamma(f g)\right\|_{L^2\left(\mathbb{R}^3\right)} & \leq C\left\|\nabla_x f\right\|_{H^1\left(\mathbb{R}^3\right)}\left\|\nabla_x g\right\|_{H^1\left(\mathbb{R}^3\right)} .
\end{align}
\end{lem} 

\smallskip 
\subsection{A priori estimates}

In this subsection, we shall establish the a priori estimates 
of strong solutions to the problem \eqref{G1}--\eqref{G5} in the whole space $\mathbb{R}^3$. We assume that $(f, \rho, u, \theta)$ is a strong solution to   \eqref{G1}--\eqref{G5} on $0\leq t<T$ for some $T>0$ and satisfies
\begin{equation}\label{G2.1}
\sup _{0<t<T}\big\{\|f\|_{H_{x, v}^2}+\|(\rho, u, \theta)(t)\|_{H_x^2}\big\} \leq \delta,
\end{equation}
where $0<\delta<1$ is a sufficiently small constant.
According to the energy norms of the solutions, we divide the estimates into two parts: 
energy estimates of $(f,\rho, u, \theta)$ for the  $x$-variable, and energy estimates of $f$ for mixed space-velocity derivatives.

\subsubsection{Energy estimates in the $x$-variable} 
In this subsection, we show the energy estimates 
with regard to the  spatial  derivatives of the solution $(f, \rho, u, \theta)$.
First, we have
\begin{lem}
 For strong solutions to the problem \eqref{G1}--\eqref{G5}, there exists a positive constant $\lambda_1>0$, such that
\begin{align}\label{G2.5}
& \frac{1}{2} \frac{\mathrm{d}}{\mathrm{d} t}\|(f, \rho, u, \theta)\|^2+\lambda_1\left(\|\{\mathbf{I}-\mathbf{P}\} f\|_v^2+\|b-u\|^2+\|\sqrt{2} \omega-\sqrt{3} \theta\|^2+\|(\nabla_x u, \nabla_x \theta)\|^2\right) \nonumber\\
&\qquad \leq  C\left(\|(\rho, u, \theta)\|_{H^2}+\|(\rho, u, \theta)\|_{H^2}^2\right) \nonumber\\
& \qquad \quad  \times\left(\left\|\nabla_x(a, b, \omega, \rho, u, \theta)\right\|^2+\|u-b\|
^2+ \|\sqrt{2} \omega-\sqrt{3} \theta\|^2+\|\{\mathbf{I}-\mathbf{P}\} f\|_\nu^2\right),
\end{align}
holds for all $0 \leq t<T$. 
\end{lem}
\begin{proof}
Multiplying \eqref{G1}--\eqref{G4} by $f, \rho, u$, and $\theta$ respectively, taking integration and summation, we achieve
\begin{align}\label{G2.6}
&\frac{1}{2} \frac{\mathrm{d}}{\mathrm{d} t}\|(f, \rho, u, \theta)\|^2+\int_{\mathbb{R}^3} \langle-\mathbf{L}\{\mathbf{I}-\mathbf{P}\} f, f\rangle \mathrm{d} x \nonumber\\
&+\mu_1\int_{\mathbb{R}^3} \frac{|\nabla_x u|^2}{1+\rho} \mathrm{d} x+\kappa\int_{\mathbb{R}^3} \frac{|\nabla_x \theta|^2}{1+\rho} \mathrm{d} x+\|b-u\|^2+\|\sqrt{2} \omega-\sqrt{3} \theta\|^2 \nonumber\\
=\,&\int_{\mathbb{R}^3} \Big\langle\frac{1}{2} u\cdot v f, f\Big\rangle \mathrm{d} x-\int_{\mathbb{R}^3}\frac{a|u|^2}{1+\rho}\mathrm{d}x-\int_{\mathbb{R}^3}\frac{\theta u \cdot b}{1+\rho}\mathrm{d}x \nonumber\\
&-\int_{\mathbb{R}^3} \theta\left(\nabla_v f-\frac{v}{2} f\right)\left(\nabla_v f+\frac{v}{2} f\right) \mathrm{d} x \mathrm{d} v-\int_{\mathbb{R}^3} \frac{3a|\theta|^2}{1+\rho}\mathrm{d}x \nonumber\\
&-\int_{\mathbb{R}^3}(u \cdot \nabla_x u) \cdot u \mathrm{d} x-\int_{\mathbb{R}^3} \rho \nabla_x \rho \cdot u \mathrm{d} x-\int_{\mathbb{R}^3} \frac{\theta-\rho}{1+\rho} \nabla_x \rho \cdot u \mathrm{d} x \nonumber\\
&+\int_{\mathbb{R}^3} \frac{\nabla_x \rho \cdot \nabla_x u}{(1+\rho)^2} \cdot u \mathrm{d} x+\int_{\mathbb{R}^3} \frac{\theta\nabla_x \rho}{(1+\rho)^2} \cdot \nabla_x \theta  \mathrm{d} x-\int_{\mathbb{R}^3}\rho^2 \operatorname{div} u \mathrm{d} x \nonumber\\
&+\int_{\mathbb{R}^3}\frac{\theta[(u-b) \cdot u+a|u|^2]}{1+\rho}\mathrm{d}x-\int_{\mathbb{R}^3} \frac{\rho}{1+\rho}(u-b)\cdot u\mathrm{d}x-\int_{\mathbb{R}^3}\frac{\rho\theta}{1+\rho}(\sqrt{6} \omega-3 \theta)\mathrm{d}x\nonumber\\
&+\int_{\mathbb{R}^3}\frac{2\mu_1\theta |D(u)|^2}{1+\rho}+\int_{\mathbb{R}^3} \frac{(\mu_1+\mu_2)\nabla_x({\rm div}_x  u)\cdot u}{1+\rho} \mathrm{d} x+\int_{\mathbb{R}^3}\frac{\theta |{\rm div}_x  u|^2}{1+\rho}\nonumber\\
=:\,&\sum_{j=1}^{17}  {\mathfrak{I}}_j.
\end{align}
 The terms $ {\mathfrak{I}}_1,\dots, {\mathfrak{I}}_{14}$ are as same as those in \cite{MW-cvpde-2020}.
Since no more than two spatial derivatives are needed in the arguments, the estimates are also as same as those in 
\cite{MW-cvpde-2020},  we state the results and omit the details for brevity (see pages 11--13 in  \cite{MW-cvpde-2020}).
 We have 
 \begin{align*}
  \mathfrak{I}_1+\mathfrak{I}_2+\mathfrak{I}_3 
\leq\,&  C\left(\|u\|_{H^1}+\|u\|_{H^2}+\|\rho\|_{H^1}\left(\|u\|_{H^1}+\|\theta\|_{H^1}\right)\right) \\
& \times\big(\|(\mathbf{I}-\mathbf{P}) f\|_v^2+\|u-b\|_{L^2}^2+\|\sqrt{2} \omega-\sqrt{3} \theta\|_{L^2}^2+\|\nabla_x(a, b, \omega, u)\|_{L^2}^2\big),\\
\mathfrak{I}_4+\mathfrak{I}_5 
\leq\, & C\big(\|\theta\|_{H^2}+\|\theta\|_{H^1}\|\rho\|_{H^1}\big) \big(\|\{\mathbf{I}-\mathbf{P}\} f\|_\nu^2+\left\|\nabla_x(a, b, \omega)\right\|_{L^2}^2+\|\sqrt{2} \omega-\sqrt{3} \theta\|_{L^2}^2\big),\\
\mathfrak{I}_6+\mathfrak{I}_7\leq\,&  C\|u\|_{H^1}(\|\nabla_x u\|_{L^2}+\|\nabla_x \rho\|_{L^2}),\\
\mathfrak{I}_8\leq\,& C\|u\|_{H^1}\|\nabla_x \rho\|_{L^2}+ C\|u\|_{H^1}\|\nabla_x \theta\|_{L^2},\\
\mathfrak{I}_9+\mathfrak{I}_{10}\leq\, &C\|(u,\theta)\|_{H^2}\left(\|\nabla_x \rho\|_{L^2}^2+\|\nabla_x u\|_{L^2}^2+\|\nabla_x a\|_{L^2}^2\right),\\
\mathfrak{I}_{11}\leq\,& C\|\rho\|_{H^1}\left(\|\nabla_x \rho\|_{L^2}^2+\|\nabla_x u\|_{L^2}^2\right),\\
\mathfrak{I}_{12}\leq\,& C\|\theta\|_{H^1}\left(\|u-b\|_{L^2}^2+\|\nabla_x u\|_{L^2}^2\right)+ C\|\theta\|_{H^1}\left(\|\nabla_x a\|_{L^2}^2+\|\nabla_x u\|_{L^2}^2\right),\\
\mathfrak{I}_{13}\leq\,& C\|u\|_{H^1}\left(\|u-b\|_{L^2}^2+\|\nabla_x \rho\|_{L^2}^2\right),\\
\mathfrak{I}_{14}\leq\,& C\|\rho\|_{H^1}\left(\|\sqrt{2}\omega-\sqrt{3}\theta\|_{L^2}^2+\|\nabla_x \theta\|_{L^2}^2\right). 
\end{align*}
Using the H\"{o}lder's, Sobolev's inequalities and Lemma \ref{L2.1}, we easily get
\begin{align*}
\mathfrak{I}_{15}\leq\,& C\|\theta\|_{H^2}\|\nabla_x u\|^2,\\
\mathfrak{I}_{16}\leq\,& C\|\rho\|_{H^2}\|\nabla_x u\|^2,\\
\mathfrak{I}_{17}\leq\,& C\|\theta\|_{H^2}\|\nabla_x u\|^2.\\
\end{align*}
Applying all the above estimates to \eqref{G2.6} and noticing the assumption \eqref{G2.1}, we achieve \eqref{G2.5}.
\end{proof}
\begin{lem}\label{L2.3}
For strong solutions to the problem \eqref{G1}--\eqref{G5},  there exists a positive constant $\lambda_2>0$, such that
\begin{align}\label{G2.7}
&\frac{1}{2} \frac{\mathrm{d}}{\mathrm{d} t}  \sum_{1 \leq|\alpha| \leq 2}\left\|\partial^\alpha (f, \rho, u,\theta )\right\|^2\nonumber \\
&\quad \quad + \lambda_2 \sum_{1 \leq|\alpha| \leq 2}\Big(\|\{\mathbf{I}-\mathbf{P}\} \partial^\alpha f\|_\nu^2+\left\|\partial^\alpha(b-u)\right\|^2+\|\partial^\alpha(\sqrt{2} \omega-\sqrt{3} \theta)\|^2+\left\|\nabla_x \partial^\alpha(u, \theta)\right\|^2\Big)\nonumber \\
 &\quad \leq C\Big(\|(\rho, u, \theta)\|_{H^2}+\|(\rho, u, \theta)\|_{H^2}^2+\|(\rho, u, \theta)\|_{H^2}^4\Big)\times\Big(\|\nabla_x(a, b, \omega, \rho)\|_{H^1}^2+\|\nabla_x(u,\theta)\|_{H^2}^2 \nonumber\\
& \quad \quad  +\sum_{1 \leq|\alpha| \leq 2}\big(\|\partial^{\alpha}\{\mathbf{I}-\mathbf{P}\}  f\|_\nu^2+\|\partial^{\alpha} (b-u)\|^2+\|\partial^\alpha(\sqrt{2} \omega-\sqrt{3} \theta)\|^2\big)\Big), 
\end{align}
hold for all $0 \leq t<T$. 
\end{lem}
\begin{proof}
Applying  $\partial^\alpha(1 \leq|\alpha| \leq 2)$ to the system \eqref{G1}-\eqref{G4}, one has
\begin{align}\label{G2.8}
& \partial_t\left(\partial^\alpha f\right)+v \cdot \nabla_x\left(\partial^\alpha f\right)+u \cdot \nabla_v\left(\partial^\alpha f\right)-\partial^\alpha u \cdot v \sqrt{M}-\left(|v|^2-3\right) \sqrt{M} \partial^\alpha \theta-\mathbf{L} \partial^\alpha f\nonumber \\
 &\quad =\frac{1}{2} \partial^\alpha(u \cdot v f)+\left[-\partial^\alpha, u \cdot \nabla_v\right] f+\partial^\alpha\left(\frac{\theta}{\sqrt{M}} \cdot  \Delta_v\big(\sqrt{M} f\big)\right), \\
& \partial_t\left(\partial^\alpha \rho\right)+u \cdot \nabla_x \partial^\alpha \rho+(1+\rho) \operatorname{div} \partial^\alpha u=\left[-\partial^\alpha, \rho {\rm div}_x \right] u+\left[-\partial^\alpha, u \cdot {\rm div} _x\right] \rho ,\\
& \partial_t\left(\partial^\alpha u\right)+u \cdot \nabla_x\left(\partial^\alpha u\right)+\nabla_x \partial^\alpha \theta+\nabla_x \partial^\alpha \rho-\partial^\alpha\left(\frac{\mu_1}{1+\rho} \Delta_x u\right)-\partial^\alpha(b-u) \nonumber \\
&\quad =\left[-\partial^\alpha, u \cdot \nabla_x\right] u+\partial^\alpha\left(\frac{\rho-\theta}{1+\rho} \nabla_x \rho\right)-\partial^\alpha\left(\frac{1}{1+\rho} u a\right)+\partial^{\alpha}\left(
\frac{(\mu_1+\mu_2)\nabla_x({\rm div}_x  u)}{1+\rho}        \right), \\
\label{G2.11}
& \partial_t\left(\partial^\alpha \theta\right)+\operatorname{div}\left(\partial^\alpha u\right)-\partial^\alpha\left(\frac{\kappa}{1+\rho} \Delta_x \theta\right)-\sqrt{3}\left(\sqrt{2} \partial^\alpha \omega-\sqrt{3} \partial^\alpha \theta\right) \nonumber \\
&\quad =-\partial^\alpha(u \cdot \nabla_x \theta+\theta \operatorname{div} u)+\partial^\alpha\left(\frac{1}{1+\rho}\left((u-b) \cdot u+a|u|^2-u \cdot b\right)\right)+\partial^{\alpha}\left(\frac{|{\rm div}_x  u|^2}{1+\rho}         \right)\nonumber \\
&\qquad -\sqrt{3}\left(\frac{\rho}{1+\rho}\big(\sqrt{2} \partial^\alpha \omega-\sqrt{3} \partial^\alpha \theta\big)\right)-\partial^\alpha\left(\frac{3}{1+\rho} a \cdot \theta\right)+2\mu_1\partial^\alpha\left( \frac{|D(u)|^2}{1+\rho}    \right),
\end{align}
where $[X, Y]$ := $X Y-Y X$ for two operators $X$ and $Y$.

Multiplying \eqref{G2.8}--\eqref{G2.11} by $\partial^\alpha f,
{\partial^\alpha \rho}, \partial^\alpha u,$ and $ \partial^\alpha \theta$ respectively, and taking integration and summation, we have
\begin{align}\label{G2.12}
&\frac{1}{2} \frac{\mathrm{d}}{\mathrm{d} t} \left\|\partial^\alpha (f, \rho, u, \theta)\right\|^2+\mu_1\int_{\mathbb{R}^3}\frac{|\nabla_x\left(\partial^\alpha u\right)|^2}{1+\rho}\mathrm{d} x+\kappa\int_{\mathbb{R}^3}\frac{\left|\nabla_x\left(\partial^\alpha \theta\right)\right|^2 }{1+\rho}\mathrm{d} x\nonumber\\
& +\int_{\mathbb{R}^3}\left\langle-\mathbf{L}\{\mathbf{I}-\mathbf{P}\} \partial^\alpha f,\{\mathbf{I}-\mathbf{P}\} \partial^\alpha f\right\rangle \mathrm{d} x+\left\|\partial^\alpha(b-u)\right\|^2+\|\partial^\alpha(\sqrt{2} \omega-\sqrt{3} \theta)\|^2 \nonumber\\
=\, & \int_{\mathbb{R}^3} \frac{1}{2}\left\langle\partial^\alpha(u \cdot v f), \partial^\alpha f\right\rangle \mathrm{d} x+\int_{\mathbb{R}^3}\left\langle\partial^\alpha\left(\frac{\theta }{\sqrt{M}} \Delta_v(\sqrt{M}f)\right), \partial^\alpha f\right\rangle \mathrm{d} x. \nonumber\\
& +\int_{\mathbb{R}^3}\left\langle\left[-\partial^\alpha, u \cdot \nabla_v\right] f, \partial^\alpha f\right\rangle \mathrm{d} x+\int_{\mathbb{R}^3}\left[-\partial^\alpha, u \cdot \nabla_x\right] \rho \partial^\alpha \rho \mathrm{d} x \nonumber\\
& +\int_{\mathbb{R}^3}\left[-\partial^\alpha, u \cdot \nabla_x\right] u \partial^\alpha u \mathrm{d} x+\int_{\mathbb{R}^3}\left[-\partial^\alpha, \rho {\rm div}_x  \right] u \partial^\alpha \rho \mathrm{d} x \nonumber\\
& -\int_{\mathbb{R}^3} \partial^\alpha\left(\frac{a u}{1+\rho}\right) \cdot \partial^\alpha u \mathrm{d} x-\int_{\mathbb{R}^3} \partial^\alpha\left(\frac{u \cdot b}{1+\rho}\right) \partial^\alpha \theta \mathrm{d} x-\int_{\mathbb{R}^3} \partial^\alpha\left(\frac{3 a \theta}{1+\rho}\right)  \partial^\alpha \theta \mathrm{d} x \nonumber \\
& +\int_{\mathbb{R}^3} \partial^\alpha\left(\frac{\rho-\theta}{1+\rho} \nabla_x \rho\right)\cdot \partial^\alpha u \mathrm{d} x-\int_{\mathbb{R}^3} \rho \operatorname{div} \partial^\alpha u \partial^\alpha \rho \mathrm{d} x-\int_{\mathbb{R}^3} \partial^\alpha(\theta \operatorname{div} u) \partial^\alpha \theta \mathrm{d} x \nonumber\\
& -\int_{\mathbb{R}^3} u \cdot \nabla_x \partial^\alpha \rho \partial^\alpha \rho \mathrm{d} x-\int_{\mathbb{R}^3} u \cdot \nabla_x \partial^\alpha u \cdot \partial^\alpha u \mathrm{d} x-\int_{\mathbb{R}^3} \partial^\alpha(u \cdot \nabla_x \theta) \partial^\alpha \theta \mathrm{d} x \nonumber \\
& -\sum_{1 \leq |\beta| \leq |\alpha|} C_{\alpha, \beta} \int_{\mathbb{R}^3} \partial^\beta\left(\frac{1}{1+\rho}\right) \partial^{\alpha-\beta} \Delta_x u \partial^\alpha u \mathrm{d} x- \int_{\mathbb{R}^3} \nabla_x \frac{1}{1+\rho} \nabla_x \partial^\alpha u\partial^\alpha u \mathrm{d} x \nonumber \\
& -\sum_{1 \leq |\beta| \leq |\alpha|} C_{\alpha, \beta} \int_{\mathbb{R}^3} \partial^\beta\left(\frac{1}{1+\rho}\right) \partial^{\alpha-\beta} \Delta_x \theta \partial^\alpha \theta \mathrm{d} x-\int_{\mathbb{R}^3} \nabla_x \frac{1}{1+\rho} \nabla_x\partial^\alpha \theta \partial^\alpha \theta \mathrm{d} x  \nonumber\\
& -\int_{\mathbb{R}^3} \partial^\alpha\left(\frac{\rho}{1+\rho}(b-u)\right) \cdot \partial^\alpha u \mathrm{d} x-\sqrt{3} \int_{\mathbb{R}^3} \partial^\alpha\left(\frac{\rho}{1+\rho}(\sqrt{2} \omega-\sqrt{3} \theta)\right)  \partial^\alpha \theta \mathrm{d} x \nonumber \\
& +\int_{\mathbb{R}^3} \partial^\alpha\left(\frac{1}{1+\rho}(u-b) \cdot u\right) \partial^\alpha \theta \mathrm{d} x+\int_{\mathbb{R}^3} \partial^\alpha\left(\frac{1}{1+\rho} a|u|^2\right) \partial^\alpha \theta \mathrm{d} x \nonumber \\
&+2\mu_1\int_{\mathbb{R}^3}\partial^\alpha\left( \frac{|D(u)|^2}{1+\rho}\right)\partial^\alpha \theta\mathrm{d}x+\mu_2\int_{\mathbb{R}^3}\partial^{\alpha}\left(\frac{|{\rm div}_x  u|^2}{1+\rho}\right)\partial^{\alpha}\theta\mathrm{d} x\nonumber\\
&+\int_{\mathbb{R}^3}\partial^{\alpha}\left(\frac{(\mu_1+\mu_2)\nabla_x({\rm div}_x  u)}{1+\rho}        \right)\partial^{\alpha} u\mathrm{d}x\nonumber\\
=:\,&\sum_{i=1}^{26} \mathfrak{J}_i,
\end{align}
where $C_{\alpha, \beta}$ is a positive constant depending only on $\alpha$ and $\beta$.
First, we calculate $\mathfrak{J}_1$ and $\mathfrak{J}_2$ to get 
\begin{align*}
\mathfrak{J}_1 & \leq C\left\|\partial^\alpha(u f)\right\|\left\|v \partial^\alpha f\right\| \\
& \leq C\|\nabla_x u\|_{H^1}\left\|\nabla_x f\right\|_{L_v^2\left(H_x^1\right)}\left\|v \partial^\alpha f\right\| \\
&\leq C\|u\|_{H^2}\bigg(\sum_{1\leq |\alpha|\leq 2} \|\partial^{\alpha}\{\mathbf{I}-\mathbf{P}\}f\|_\nu^2+\|\nabla_x(a,b,\omega)\|_{H^1}^2    \bigg),\\
\mathfrak{J}_2 & =-\int_{\mathbb{R}^3}\left\langle \partial^\alpha\left(\theta\left(\nabla_v f-\frac{v}{2} f\right)\right) ,\left(\nabla_v \partial^\alpha f+\frac{v}{2} \partial^\alpha f\right)\right\rangle\mathrm{d}x \\
& \leq\left\|\partial^\alpha\left(\theta\left(\nabla_v f-\frac{v}{2} f\right)\right)\right\| \left\| \nabla_v \partial^\alpha f+\frac{v}{2} \partial^\alpha f \right\| \\
& \leq C\|\nabla_x \theta\|_{H^1}\left\|\left(\nabla_v-\frac{v}{2}\right) \nabla_x f\right\|_{L_v^2\left(H_x^1\right)}\left\|\left(\nabla_v+\frac{v}{2}\right) \nabla_x f\right\|_{L_v^2\left(H_x^1\right)}\\
&\leq C\|\nabla_x \theta\|_{H^1}\bigg(\sum_{1 \leq\left|\alpha\right| \leq 2}\left\|\partial^{\alpha}\{\mathbf{I}-\mathbf{P}\}  f\right\|_\nu^2+\|\nabla_x(a, b, \omega)\|_{H^1}^2\bigg).
\end{align*}
Using the H\"{o}lder's, Sobolev's and Young's inequalities, we get the following bounds
\begin{align*}
\mathfrak{J}_3=\,&\int_{\mathbb{R}^3}\left\langle\partial^\alpha(u f), \partial_v \partial^\alpha f\right\rangle \mathrm{d} x \\
\leq\,& \left\|\partial^\alpha(u f)\right\|\left\|\nabla_v \partial^\alpha f\right\| \\
\leq \,&C\|\nabla_x u\|_{H^1}\left\|\nabla_x f\right\|_{L_v^2\left(H_x^1\right)}\left\|\nabla_v \partial^\alpha f\right\|\\
\leq\,& C\|u\|_{H^2}\bigg( \sum_{1 \leq\left|\alpha\right| \leq 2}\left\|\partial^{\alpha}\{\mathbf{I}-\mathbf{P}\}  f\right\|_\nu^2+\|\nabla_x(a, b, \omega)\|_{H^1}^2   \bigg), \\
\mathfrak{J}_4 \leq\,& C\|\rho\|_{H^2}\left(\|\nabla_x \rho\|_{H^1}^2+\left\|\nabla_x u\right\|_{H^2}^2\right), \\
\mathfrak{J}_5 \leq\,& C\|u\|_{H^2}\|\nabla_x u\|_{H^2}^2, \\
\mathfrak{J}_6 \leq \,&C\|\rho\|_{H^2}\left(\|\nabla_x \rho\|_{H^1}^2+\left\|\nabla_x u\right\|^2\right), \\
\mathfrak{J}_7 \leq\,&\left\|\partial^\alpha\left(\frac{a u}{1+\rho}\right)\right\|\left\|\partial^\alpha u\right\| \\
\leq\, & C\left\|\nabla_x\left(\frac{u}{1+\rho}\right)\right\|_{H^1}\|\nabla_x a\|_{H^1}\left\|\partial^\alpha u\right\| \\
\leq\,& C\|u\|_{H^2}\left(\|\nabla_x a\|_{H^1}^2+\left\|\nabla_x u\right\|_{H^2}^2\right)+ C\|u\|_{H^2}\left(\|\nabla_x \rho\|_{H^1}^2+\left\|\nabla_x a\right\|_{H^1}^2\right), \\
\mathfrak{J}_8 \leq\,&\left\|\partial^\alpha\left(\frac{u \cdot b}{1+\rho}\right)\right\|\left\|\partial^\alpha \theta\right\| \\
\leq\,& C\left\|\nabla_x\left(\frac{u}{1+\rho}\right)\right\|_{H^1}\|\nabla_x b\|_{H^1}\left\|\partial^\alpha \theta\right\| \\
\leq\,& C\|\theta\|_{H^2}\left(\|\nabla_x b\|_{H^1}^2+\left\|\nabla_x u\right\|_{H^2}^2\right)+ C\|\rho\|_{H^2}\|\theta\|_{H^2}\left(\|\nabla_x b\|_{H^1}^2+\left\|\nabla_x u\right\|_{H^2}^2\right), \\
\mathfrak{J}_9 \leq\,&\left\|\partial^\alpha\left(\frac{a \theta}{1+\rho}\right)\right\|\left\|\partial^\alpha \theta\right\| \\
\leq\,& C\left\|\nabla_x\left(\frac{\theta}{1+\rho}\right)\right\|_{H^1}\|\nabla_x a\|_{H^1}\left\|\partial^\alpha \theta\right\| \\
\leq\,& C\|\theta\|_{H^2}\left(\|\nabla_x a\|_{H^1}^2+\left\|\nabla_x \theta\right\|_{H^2}^2\right)+ C\|\rho\|_{H^2}\|\theta\|_{H^2}\left(\|\nabla_x a\|_{H^1}^2+\left\|\nabla_x \theta\right\|_{H^2}^2\right), \\
\mathfrak{J}_{10}=&\sum_{1 \leq |\beta| \leq |\alpha|} C_{\alpha, \beta} \int_{\mathbb{R}^3} \partial^\beta\left(\frac{\rho-\theta}{1+\rho} \partial^{\alpha-\beta} \nabla_x \rho\right)\cdot \partial^\alpha u \mathrm{d} x+\int_{\mathbb{R}^3} \frac{\rho-\theta}{1+\rho} \nabla_x \partial^\alpha \rho \partial^\alpha u \mathrm{d} x \\
=\,&\sum_{1 \leq |\beta| \leq |\alpha|} C_{\alpha, \beta} \int_{\mathbb{R}^3} \partial^\beta\left(\frac{\rho-\theta}{1+\rho} \partial^{\alpha-\beta} \nabla_x \rho\right)\cdot \partial^\alpha u \mathrm{d} x \\
&-\int_{\mathbb{R}^3} \frac{\rho-\theta}{1+\rho} \partial^\alpha \rho \partial^\alpha \operatorname{div} u \mathrm{d} x-\int_{\mathbb{R}^3} \nabla_x \frac{\rho-\theta}{1+\rho} \partial^\alpha \rho\cdot \partial^\alpha u \mathrm{d} x \\
\leq\,& C\left(\|(\rho,\theta)\|_{H^2}+\| (\rho,\theta)\|_{H^2}^2\right)\left(\|\nabla_x \rho\|_{H^1}^2+\left\|\nabla_x u\right\|_{H^2}^2\right). 
\end{align*}
Similarly, by direct calculations, one has
\begin{align*}
\mathfrak{J}_{11}  \leq\,& C\|\rho\|_{H^2}\left(\|\nabla_x \rho\|_{H^1}^2+\left\|\nabla_x u\right\|_{H^2}^2\right), \\
\mathfrak{J}_{12}  \leq\,& C\|\theta\|_{H^2}\left(\|\nabla_x \theta\|_{H^2}^2+\left\|\nabla_x u\right\|_{H^2}^2\right), \\
\mathfrak{J}_{13} \leq\,& C\|\rho\|_{H^2}\left(\|\nabla_x u\|_{H^2}^2+\left\|\nabla_x \rho\right\|_{H^1}^2\right), \\
\mathfrak{J}_{14} \leq \,& C\|u\|_{H^2}\left\|\nabla_x u\right\|_{H^2}^2, \\
\mathfrak{J}_{15} \leq \,&C\|\theta\|_{H^2}\left(\|\nabla_x u\|_{H^2}^2+\left\|\nabla_x \theta\right\|_{H^2}^2\right), \\
\mathfrak{J}_{16} + \mathfrak{J}_{17}\leq \,& C\|\rho\|_{H^2}\left\|\nabla_x u\right\|_{H^2}^2, \\
\mathfrak{J}_{18} +\mathfrak{J}_{19} \leq \,& C\|\rho\|_{H^2}\left\|\nabla_x \theta\right\|_{H^2}^2.
\end{align*}
For the remaining terms, we have
\begin{align*}
\mathfrak{J}_{20} \leq\,&\left\|\partial^\alpha\left( \frac{\rho(b-u)}{1+\rho}\right)\right\|\left\|\partial^\alpha u\right\|\\
 \leq\,& \left\|\nabla_x \frac{\rho}{1+\rho}\right\|_{H^1}\|\nabla_x(b-u)\|_{H^1}\left\|\partial^\alpha u\right\| \\
\leq\,& C\|\rho\|_{H^2}\bigg(\sum_{1\leq |\alpha|\leq 2}\|\partial^{\alpha} (b-u)\|^2+\|\nabla_x u\|_{H^2}^2\bigg), \\
\mathfrak{J}_{21} \leq \,&C\|\rho\|_{H^2}\bigg(\sum_{1\leq |\alpha|\leq 2}\|\partial^{\alpha} (\sqrt{2}\omega-\sqrt{3}\theta)\|^2+\|\nabla_x u\|_{H^2}^2\bigg), \\
\mathfrak{J}_{22}  \leq \,&C\|u\|_{H^2}\bigg(\sum_{1\leq |\alpha|\leq 2}\|\partial^{\alpha} (b-u)\|^2+\|\nabla_x \theta\|_{H^2}^2\bigg) \\
&+C\left(\|u\|_{H^2}^2+\|\rho\|_{H^2}^2\right)\bigg(\sum_{1\leq |\alpha|\leq 2}\|\partial^{\alpha} (b-u)\|^2+\|\nabla_x \theta\|_{H^2}^2\bigg),\\
\mathfrak{J}_{23} \leq\,& C\left\|\nabla_x \frac{|u|^2}{1+\rho}\right\|_{H^1}\|\nabla_x a\|_{H^1}\left\|\partial^\alpha \theta\right\| \\
\leq\,& C\|u\|_{H^2}^2\left(\|\nabla_x a\|_{H^1}^2+\|\nabla_x\theta\|_{H^2}^2\right)+ C\left(\|\rho\|_{H^2}^2+\|u\|_{H^2}^4\right)\left(\|\nabla_x a\|_{H^1}^2+\|\nabla_x\theta\|_{H^2}^2\right), \\
\mathfrak{J}_{24}\leq\, &C\left(\|\rho\|_{H^2}^2+\|\theta\|_{H^2}^2 \right)\|\nabla_x u\|_{H^2}^2,\\
\mathfrak{J}_{25}\leq \,&C\left(\|\rho\|_{H^2}^2+\|\theta\|_{H^2}^2 \right)\|\nabla_x u\|_{H^2}^2,\\
\mathfrak{J}_{26} \leq\,&  C\|\rho\|_{H^2}\left\|\nabla_x u\right\|_{H^2}^2. 
\end{align*}
Applying all the above estimates on $\mathfrak{J}_i(1 \leq i \leq 26)$ to \eqref{G2.12},   taking summation over $1 \leq|\alpha| \leq 2$, and noticing the assumption \eqref{G2.1}, we easily achieve \eqref{G2.7}.
\end{proof}

In order to get the estimates of the energy dissipation of $\left\|\nabla_x(a, b, \omega)\right\|_{H^1}$, we shall study the equations of $a, b,$ and $ \omega$ involving
\begin{align*}
\Gamma_{i, j} g=\big\langle\left(v_i v_j-1\right) \sqrt{M(v)}, g\big\rangle, \quad 
\Upsilon_i g=\left\langle\frac{1}{\sqrt{6}} v_i\left(|v|^2-3\right) \sqrt{M(v)}, g\right\rangle,
\end{align*}
for any $g=g(v), 1 \leq i, j \leq 3$. By direct calculations, we get the following  system of   $(a, b, \omega)$: 
\begin{align}\label{G2.13}
&\partial_t a+{\rm div}_x  b=0,\\
\label{G2.14}
& \partial_t b_i+\partial_i a+\frac{2}{\sqrt{6}} \partial_i \omega+\sum_{j=1}^3 \partial_{x_j} \Gamma_{i, j}(\{\mathbf{I}-\mathbf{P}\} f)=-b_i+u_i+u_i a, \\
\label{G2.15}
& \partial_t \omega+\sqrt{2}(\sqrt{2} \omega-\sqrt{3} \theta)-\sqrt{6} a \theta+\frac{2}{\sqrt{6}} {\rm div}_x  b-\frac{2}{\sqrt{6}} u \cdot b \nonumber\\
& \qquad  +\sum_{i=1}^3 \partial_{x_i} \Upsilon_i(\{\mathbf{I}-\mathbf{P}\} f)=0, \\
\label{G2.16}
& \partial_j b_i+\partial_i b_j-\left(u_i b_j+u_j b_i\right)\nonumber \\
& \qquad  -\frac{2}{\sqrt{6}} \delta_{i j}\left(\frac{2}{\sqrt{6}} {\rm div}_x  b-\frac{2}{\sqrt{6}} u \cdot b+\sum_{i=1}^3 \partial_{x_i} \Upsilon_i(\{\mathbf{I}-\mathbf{P}\} f)\right)\nonumber \\
&\quad =-\partial_t \Gamma_{i, j}(\{\mathbf{I}-\mathbf{P}\} f)+\Gamma_{i, j}(l+r+s), \\
\label{G2.17}
& \frac{5}{3}\left(\partial_i \omega-\omega u_i-\sqrt{6} \theta b_i\right)-\frac{2}{\sqrt{6}} \sum_{j=1}^3 \partial_{x_j} \Gamma_{i, j}(\{\mathbf{I}-\mathbf{P}\} f)\nonumber \\
&\quad  =-\partial_t \Upsilon_i(\{\mathbf{I}-\mathbf{P}\} f)+\Upsilon_i(\mathfrak{l}+\mathfrak{r}+\mathfrak{s}).
\end{align}
Here, $\mathfrak{l}, \mathfrak{r},$ and $ \mathfrak{s}$ in \eqref{G2.16}--\eqref{G2.17} are defined as
\begin{align}
& \mathfrak{l}=-v \cdot \nabla_x\{\mathbf{I}-\mathbf{P}\} f+\mathbf{L}\{\mathbf{I}-\mathbf{P}\} f,\label{NJK2.181} \\
& \mathfrak{r}=-u \cdot \nabla_v\{\mathbf{I}-\mathbf{P}\} f+\frac{1}{2} u \cdot v\{\mathbf{I}-\mathbf{P}\} f, \label{NJK2.182}\\
& \mathfrak{s}=\frac{\theta}{\sqrt{M}} \operatorname{div}_v\left(\sqrt{M}\left(\nabla_v-\frac{v}{2}\right)\{\mathbf{I}-\mathbf{P}\} f\right) .\label{NJK2.18}
\end{align}
To get the above equations, we need to represent \eqref{G1} as
\begin{align}\label{G2.19}
& \partial_t \mathbf{P} f+v \cdot \nabla_x \mathbf{P} f+u \cdot \nabla_v \mathbf{P} f-\frac{1}{2} u \cdot v \mathbf{P} f-u \cdot v \sqrt{M}-\theta\left(|v|^2-3\right) \sqrt{M} \nonumber\\
&\quad +\mathbf{P}_1 f+2 \mathbf{P}_2 f-\frac{\theta}{\sqrt{M}} \operatorname{div}_v\left(\sqrt{M}\left(\nabla_v-\frac{v}{2}\right) \mathbf{P} f\right) b =-\partial_t\{\mathbf{I}-\mathbf{P}\} f+\mathfrak{l}+\mathfrak{r}+\mathfrak{s }.
\end{align}
Applying $\Gamma_{i j}$ to \eqref{G2.19} and combining \eqref{G2.13} with \eqref{G2.15}, we obtain \eqref{G2.16}. Besides, applying $\Upsilon_i$ to \eqref{G2.19} and combing \eqref{G2.14}, we  get \eqref{G2.17}.

Now we  define the temporal functional
\begin{align}
\mathcal{E}_0(t):= & \sum_{|\alpha| \leq 1} \sum_{i, j} \int_{\mathbb{R}^3} \partial^\alpha\left(\partial_j b_i+\partial_i b_j\right) \partial_x^\alpha \Gamma_{i, j}(\{\mathbf{I}-\mathbf{P}\} f) \mathrm{d} x \nonumber\\
& +\sum_{|\alpha| \leq 1} \sum_i \int_{\mathbb{R}^3} \partial^\alpha \partial_i \omega \partial^\alpha \Upsilon_i(\{\mathbf{I}-\mathbf{P}\} f) \mathrm{d} x \nonumber\\
& +\frac{2}{21} \sum_{|\alpha| \leq 1} \int_{\mathbb{R}^3} \partial_x^\alpha a \partial^\alpha\left(\frac{\sqrt{6}}{5} \sum_i \partial_i \Upsilon_i(\{\mathbf{I}-\mathbf{P}\} f)-{\rm div}_x  b\right) \mathrm{d} x .
\end{align} 
We   have    the following lemma.

\begin{lem}\label{L2.4}
For strong solutions to the problem  \eqref{G1}--\eqref{G5}, there exists a positive constant $\lambda_3>0$, such that
\begin{align}\label{G2.21}
& \frac{\mathrm{d}}{\mathrm{d} t} \mathcal{E}_0(t)+\lambda_3\left\|\nabla_x(a, b, \omega)\right\|_{H^1}^2 \nonumber\\
&\quad\leq C\left(\|\{\mathbf{I}-\mathbf{P}\} f\|_{L_v^2\left(H_x^2\right)}^2+\|u-b\|_{H^1}^2+\|\sqrt{2} \omega-\sqrt{3} \theta\|_{H^1}^2\right) \nonumber\\
& \qquad +C\left(\|u, \theta\|_{H^1}+\|u, \theta\|_{H^1}^2\right)\times\left(\left\|\nabla_x\{\mathbf{I}-\mathbf{P}\} f\right\|_{L_v^2\left(H_x^1\right)}^2+\left\|\nabla_x(a, b, \omega)\right\|_{H^1}^2\right), 
\end{align}
holds for all $0 \leq t<T$.
\end{lem}

\begin{proof}
Through \eqref{G2.16}, it's easy to get 
\begin{align}\label{G2.22}
& \sum_{i, j} \left\|\partial^\alpha\left(\partial_i b_j+\partial_j b_i\right)\right\|^2-\sum_{i, j} \int_{\mathbb{R}^3} \frac{2}{3} \delta_{i j} \partial^\alpha ({\rm div}_x  b)\cdot \partial^\alpha\left(\partial_i b_j+\partial_j b_i\right) \mathrm{d} x \nonumber\\
&\quad = -\frac{\mathrm{d}}{\mathrm{d} t} \sum_{i, j} \int_{\mathbb{R}^3} \partial^\alpha\left(\partial_i b_j+\partial_j b_i\right) \partial^\alpha \Gamma_{i, j}(\{\mathbf{I}-\mathbf{P}\} f) \mathrm{d} x \nonumber\\
&\qquad  +\sum_{i, j} \int_{\mathbb{R}^3} \partial^\alpha\left(\partial_i \partial_t b_j+\partial_j \partial_t b_i\right) \partial^\alpha \Gamma_{i, j}(\{\mathbf{I}-\mathbf{P}\} f) \mathrm{d} x \nonumber\\
& \qquad +\sum_{i, j} \int_{\mathbb{R}^3} \partial^\alpha\left(\partial_i b_j+\partial_j b_i\right) \partial^\alpha\Big[\left(u_i b_j+u_j b_i\right)-\delta_{i j}\sum_i \partial_i \Upsilon_i\{\mathbf{I}-\mathbf{P}\} f\big.\nonumber \\
&\qquad \qquad\qquad\big.-\frac{2}{3}\delta_{i j} (u \cdot b)+\Gamma_{i, j}(\mathfrak{l}+\mathfrak{r}+\mathfrak{s})\Big] \mathrm{d} x .
\end{align}
By \eqref{G2.14}, Lemma \ref{L2.1} and Young's inequality, one has
$$
\begin{aligned}
& \sum_{i, j} \int_{\mathbb{R}^3} \partial^\alpha\left(\partial_i \partial_t b_j+\partial_j \partial_t b_i\right) \partial^\alpha \Gamma_{i, j}(\{\mathbf{I}-\mathbf{P}\} f) \mathrm{d} x \\
\, &\quad \leq\varepsilon_1\left\|\nabla_x(a, \omega)\right\|_{H^1}^2+C_{\varepsilon_1}\left\|\nabla_x\{\mathbf{I}-\mathbf{P}\} f\right\|_{L_v^2\left(H_x^1\right)}^2 \\
&\qquad +C\left(\|u-b\|_{H^1}^2+\|u\|_{H^1}^2\left\|\nabla_x a\right\|_{H^1}^2\right),
\end{aligned}
$$
where $\varepsilon_1>0$ is sufficient small.

By a direct calculation, we have
\begin{align*}
& \sum_{i, j} \int_{\mathbb{R}^3} \partial^\alpha\left(\partial_i b_j+\partial_j b_i\right) \partial^\alpha\Big[\left(u_i b_j+u_j b_i\right)- \delta_{i j}\sum_i \partial_i \Upsilon_i\{\mathbf{I}-\mathbf{P}\} f\big. \\
& \qquad\quad \big.-\frac{2}{3}\delta_{i j}( u \cdot b)+\Gamma_{i, j}(\mathfrak{l}+\mathfrak{r}+\mathfrak{s})\Big] \mathrm{d} x \\
\leq \,& \frac{1}{4} \sum_{i, j}\left\|\partial^\alpha\left(\partial_i b_j+\partial_j b_i\right)\right\|^2+C \sum_{i, j}\left(\left\|\partial^\alpha\left(u_i b_j+u_j b_i\right)\right\|^2+\left\|\partial^\alpha \partial_i \Upsilon_i\{\mathbf{I}-\mathbf{P}\} f\right\|^2\right. \\
& \left.+\left\|\partial^\alpha(u \cdot b)\right\|^2+\left\|\partial^\alpha \Gamma_{i, j}(\mathfrak{l}+\mathfrak{r}+\mathfrak{s})\right\|^2\right) \\
\leq\,& \frac{1}{4} \sum_{i, j}\left\|\partial^\alpha\left(\partial_i b_j+\partial_j b_i\right)\right\|^2+C\|\{\mathbf{I}-\mathbf{P}\} f\|_{L_v^2\left(H_x^2\right)}^2\\
\quad &+C\|(u, \theta)\|_{H^1}^2\Big(\|\nabla_x b\|_{H^1}^2+\left\|\nabla_x\{\mathbf{I}-\mathbf{P}\} f\right\|_{L_v^2\left(H_x^1\right)}^2\Big) .
\end{align*}
Here we have used the following facts:
\begin{align*}
\sum_{i, j}\left\|\partial^\alpha \Gamma_{i, j}(\mathfrak{l})\right\|^2  \leq\,& C\|\{\mathbf{I}-\mathbf{P}\} f\|_{L_v^2\left(H_x^2\right)}^2, \\
\sum_{i, j}\left\|\partial^\alpha \Gamma_{i, j}(\mathfrak{r})\right\|^2  \leq\,& C\|u\|_{H^1}^2\left\|\nabla_x\{\mathbf{I}-\mathbf{P}\} f\right\|_{L_v^2\left(H_x^1\right)}^2, \\
\sum_{i, j}\left\|\partial^\alpha \Gamma_{i, j}(\mathfrak{s})\right\|^2  \leq \,&C\|\theta\|_{H^1}^2\left\|\nabla_x\{\mathbf{I}-\mathbf{P}\} f\right\|_{L_v^2\left(H_x^1\right)}^2,\\
\left\|\partial^\alpha \partial_i \Upsilon_i\{\mathbf{I}-\mathbf{P}\} f\right\|^2 \leq\,& C\left\|\nabla_x\{\mathbf{I}-\mathbf{P}\} f\right\|_{L_v^2\left(H_x^1\right)}^2, \\
\sum_{i, j}\left\|\partial^\alpha\left(u_i b_j+u_j b_i\right)\right\|^2+\left\|\partial^\alpha(u \cdot b)\right\|^2  \leq\,& C\|u\|_{H^1}^2\left\|\nabla_x b\right\|_{H^1}^2.
\end{align*}
It's easy to verify that
\begin{gather*}
\sum_{i, j}\left\|\partial^\alpha\left(\partial_i b_j+\partial_j b_i\right)\right\|^2 =2\left\|\nabla_x \partial^\alpha b\right\|^2+2\left\|\partial^\alpha{\rm div} b\right\|^2, \\
-\sum_{i, j} \frac{2}{3} \delta_{i j} \int_{\mathbb{R}^3} \partial^\alpha {\rm div}_x  b \partial^\alpha\left(\partial_i b_j+\partial_j b_i\right) \mathrm{d} x 
=\frac{4}{3}\left\|\partial^\alpha {\rm div}_x  b\right\|^2.
\end{gather*}
Applying the above estimates to \eqref{G2.22} and taking the summation over $|\alpha| \leq 1$, one has
\begin{align}\label{G2.23}
& \frac{\mathrm{d}}{\mathrm{d} t} \sum_{|\alpha| \leq 1} \sum_{i, j} \int_{\mathbb{R}^3} \partial^\alpha\left(\partial_i b_j+\partial_j b_i\right) \partial^\alpha \Gamma_{i, j}(\{\mathbf{I}-\mathbf{P}\} f) \mathrm{d} x \nonumber\\
&\qquad +\sum_{|\alpha| \leq 1}\bigg(\frac{3}{2}\left\|\nabla_x \partial^\alpha b\right\|^2+\frac{1}{6}\left\|\partial^\alpha \nabla_x \cdot b\right\|^2\bigg) \nonumber\\
\, & \quad \leq\varepsilon_2 \|\nabla_x(a, \omega)\|_{H^1}^2+C_{\varepsilon_2}\|\{\mathbf{I}-\mathbf{P}\} f\|_{L_x^2\left(H_x^2\right)}^2+C\|u-b\|_{H^1}^2 \nonumber\\
&\qquad +C\|(u, \theta)\|_{H^1}^2 \left(\left\|\nabla_x\{\mathbf{I}-\mathbf{P}\} f\right\|_{L_v^2\left(H_x^1\right)}^2+\left\|\nabla_x(a, b)\right\|_{H^1}^2\right),
\end{align}
where $\varepsilon_2>0$ is sufficient small.

By \eqref{G2.17}, we get
\begin{align}\label{G2.24}
\left\|\partial^\alpha \partial_i \omega\right\|^2= & -\frac{3}{5} \frac{\mathrm{d}}{\mathrm{d} t} \int_{\mathbb{R}^3} \partial^\alpha \partial_i \omega \partial^\alpha \Upsilon_i\{\mathbf{I}-\mathbf{P}\} f \mathrm{d} x+\frac{3}{5} \int_{\mathbb{R}^3} \partial^\alpha \partial_i \partial_t \omega \partial^\alpha \Upsilon_i\{\mathbf{I}-\mathbf{P}\} f \mathrm{d} x \nonumber\\
& +\int_{\mathbb{R}^3} \partial^\alpha \partial_i \omega \partial^\alpha\left(\omega u_i+\sqrt{6} \theta b_i\right) \mathrm{d} x+\frac{3}{5}\int_{\mathbb{R}^3} \partial^\alpha \partial_i \omega \partial^\alpha
 \Upsilon_i(\mathfrak{l}+\mathfrak{r}+\mathfrak{s}) \mathrm{d} x \nonumber\\
& +\frac{\sqrt{6}}{5} \int_{\mathbb{R}^3} \partial^\alpha \partial_i \omega \partial^\alpha\bigg(\sum_j \partial_j \Gamma_{i, j}\{\mathbf{I}-\mathbf{P}\} f\bigg) \mathrm{d} x.
\end{align}
Through direct calculations, we arrive at
\begin{align*}
&\int_{\mathbb{R}^3} \partial^\alpha \partial_i \partial_t \omega \partial^\alpha \Upsilon_i\{\mathbf{I}-\mathbf{P}\} f \mathrm{d} x\\
= & \int_{\mathbb{R}^3} \partial_i \partial^\alpha\Big(\sqrt{6} a \theta-\sqrt{2}(\sqrt{2} \omega-\sqrt{3} \theta)-\frac{2}{\sqrt{6}} {\rm div}_x  b+\frac{2}{\sqrt{6}} u \cdot b\bigg. \\
&  \quad\bigg.-\sum_i \partial_i \Upsilon_i\{\mathbf{I}-\mathbf{P}\} f\Big)  \partial^\alpha \Upsilon_i\{\mathbf{I}-\mathbf{P}\} f \mathrm{d} x \\
\leq \,& \varepsilon_3\|\nabla_x b\|_{H^1}^2+C_{\varepsilon_3}\|\{\mathbf{I}-\mathbf{P}\} f\|_{L_v^2\left(H_x^2\right)}^2+C\|\sqrt{2} \omega-\sqrt{3} \theta\|_{H^1}^2 \\
& +C\|(u, \theta)\|_{H^1}\left(\|\nabla_x(a, b)\|_{H^1}^2+\left\|\nabla_x\{\mathbf{I}-\mathbf{P}\} f\right\|_{L_v^2\left(H_x^1\right)}^2\right),
\end{align*}
where $\varepsilon_3>0$ is sufficient small.
Similarly,
\begin{align*}
&\int_{\mathbb{R}^3} \partial^\alpha \partial_i \omega \partial^\alpha\left(\omega u_i+\sqrt{6} \theta b_i\right) \mathrm{d} x+\frac{3}{5}\int_{\mathbb{R}^3} \partial^\alpha \partial_i \omega \partial^\alpha \Upsilon_i(\mathfrak{l}+\mathfrak{r}+\mathfrak{s}) \mathrm{d} x \\
&\qquad +\frac{\sqrt{6}}{5} \int_{\mathbb{R}^3} \partial^\alpha \partial_i \omega \partial^\alpha\bigg(\sum_j \partial_j \Gamma_{i, j}\{\mathbf{I}-\mathbf{P}\} f\bigg) \mathrm{d} x \\
\,&\quad \leq\frac{2}{5}\left\|\partial^\alpha \partial_i \omega\right\|^2+C\|\{\mathbf{I}-\mathbf{P}\} f\|_{L_v^2\left(H_x^2\right)}^2 \\
&\qquad +C\|(u, \theta)\|_{H^1}^2\left(\|\nabla_x(b, \omega)\|_{H^1}^2+\left\|\nabla_x\{\mathbf{I}-\mathbf{P}\} f\right\|_{L_v^2\left(H_x^1\right)}^2\right) .
\end{align*}
Therefore, applying the above estimates to \eqref{G2.24}, and taking summation over $|\alpha| \leq 1$, $1 \leq i \leq 3$, we arrive at
\begin{align}\label{G2.25}
& \frac{\mathrm{d}}{\mathrm{d} t}\sum_{|\alpha|\leq 1}  \sum_{ i} \int_{\mathbb{R}^3} \partial^\alpha \partial_i \omega \partial^\alpha \Upsilon_i\{\mathbf{I}-\mathbf{P}\} f \mathrm{d} x+\sum_{|\alpha|\leq 1}\left\|\partial^\alpha \nabla_x \omega\right\|^2 \nonumber\\
\, &\quad  \leq\varepsilon_4\|\nabla_x b\|_{H^1}^2+C_{\varepsilon_4}\|\{\mathbf{I}-\mathbf{P}\} f\|_{L_v^2\left(H_x^2\right)}^2+C\|\sqrt{2} \omega-\sqrt{3} \theta\|_{H^1}^2 \nonumber\\
& \qquad +C\left(\|(u, \theta)\|_{H^1}+\|(u, \theta)\|_{H^1}^2\right)\left(\|\nabla_x(a, b, \omega)\|_{H^1}^2+\left\|\nabla_x\{\mathbf{I}-\mathbf{P}\} f\right\|_{L_v^2\left(H^1\right)}^2\right), 
\end{align}
where $\varepsilon_4>0$ is small enough.

We    achieve the following equality
by direct calculations from   \eqref{G2.14} to \eqref{G2.17},
\begin{align}\label{G2.26}
\left\|\partial^\alpha \partial_i a\right\|^2=\, & \frac{\mathrm{d}}{\mathrm{d} t} \int_{\mathbb{R}^3} \partial^\alpha \partial_i a \partial^\alpha\bigg(\frac{\sqrt{6}}{5} \Upsilon_i\{\mathbf{I}-\mathbf{P}\} f-b_i\bigg) \mathrm{d} x \nonumber\\
& -\int_{\mathbb{R}^3} \partial^\alpha \partial_i a_t \partial^\alpha\bigg(\frac{\sqrt{6}}{5} \Upsilon_i\{\mathbf{I}-\mathbf{P}\} f-b_i\bigg) \mathrm{d} x \nonumber\\
& +\int_{\mathbb{R}^3} \partial^\alpha \partial_i a \partial^\alpha\bigg(\left(u_i-b_i\right)+u_i a-\frac{2}{\sqrt{6}} \omega u_i-2 \theta b_i\nonumber \\
& \qquad -\frac{7}{5} \sum_j \partial_j \Gamma_{i j}\{\mathbf{I}-\mathbf{P}\} f-\frac{\sqrt{6}}{5} \Upsilon_i(\mathfrak{l}+\mathfrak{r}+\mathfrak{s})\bigg) \mathrm{d} x .
\end{align}
Next, we    obtain the following estimate by \eqref{G2.13}, 
\begin{align*}
& -\int_{\mathbb{R}^3} \partial^\alpha \partial_i a_t \partial^\alpha\bigg(\frac{\sqrt{6}}{5} \Upsilon_i\{\mathbf{I}-\mathbf{P}\} f-b_i\bigg) \mathrm{d} x \\
=&\int_{\mathbb{R}^3} \partial^\alpha\partial_i ({\rm div}_x  b) \partial^\alpha\bigg(b_i-\frac{\sqrt{6}}{5} \Upsilon_i\{\mathbf{I}-\mathbf{P}\} f\bigg) \mathrm{d} x\\
\leq& \int_{\mathbb{R}^3} \partial^\alpha ({\rm div}_x  b) \partial^\alpha \partial_i b_i \mathrm{d} x+\frac{1}{4}\left\|\partial^\alpha {\rm div}_x  b\right\|^2+C\|\{\mathbf{I}-\mathbf{P}\} f\|_{L_v^2\left(H_x^2\right)}^2 .
\end{align*}
We also have 
\begin{align*}
& \int_{\mathbb{R}^3} \partial^\alpha \partial_i a \partial^\alpha\bigg(\left(u_i-b_i\right)+u_i a-\frac{2}{\sqrt{6}} \omega u_i-2 \theta b_i  \\
& \quad -\frac{7}{5} \sum_j \partial_j \Gamma_{i j}\{\mathbf{I}-\mathbf{P}\} f-\frac{\sqrt{6}}{5} \Upsilon_i(\mathfrak{l}+\mathfrak{r}+\mathfrak{s})\bigg) \mathrm{d} x \\
\leq \,&\frac{1}{4}\|\partial^\alpha \partial_i a\|^2+C\Big(\left\|\partial^\alpha\left(u_i-b_i\right)\right\|^2+\|\partial^\alpha(u_i a-\frac{2}{\sqrt{6}} \omega u_i-2 \theta b_i)\|^2\Big. \\
& \Big.+\sum_j\left\|\partial^\alpha \partial_j \Gamma_{i j}\{\mathbf{I}-\mathbf{P}\} f\right\|^2+\left\|\partial^\alpha \Upsilon_i(\mathfrak{l}+\mathfrak{r}+\mathfrak{s})\right\|^2\Big) \\
\leq\, &\frac{1}{4}\left\|\partial^\alpha \partial_i a\right\|^2+C\|u-b\|_{H^1}^2+C\|\{\mathbf{I}-\mathbf{P}\} f\|_{L_v^2\left(H_x^2\right)}^2 \\
&+C\|(u, \theta)\|_{H^1}^2 \left(\|\nabla_x(a, b, \omega)\|_{H^1}^2+\left\|\nabla_x\{\mathbf{I}-\mathbf{P}\} f\right\|_{L_v^2\left(H_x^1\right)}^2\right) .
\end{align*}
Therefore, applying the above estimates to \eqref{G2.26}, and taking summation over $|\alpha| \leq 1$, one has
\begin{align}\label{G2.27}
& \frac{\mathrm{d}}{\mathrm{d} t} \sum_{|\alpha| \leq 1} \int_{\mathbb{R}^3} \frac{2}{21} \partial^\alpha a \partial^\alpha\left(\frac{\sqrt{6}}{5}\sum_i \Upsilon_i\{\mathbf{I}-\mathbf{P}\} f-{\rm div}_x  b\right) \mathrm{d} x+\frac{1}{14} \sum_{|\alpha| \leq 1}\left\|\partial^\alpha \nabla_x a\right\|^2\nonumber \\
\,&\quad \leq \frac{5}{42} \sum_{|\alpha| \leq 1}\left\|\partial^\alpha {\rm div}_x  b\right\|^2+C\|u-b\|_{H^1}^2+C\|\{\mathbf{I}-\mathbf{P}\} f\|_{L_v^2\left(H_x^2\right)}^2 \nonumber\\
&\qquad  +C\|(u, \theta)\|_{H^1}^2 \left(\|\nabla_x(a, b, \omega)\|_{H^1}^2+\left\|\nabla_x\{\mathbf{I}-\mathbf{P}\} f\right\|_{L_v^2\left(H_x^1\right)}^2\right) . 
\end{align}
Combining the estimates   
\eqref{G2.23} and \eqref{G2.25} with \eqref{G2.27},  and noticing the assumption \eqref{G2.1},
we    infer that (2.21) holds.
\end{proof}

\begin{lem}\label{L2.5}
For strong solutions to the problem \eqref{G1}--\eqref{G5}, there exists a positive constant $\lambda_4>0$, such that
\begin{align}\label{G2.28}
& \frac{\mathrm{d}}{\mathrm{d} t} \sum_{|\alpha| \leq 1} \int_{\mathbb{R}^3} \partial^\alpha u \cdot \partial^\alpha \nabla_x \rho \mathrm{d} x+\lambda_4\|\nabla_x \rho\|_{H^1}^2\nonumber \\
\, &\quad \leq C\left(\|u-b\|_{H^1}^2+\|\nabla_x u\|_{H^2}^2+\|\nabla_x \theta\|_{H^1}^2\right)\nonumber \\
&\qquad +C\left(\|(\rho,u,\theta)\|_{H^2}+\|(\rho,u,\theta)\|_{H^2}^2\right)\left(\|\nabla_x(a, \rho)\|_{H^1}^2+\|\nabla_x u\|_{H^2}^2\right),
\end{align}
holds for all $0 \leq t<T$. 
\end{lem}
\begin{proof}
Applying $\partial^\alpha(|\alpha| \leq 1)$ to \eqref{G3}, it holds
\begin{align}\label{G2.29}
\left\|\nabla_x \partial^\alpha \rho\right\|^2=\, & -\int_{\mathbb{R}^3} \nabla_x \partial^\alpha \rho \partial^\alpha \partial_t u \mathrm{d} x+\int_{\mathbb{R}^3} \nabla_x \partial^\alpha \rho \partial^\alpha\left(\frac{1}{1+\rho}(b-u)-\nabla_x \theta\right) \mathrm{d} x\nonumber \\
& +\int_{\mathbb{R}^3} \nabla_x \partial^\alpha \rho \partial^\alpha\left(-u \cdot \nabla_x u+\frac{\mu_1}{1+\rho} \Delta_x u+\frac{\rho-\theta}{1+\rho} \nabla_x \rho-\frac{1}{1+\rho} u a\right)\ \mathrm{d} x\nonumber\\
&+(\mu_1+\mu_2)\int_{\mathbb{R}^3}\nabla_x\partial^{\alpha}\rho\partial^{\alpha}\nabla_x({\rm div}u)     \mathrm{d} x\nonumber\\  
=:\, & \mathfrak{X}_1+\mathfrak{X}_2+\mathfrak{X}_3+\mathfrak{X}_4.
\end{align}
Using \eqref{G2} and the H\"{o}lder's, Sobolev's and Young's inequalities, we get
\begin{align*}
\mathfrak{X}_1= & -\frac{\mathrm{d}}{\mathrm{d} t} \int_{\mathbb{R}^3} \nabla_x \partial^\alpha \rho \partial^\alpha u \mathrm{d} x+\int_{\mathbb{R}^3} \partial^\alpha {\rm div}_x  u \partial^\alpha((1+\rho) {\rm div}_x  u+u \cdot \nabla_x \rho) \mathrm{d} x \\
& \leq-\frac{\mathrm{d}}{\mathrm{d} t} \int_{\mathbb{R}^3} \nabla_x \partial^\alpha \rho \partial^\alpha u \mathrm{d} x+C\left\|\partial^\alpha {\rm div}_x  u\right\|^2+C\|\rho\|_{H^2}\|\nabla u\|_{H^1}^2, \\
\mathfrak{X}_2 \leq \, & \frac{1}{16}\left\|\nabla_x \partial^\alpha \rho\right\|^2+C\left(\|u-b\|_{H^1}^2+\|\nabla_x \theta\|_{H^1}^2\right)+C\|\rho\|_{H^2}^2\|b-u\|_{H^1}^2, \\
\mathfrak{X}_3 \leq\, & \frac{1}{16}\left\|\nabla_x \partial^\alpha \rho\right\|^2+C\|u \cdot \nabla_x u\|_{H^1}^2+C\left\|\frac{1}{1+\rho} \Delta_x u\right\|_{H^1}^2 \\
& +C\left\|\frac{\rho-\theta}{1+\rho} \nabla_x \rho\right\|_{H^1}^2+C\left\|\frac{1}{1+\rho} u a\right\|_{H^1}^2 \\
\leq\,  & \frac{1}{16}\left\|\nabla_x \partial^\alpha \rho\right\|^2+C\|\nabla_x u\|_{H^2}^2+C\left(\|(\rho,u,\theta)\|_{H^2}+\|(\rho,u,\theta)\|_{H^2}^2\right)\|\nabla_x(a,\rho)\|_{H^1}^2\\
&+C\left(\|(\rho,u,\theta)\|_{H^2}+\|(\rho,u,\theta)\|_{H^2}^2\right)\|\nabla_x u\|_{H^2}^2,\\
\mathfrak{X}_4 \leq\, & \frac{1}{16}\left\|\nabla_x \partial^\alpha \rho\right\|^2+C\left\|\frac{1}{1+\rho} \nabla_x({\rm div}u)\right\|_{H^1}^2\\
\leq &\frac{1}{16}\left\|\nabla_x \partial^\alpha \rho\right\|^2+C\| \nabla_x u\|_{H^2}^2.
\end{align*}
According to the assumption \eqref{G2.1}, applying the above 
estimates to \eqref{G2.29}, we obtain \eqref{G2.28}.
\end{proof}
\subsubsection{Energy estimates for mixed space-velocity derivatives} 
In this subsection, we show the energy estimates 
with regard to the mixed space-velocity derivatives
of $f$, i.e.  $\partial_\beta^\alpha f$. 
Using the properties of $\mathbf{P}$, 
we obtain 
\begin{align*}
 \{\mathbf{I}-\mathbf{P}\}(v \sqrt{M})=\,&0, \\
\{\mathbf{I}-\mathbf{P}\}\Big(\left(|v|^2-3\right) \sqrt{M}\Big)=\,&0, \\
\mathbf{L}\{\mathbf{I}-\mathbf{P}\} f-\{\mathbf{I}-\mathbf{P}\} \mathbf{L} f=\,&0,\\
g \cdot\{\mathbf{I}-\mathbf{P}\} f-\mathbf{P}(g \cdot\{\mathbf{I}-\mathbf{P}\} f)+\{\mathbf{I}-\mathbf{P}\}(g \cdot \mathbf{P} f)=\,&\{\mathbf{I}-\mathbf{P}\}(g \cdot f).
\end{align*}
Thus, applying $ \mathbf{I}-\mathbf{P} $ to both sides of \eqref{G1}, one gets
\begin{align}\label{G2.30}
&\partial_t\{\mathbf{I} -\mathbf{P}\} f+v \cdot \nabla_x\{\mathbf{I}-\mathbf{P}\} f+u \cdot \nabla_v\{\mathbf{I}-\mathbf{P}\} f-\frac{1}{2} u \cdot v\{\mathbf{I}-\mathbf{P}\} f \nonumber\\
& \qquad -\frac{\theta }{\sqrt{M}} \Delta_v\left(\sqrt{M}\{\mathbf{I}-\mathbf{P}\} f\right)-\mathbf{L}\{\mathbf{I}-\mathbf{P}\} f\nonumber \\
&\quad = \mathbf{P}\Big(v \cdot \nabla_x\{\mathbf{I}-\mathbf{P}\} f+u \cdot \nabla_v\{\mathbf{I}-\mathbf{P}\} f-\frac{1}{2} u \cdot v\{\mathbf{I}-\mathbf{P}\} f\Big. \nonumber\\
&\qquad  \Big.-\frac{\theta}{\sqrt{M}} \Delta_v(\sqrt{M}\{\mathbf{I}-\mathbf{P}\} f)\Big)-\{\mathbf{I}-\mathbf{P}\}\Big(v \cdot \nabla_x \mathbf{P} f+u \cdot \nabla_v \mathbf{P} f\Big.\nonumber \\
& \qquad \Big.-\frac{1}{2} u \cdot v \mathbf{P} f-\frac{\theta }{\sqrt{M}} \Delta_v(\sqrt{M} \mathbf{P} f)\Big) .
\end{align}

\begin{lem}\label{L2.6}
Let $1 \leq k \leq 2$. For strong solutions to the problem \eqref{G1}--\eqref{G5}, there exists a  positive constant     $\lambda_5$, such that
\begin{align}\label{G2.31}
& \frac{\mathrm{d}}{\mathrm{d} t} \sum_{\substack{|\beta|=k \\
|\alpha|+|\beta| \leq 2}}\left\|\partial_\beta^\alpha\{\mathbf{I}-\mathbf{P}\} f\right\|^2+\lambda_5 \sum_{\substack{|\beta|=k \\
|\alpha|+|\beta| \leq 2}}\left\|\partial_\beta^\alpha\{\mathbf{I}-\mathbf{P}\} f\right\|_\nu^2 \nonumber\\
\,&\quad \leq C \big(\|(u, \theta)\|_{H^2}+\|(u, \theta)\|_{H^2}^2\big)\nonumber\\
&\qquad \times \bigg( \sum_{\left|\alpha\right| \leq 2}\left\|\partial^{\alpha}\{\mathbf{I}-\mathbf{P}\} f\right\|_\nu^2+ \sum_{\substack{1 \leq\left|\beta\right| \leq 2 \\
\left|\alpha\right|+\left|\beta\right| \leq 2}}\left\|\partial_{\beta}^{\alpha}\{\mathbf{I}-\mathbf{P}\} f\right\|_\nu^2+\|\nabla_x(a,b, \omega)\|_{H^{1}}^2   \bigg) \nonumber\\
& \qquad + C \chi_{\{k=2\}} \sum_{\substack{1 \leq\left|\beta^{\prime}\right| \leq k-1 \\
\left|\alpha^{\prime}\right|+\left|\beta^{\prime}\right| \leq 2}}\|\partial_{\beta^{\prime}}^{\alpha^{\prime}}\{\mathbf{I}-\mathbf{P}\} f\|_\nu^2
+C\|\nabla_x(a,b, \omega)\|_{H^{1}}^2+C\sum_{\left|\alpha\right| \leq 2}\left\|\partial^{\alpha}\{\mathbf{I}-\mathbf{P}\} f\right\|_\nu^2, 
\end{align}
holds for all $0 \leq t<T$. 
Here $\chi_B$ is the characteristic function on the set $B$.
\end{lem}
\begin{proof}
Fix $k$ with $1 \leq k \leq 2$. Select $\alpha$ and  $\beta$ 
  satisfying $|\beta|=k$ and $|\alpha|+|\beta| \leq 2$. 
Multiplying \eqref{G2.30} by $\partial_\beta^\alpha\{\mathbf{I}-\mathbf{P}\} f$ and taking integration, we have
\begin{equation}\label{G2.32}
\frac{1}{2} \frac{\mathrm{d}}{\mathrm{d} t}\left\|\partial_\beta^\alpha\{\mathbf{I}-\mathbf{P}\} f\right\|^2+\int_{\mathbb{R}^3}\left\langle-L \partial_\beta^\alpha\{\mathbf{I}-\mathbf{P}\} f, \partial_\beta^\alpha\{\mathbf{I}-\mathbf{P}\} f\right\rangle \mathrm{d} x =:\sum_{i=1}^7 \mathfrak{M}_i,
\end{equation}
where
\begin{align*}
 \mathfrak{M}_1=\,&\int_{\mathbb{R}^3}\left\langle-\partial_x^\alpha[\partial_v^\beta, v \cdot \nabla_x]\{\mathbf{I}-\mathbf{P}\} f, \partial_\beta^\alpha\{\mathbf{I}-\mathbf{P}\} f\right\rangle \mathrm{d} x, \\
 \mathfrak{M}_2=\,&\int_{\mathbb{R}^3}\left\langle\partial_x^\alpha[\partial_v^\beta,-|v|^2]\{\mathbf{I}-\mathbf{P}\} f, \partial_\beta^\alpha\{\mathbf{I}-\mathbf{P}\} f\right\rangle \mathrm{d} x, \\
  \mathfrak{M}_3=\,&\int_{\mathbb{R}^3}\left\langle-\partial_\beta^\alpha\left(u \cdot \nabla_v\{\mathbf{I}-\mathbf{P}\} f\right), \partial_\beta^\alpha\{\mathbf{I}-\mathbf{P}\} f\right\rangle \mathrm{d} x, \\
  \mathfrak{M}_4=\,&\int_{\mathbb{R}^3}\left\langle\frac{1}{2} \partial_\beta^\alpha(u \cdot v\{\mathbf{I}-\mathbf{P}\} f), \partial_\beta^\alpha\{\mathbf{I}-\mathbf{P}\} f\right\rangle \mathrm{d} x,\\
 \mathfrak{M}_5=\,&\int_{\mathbb{R}^3}\left\langle\partial_\beta^\alpha\left(\frac{\theta}{\sqrt{M}}  \Delta_v(\sqrt{M}\{\mathbf{I}-\mathbf{P}\} f)\right), \partial_\beta^\alpha\{\mathbf{I}-\mathbf{P}\} f\right\rangle \mathrm{d} x ,\\
  \mathfrak{M}_6=\,& \int_{\mathbb{R}^3}\left\langle\partial _ { \beta } ^ { \alpha } \mathbf { P } \bigg( v \cdot \nabla_x\{\mathbf{I}-\mathbf{P}\} f+u \cdot \nabla_v\{\mathbf{I}-\mathbf{P}\} f\right.\bigg. \\
& \quad \bigg.\left.-\frac{1}{2} u \cdot v\{\mathbf{I}-\mathbf{P}\} f-\frac{\theta}{\sqrt{M}}\Delta_v(\sqrt{M}\{\mathbf{I}-\mathbf{P}\} f)\bigg), \partial_\beta^\alpha\{\mathbf{I}-\mathbf{P}\} f\right\rangle \mathrm{d} x, \\
 \mathfrak{M}_7=\,&  \int_{\mathbb{R}^3}\left\langle-\partial_\beta^\alpha\{\mathbf{I}-\mathbf{P}\}\left(v \cdot \nabla_x \mathbf{P} f+u \cdot \nabla_v \mathbf{P} f-\frac{1}{2} u \cdot v \mathbf{P} f\right.\right. \\
& \quad \left.\left.-\frac{\theta}{\sqrt{M}} \Delta_v(\sqrt{M} \mathbf{P} f)\right), \partial_\beta^\alpha\{\mathbf{I}-\mathbf{P}\} f\right\rangle \mathrm{d} x.
\end{align*}
Here the fact that  $[\partial_v^\beta, \mathbf{L}]=[\partial_v^\beta,-|v|^2]$ has been used.

Below we   estimate the  terms $\mathfrak{M}_i ( i=1,\dots, 7)$ in (2.32) one by one as follows:
\begin{align*}
\mathfrak{M}_1 \leq\, & \varepsilon_5\left\|\partial_\beta^\alpha\{\mathbf{I}-\mathbf{P}\} f\right\|^2+C_{\varepsilon_5}\|[\partial_v^\beta, v \cdot \nabla_v] \partial_x^\alpha\{\mathbf{I}-\mathbf{P}\} f\|^2 \\
 \leq \, & \varepsilon_5\left\|\partial_\beta^\alpha\{\mathbf{I}-\mathbf{P}\} f\right\|^2+C_{\varepsilon_5} \sum_{\left|\alpha\right| \leq 2-k}\left\|\partial^{\alpha} \nabla_x\{\mathbf{I}-\mathbf{P}\} f\right\|^2 \\
& +\chi_{\{k=2\}} C_{\varepsilon_5} \sum_{\substack{1 \leq\left|\beta^{\prime}\right| \leq k-1 \\
\left|\alpha^{\prime}\right|+\left|\beta^{\prime}\right| \leq 2}}\|\partial_{\beta^{\prime}}^{\alpha^{\prime}}\{\mathbf{I}-\mathbf{P}\} f\|^2, \\
\mathfrak{M}_2 \leq\, &\varepsilon_6\left\|\partial_\beta^\alpha\{\mathbf{I}-\mathbf{P}\} f\right\|^2+C_{\varepsilon_6}\|[\partial_v^\beta,-|v|^2] \partial_x^\alpha\{\mathbf{I}-\mathbf{P}\} f\|^2 \\
\quad\leq \,& \varepsilon_6\left\|\partial_\beta^\alpha\{\mathbf{I}-\mathbf{P}\} f\right\|^2+C_{\varepsilon_6} \sum_{\left|\alpha\right| \leq 2-k}\left\|\partial^{\alpha}\{\mathbf{I}-\mathbf{P}\} f\right\|_\nu^2 \\
& +\chi_{\{k=2\}} C_{\varepsilon_6}\sum_{\substack{1 \leq | \beta^{\prime}|\leq k-1\\
| \alpha^{\prime}|+| \beta^{\prime} \mid \leq 2}}\|\partial_{\beta^{\prime}}^{\alpha^{\prime}}\{\mathbf{I}-\mathbf{P}\} f\|_\nu^2, \\
\mathfrak{M}_3 \leq \,& \varepsilon_7\left\|\partial_\beta^\alpha\{\mathbf{I}-\mathbf{P}\} f\right\|^2+C_{\varepsilon_7}\|\partial_x^\alpha(u \cdot \nabla_v \partial_v^\beta\{\mathbf{I}-\mathbf{P}\} f)\|^2 \\
\leq\,  &  \varepsilon_7\left\|\partial_\beta^\alpha\{\mathbf{I}-\mathbf{P}\} f\right\|^2+C_{\varepsilon_7}\|u\|_{H^1}^2 \sum_{\substack{1 \leq\left|\beta^{\prime}\right| \leq 2 \\
\left|\alpha^{\prime}\right|+\left|\beta^{\prime}\right| \leq 2}}\|\partial_{\beta^{\prime}}^{\alpha^{\prime}}\{\mathbf{I}-\mathbf{P}\} f\|^2, \\
\mathfrak{M}_4 \leq \, &\varepsilon_8\left\|\partial_\beta^\alpha\{\mathbf{I}-\mathbf{P}\} f\right\|^2+C_{\varepsilon_8}\|\partial_x^\alpha\big(u \cdot \partial_v^\beta(v\{\mathbf{I}-\mathbf{P}\} f)\big)\|^2 \\
\quad\leq \, &\varepsilon_8\left\|\partial_\beta^\alpha\{\mathbf{I}-\mathbf{P}\} f\right\|^2+C_{\varepsilon_8}\|u\|_{H^1}^2 \sum_{\substack{1 \leq\left|\beta^{\prime}\right| \leq 2 \\
\left|\alpha^{\prime}\right|+\left|\beta^{\prime}\right| \leq 2}}\|\partial_{\beta^{\prime}}^{\alpha^{\prime}}\{\mathbf{I}-\mathbf{P}\} f\|_\nu^2 \\
& +C_{\varepsilon_8}\|u\|_{H^1}^2 \sum_{\left|\alpha^{\prime}\right| \leq 2-k}\|\partial^{\alpha^{\prime}}\{\mathbf{I}-\mathbf{P}\} f\|^2, \\
\mathfrak{M}_5 \leq\, & \varepsilon_9\left\|\partial_\beta^\alpha\{\mathbf{I}-\mathbf{P}\} f\right\|^2+C_{\varepsilon_9}\|\theta\|_{H^1}^2 \sum_{\substack{1 \leq\left|\beta^{\prime}\right| \leq 2 \\
\left|\alpha^{\prime}\right|+\left|\beta^{\prime}\right| \leq 2}}\|\partial_{\beta^{\prime}}^{\alpha^{\prime}}\{\mathbf{I}-\mathbf{P}\} f\|_\nu^2 \\
& + C_{\varepsilon_9}\|\theta\|_{H^3}^2 \sum_{\left|\alpha^{\prime}\right| \leq 2-k+1}\|\partial^{\alpha^{\prime}}\{\mathbf{I}-\mathbf{P}\} f\|_\nu^2, \\
\mathfrak{M}_6 \leq\, & \varepsilon_{10}\left\|\partial_\beta^\alpha\{\mathbf{I}-\mathbf{P}\} f\right\|^2+C_{\varepsilon_{10}}\| \partial_\beta^\alpha \mathbf{P}(v \cdot \nabla_x\{\mathbf{I}-\mathbf{P}\} f) \|^2\\
& +C_{\varepsilon_{10}}\| \partial_\beta^\alpha \mathbf{P}(u \cdot \nabla_v\{\mathbf{I}-\mathbf{P}\} f) \|^2+C_{\varepsilon_{10}}\| \partial_\beta^\alpha \mathbf{P}(u \cdot v\{\mathbf{I}-\mathbf{P}\} f) \|^2\\
&  +C_{\varepsilon_{10}}\Big\|\partial_\beta^\alpha \mathbf{P}\Big(\frac{\theta }{\sqrt{M}}\Delta_v(\sqrt{M}\{\mathbf{I}-\mathbf{P}\} f)\Big)\Big\|^2 \\
\leq\, &   {\varepsilon_{10}}\left\|\partial_\beta^\alpha\{\mathbf{I}-\mathbf{P}\} f\right\|^2
+C_{\varepsilon_{10}} \sum_{|\alpha^{\prime}| \leq 2-k}\|\nabla_x \partial^{\alpha^{\prime}}\{\mathbf{I}-\mathbf{P}\} f\|^2 \\
&  +C_{\varepsilon_{10}}\|(u, \theta)\|_{H^1}^2 \sum_{\left|\alpha^{\prime}\right| \leq 2-k}\|\partial^{\alpha^{\prime}}\{\mathbf{I}-\mathbf{P}\} f\|^2, \\
\mathfrak{M}_7 \leq\, &  \varepsilon_{11}\left\|\partial_\beta^\alpha\{\mathbf{I}-\mathbf{P}\} f\right\|^2+C_{\varepsilon_{11}}\left\|\partial_\beta^\alpha\{\mathbf{I}-\mathbf{P}\}\left(v \cdot \nabla_x \mathbf{P} f\right)\right\|^2 \\
&  +C_{\varepsilon_{11}}\left\|\partial_\beta^\alpha\{\mathbf{I}-\mathbf{P}\}\left(u \cdot \nabla_v \mathbf{P} f\right)\right\|^2+C_{\varepsilon_{11}}\left\|\partial_\beta^\alpha\{\mathbf{I}-\mathbf{P}\}(u \cdot v \mathbf{P} f)\right\|^2 \\
&  +C_{\varepsilon_{11}}\Big\|\partial_\beta^\alpha \{\mathbf{I}-\mathbf{P}\}\Big(\frac{\theta}{\sqrt{M}} \Delta_v(\sqrt{M} \mathbf{P} f)\Big)\Big\|^2 \\
\leq \,& \varepsilon_{11}\left\|\partial_\beta^\alpha\{\mathbf{I}-\mathbf{P}\} f\right\|^2+C_{\varepsilon_{11}}\|\nabla_x(a,b, \omega)\|_{H^{1}}^2 
+C_{\varepsilon_{11}}\|(u, \theta)\|_{H^{2}}^2\|\nabla_x(a,b, \omega)\|_{H^1}^2,
\end{align*}
where $\varepsilon_i>0 (i=5,\dots,11)$ are sufficient small.
Thanks to  \eqref{G1.121}, we have
\begin{align*}
  \int_{\mathbb{R}^3}\left\langle- \mathbf{L} \partial_\beta^\alpha\{\mathbf{I}-\mathbf{P}\} f, \partial_\beta^\alpha\{\mathbf{I}-\mathbf{P}\} f\right\rangle \mathrm{d} x  
\geq\, &\lambda_0\left\|\left\{\mathbf{I}-\mathbf{P}_0\right\} \partial_\beta^\alpha\{\mathbf{I}-\mathbf{P}\} f\right\|_\nu^2 \\
\geq\, & \frac{\lambda_0}{2}\left\|\partial_\beta^\alpha\{\mathbf{I}-\mathbf{P}\} f\right\|_\nu^2-\lambda_0\left\|\mathbf{P}_0 \partial_\beta^\alpha\{\mathbf{I}-\mathbf{P}\} f\right\|_\nu^2 \\
\geq\, &\frac{\lambda_0}{2}\left\|\partial_\beta^\alpha\{\mathbf{I}-\mathbf{P}\} f\right\|_\nu^2-C\left\|\partial^\alpha\{\mathbf{I}-\mathbf{P}\} f\right\|^2.
\end{align*}
 Applying all the above estimates to \eqref{G2.32} and choosing $\varepsilon_i (i=5,\dots,11)$ sufficient small, we infer   \eqref{G2.31} holds.
\end{proof}

\begin{rem}
 By selecting some reasonable 
constants $C_k$ together with the above lemma, we obtain  that
there exists a positive constant $\lambda_6>0$   such that
\begin{align}\label{G2.33}
& \frac{\mathrm{d}}{\mathrm{d} t} \sum_{1 \leq k \leq 2} C_k \sum_{\substack{|\beta|=k \\
|\alpha|+|\beta| \leq 2}}\left\|\partial_\beta^\alpha\{\mathbf{I}-\mathbf{P}\} f\right\|^2+\lambda_6 \sum_{\substack{1 \leq|\beta| \leq 2 \\
|\alpha|+|\beta| \leq 2}}\left\|\partial_\beta^\alpha\{\mathbf{I}-\mathbf{P}\} f\right\|_\nu^2\nonumber \\
\,&\quad \leq C\big(\|(u, \theta)\|_{H^2}+\|(u, \theta)\|_{H^2}^2\big)\nonumber \\
&\qquad\times \Bigg( \sum_{\left|\alpha\right| \leq 2}\left\|\partial^{\alpha}\{\mathbf{I}-\mathbf{P}\} f\right\|_\nu^2+ \sum_{\substack{1 \leq\left|\beta\right| \leq 2 \\
\left|\alpha\right|+\left|\beta\right| \leq 2}}\left\|\partial_{\beta}^{\alpha}\{\mathbf{I}-\mathbf{P}\} f\right\|_\nu^2+\|\nabla_x(a,b, \omega)\|_{H^{1}}^2   \Bigg)\nonumber\\
&\qquad+ C\bigg(\|\nabla_x(a,b, \omega)\|_{H^1}^2+\sum_{|\alpha| \leq 2}\left\|\partial^\alpha\{\mathbf{I}-\mathbf{P}\} f\right\|_\nu^2\bigg). 
\end{align}  
\end{rem}

\smallskip 
Now  we   define the total energy functional $\mathcal{E}(t)$ and the 
dissipation rate $\mathcal{D}(t)$ as follows:
\begin{align}
\mathcal{E}(t) & :=\|(f, \rho, u, \theta)\|^2+\sum_{1 \leq|\alpha| \leq 2}\left(\left\|\partial^\alpha f\right\|^2+\left\|\partial^\alpha \rho\right\|^2+\left\|\partial^\alpha u\right\|^2+\left\|\partial^\alpha \theta\right\|^2\right)+\tau_1\mathcal{E}_0(t)\nonumber\\
&\quad\ +\tau_2\sum_{|\alpha| \leq 1} \int_{\mathbb{R}^3} \partial^\alpha u \cdot \partial^\alpha \nabla_x \rho \mathrm{d} x+\tau_3\sum_{1 \leq k \leq 2} C_k \sum_{\substack{|\beta|=k \\
|\alpha|+|\beta| \leq 2}}\left\|\partial_\beta^\alpha\{\mathbf{I}-\mathbf{P}\} f\right\|^2,\\
\mathcal{D}(t) & :=\|\nabla_x(a,b,\rho,\omega)\|_{H^1}^2+ \sum_{|\alpha| \leq 2}\left(\left\|\{\mathbf{I}-\mathbf{P}\} \partial^\alpha f\right\|_\nu^2+\left\|\nabla_x \partial^\alpha(u, \theta)\right\|^2\right)\nonumber\\
&\quad\ +\left\|b-u\right\|_{H^2}^2+\|\sqrt{2} \omega-\sqrt{3} \theta\|_{H^2}^2+\sum_{\substack{1 \leq | \beta|\leq 2\\
| \alpha|+| \beta \mid \leq 2}}\left\|\partial_\beta^\alpha\{\mathbf{I}-\mathbf{P}\} f\right\|_\nu^2,
\end{align}
where the constants $0<\tau_1,\tau_2,\tau_3\ll 1$ will to be determined later.
Noting that $\tau_1,\tau_2,\tau_3>0 $ are sufficient small, with  the assumption \eqref{G2.1}, we know that 
$\mathcal{E}(t)\sim\|(\rho,u,\theta,f)\|_{H^2}^{2}$.

Moreover, by  appropriately selecting constants $0<\tau_3\ll\tau_1\ll 1$, and  summing the  equations \eqref{G2.5}, \eqref{G2.7}, $\tau_1\times$\eqref{G2.21}, $\tau_2\times$\eqref{G2.28} and $\tau_3\times$\eqref{G2.33}, we get 
\begin{equation}\label{G2.36}
\frac{\mathrm{d}}{\mathrm{d} t} \mathcal{E}(t)+\lambda_7 \mathcal{D}(t) \leq C\left(\mathcal{E}^{\frac{1}{2}}(t)+\mathcal{E}(t)+\mathcal{E}^2(t)\right) \mathcal{D}(t) \leq C\left(\delta+\delta^2+\delta^4\right) \mathcal{D}(t),
\end{equation}
for some $\lambda_7>0$.
Therefore, as long as $0<\delta<1$ is sufficient small, \eqref{G2.36} implies that  there exists a positive constant $\lambda_8>0$, such that
\begin{equation}\label{G2.37}
\mathcal{E}(t)+\lambda_8 \int_0^t \mathcal{D}(s) \mathrm{d} s \leq \mathcal{E}(0),
\end{equation}
for all $0 \leq t<T$. Besides, \eqref{G2.1}  can  be proved by selecting $\mathcal{E}(0) \sim\left\|f_0\right\|_{H_{x, v}^2}^2+$ $\left\|\left(\rho_0, u_0, \theta_0\right)\right\|_{H^2}^2$ appropriately.

\subsection{Local existence of strong solutions} 
In this subsection, we shall show the 
 local existence and   uniqueness of strong solutions
to the problem \eqref{G1}--\eqref{G5}. 
To archive this goal, we shall establish the uniform energy
estimates for the iteration sequence of
approximate solutions. 
We consider the following linear equations 
with regard to the sequence 
$\left\{(f^n, \rho^n, u^n, \theta^n)\right\}_{n
\geq 0}$:
\begin{align}\label{G2.38}
& \partial_t f^{n+1}+v \cdot \nabla_x f^{n+1}-\mathbf{L} f^{n+1}=-u^n \cdot\left(\nabla_v f^{n+1}-\frac{v}{2} f^{n+1}-v \sqrt{M}\right)\nonumber \\
& \quad+\frac{\theta^n}{\sqrt{M}}  \Delta_v\big(M+\sqrt{M} f^{n+1}\big), \\
\label{G2.39}
& \partial_t \rho^{n+1}+u^n \cdot \nabla_x \rho^{n+1}+\left(1+\rho^n\right) \operatorname{div} u^{n+1}=0 ,\\
\label{G2.40}
& \partial_t u^{n+1}-\frac{\mu_1}{1+\rho^n} \Delta_x u^{n+1}=-u^n \cdot \nabla_x u^{n+1}-\nabla_x \theta^{n+1}-\frac{1+\theta^n}{1+\rho^n} \nabla_x \rho^{n+1} \nonumber\\
& \quad+\frac{1}{1+\rho^n}\left(b^{n+1}-u^{n+1}-u^{n+1} a^n+(\mu_1+\mu_2)\nabla_x({\rm div}_x  u^{n+1})\right),\ \\
\label{G2.41}
& \partial_t \theta^{n+1}-\frac{\kappa}{1+\rho^n} \Delta_x \theta^{n+1}=-u^n \cdot \nabla_x \theta^{n+1}-\theta^{n+1} \operatorname{div} u^n-\frac{\sqrt{3}}{1+\rho^n}\big(\sqrt{2} \omega^{n+1}-\sqrt{3} \theta^{n+1}\big) \nonumber\\
& \quad-\operatorname{div} u^{n+1}+\frac{1}{1+\rho^n}\big(|u^n|^2-2u^n \cdot b^n+a^n\left|u^n\right|^2-3 a^n \theta^n+2\mu_1|D(u^n)|^2+\mu_2|{\rm div}_x  u^n|^2\big),
\end{align}
with $\left(f^{n+1}(t,x,v), \rho^{n+1}(t,x), u^{n+1}(t,x), \theta^{n+1}(t,x)\right)\big|_{t=0}=\left(f_0, \rho_0, u_0, \theta_0\right)$, and   $n=0,1,2, \dots$.

Next, we construct the following solutions in Banach space:
$$
\mathcal{X}(0, T ; B):=\left\{\begin{array}{c}
f \in C\left([0, T]; {H_{x, v}^2}\right), \\
\rho \in C\left([0, T], H^2\right) \cap C^1\left([0, T], H^1\right) ; \\
u \in C\left([0, T], H^2\right) \cap C^1\left([0, T], L^2\right) ; \\
\theta \in C\left([0, T], H^2\right) \cap C^1\left([0, T], L^2\right) ; \\
\sup\limits_{0 \leq t \leq T}\big\{\|f(t)\|_{H_{x, v}^2}+\|(\rho, u, \theta)(t)\|_{H^2}\big\} \leq B ; \\
\rho_1=\frac{1}{2}\Big(-1+\inf\limits_{  x \in \mathbb{R}^3} \rho(0, x)\Big)>-1 ; \\
\inf\limits_{0 \leq t \leq T,\, x \in \mathbb{R}^3} \rho(t, x) \geq \tilde\varrho>0.
\end{array}\right\}.
$$
Now we give our result on  local existence of strong solutions as follows:
\begin{thm}\label{T2.7}
There exist $B_0>0$ and $T^*>0$, such that if $\left\|f_0\right\|_{H_{x, v}^2}+\left\|\left(\rho_0, u_0, \theta_0\right)\right\|_{H^2}<\frac{B_0}{3}$, then for each $n \geq 1,\left\{(f^n, \rho^n, u^n, \theta^n\right)\}$ is well-defined with
$$
\left(f^n, \rho^n, u^n, \theta^n\right) \in \mathcal{X}\left(0, T^* ; B_0\right) .
$$
Moreover, we   achieve:

$(1)$ $\{\left(f^n, \rho^n, u^n, \theta^n\right)\}_{n \geq 0}$ is a Cauchy sequence in the Banach space $C\left(\left[0, T^*\right], H_{x,v}^1\right)$,

$(2)$ the limit of $\{\left(f^n, \rho^n, u^n, \theta^n\right)\}$, denoted by 
  $(f, \rho, u, \theta)$, belongs to $ \mathcal{X}\left(0, T^* ; B_0\right)$,

$(3)$ $(f, \rho, u, \theta)$ satisfies the system \eqref{G1}-\eqref{G5},

$(4)$ $(f, \rho, u, \theta)$ is the unique strong solution of  \eqref{G1}-\eqref{G5} in $\mathcal{X}\left(0, T^* ; B_0\right)$.
\end{thm}

\begin{proof}
Let $T^*>0$ be a constant which will be determined later.
For brevity, we    suppose that $\left(f^n, \rho^n, u^n, \theta^n\right)$ 
is smooth enough. Overwise,
we  can   discuss the following problem:
\begin{align*}
&\partial_t f^{n+1, \varepsilon}+v \cdot \nabla_x f^{n+1, \varepsilon}-\mathbf{L} f^{n+1, \varepsilon}\\
& \quad=-u^{n, \varepsilon} \cdot\left(\nabla_v f^{n+1, \varepsilon}-\frac{v}{2} f^{n+1, \varepsilon}-v \sqrt{M}\right)+\frac{\theta^{n, \varepsilon} }{\sqrt{M}} \Delta_v\left(M+\sqrt{M} f^{n+1, \varepsilon}\right), \\
&\partial_t \rho^{n+1, \varepsilon}+u^{n, \varepsilon} \cdot \nabla_x \rho^{n+1, \varepsilon}+\left(1+\rho^{n, \varepsilon}\right) \operatorname{div} u^{n+1, \varepsilon}=0, \\
&  \partial_t u^{n+1, \varepsilon}-\frac{\mu_1}{1+\rho^{n, \varepsilon}} \Delta_x u^{n+1, \varepsilon}\\
&  \quad=-u^{n, \varepsilon} \cdot \nabla_x u^{n+1, \varepsilon}-\nabla_x \theta^{n+1, \varepsilon}-\frac{1+\theta^{n, \varepsilon}}{1+\rho^{n, \varepsilon}} \nabla_x \rho^{n+1, \varepsilon}\\
& \quad \quad +\frac{1}{1+\rho^{n, \varepsilon}}\left(b^{n+1, \varepsilon}-u^{n+1, \varepsilon}-u^{n+1, \varepsilon} a^{n, \varepsilon}+(\mu_1+\mu_2)\nabla_x({\rm div}_x  u^{n+1,\varepsilon})\right), \\
&  \partial_t \theta^{n+1, \varepsilon}-\frac{\kappa}{1+\rho^{n, \varepsilon}} \Delta_x \theta^{n+1, \varepsilon}\\
& \quad  =-u^{n, \varepsilon} \cdot \nabla_x \theta^{n+1, \varepsilon}-\theta^{n+1, \varepsilon} \operatorname{div} u^{n, \varepsilon}-\operatorname{div} u^{n+1, \varepsilon}-\frac{\sqrt{3}}{1+\rho^{n, \varepsilon}}
\big(\sqrt{2} \omega^{n+1, \varepsilon}-\sqrt{3} \theta^{n+1, \varepsilon}\big)\\
&   \qquad+\frac{1}{1+\rho^{n, \varepsilon}} \big(|u^{n, \varepsilon}|^2-2u^{n, \varepsilon} \cdot b^{n, \varepsilon}+a^{n, \varepsilon}\left|u^{n, \varepsilon}\right|^2+2\mu_1 |D(u^{n, \varepsilon})|^2- 3a^{n, \varepsilon} \theta^{n, \varepsilon}+\mu_2|{\rm div}u^{n,\varepsilon}|^2\big), \\
&  (f^{n+1, \varepsilon}, \rho^{n+1, \varepsilon}, u^{n+1, \varepsilon}, \theta^{n+1, \varepsilon})|_{t=0}=\left(f_0^{\varepsilon}, \rho_0^{\varepsilon}, u_0^{\varepsilon}, \theta_0^{\varepsilon}\right), 
\end{align*}
for all $\varepsilon>0$, where $\left(f_0^{\varepsilon}, \rho_0^{\varepsilon}, u_0^{\varepsilon}, \theta_0^{\varepsilon}\right)$ 
is a smooth approximation of $\left(f_0, \rho_0, u_0, \theta_0\right)$.
Through letting $\varepsilon \rightarrow 0$, 
our supposition is reasonable.

Applying $\partial_x^\alpha$ with $|\alpha| \leq 2$ to the equation \eqref{G2.38}, multiplying the result by $\partial_x^\alpha f^{n+1}$, and then taking integration, we get
\begin{align}\label{G2.42}
&\ \quad \frac{1}{2} \frac{\mathrm{d}}{\mathrm{d} t}\|\partial^\alpha f^{n+1}\|^2+\int_{\mathbb{R}^3}\left\langle-\mathbf{L} \partial^\alpha f^{n+1}, \partial^\alpha f^{n+1}\right\rangle \mathrm{d} x \nonumber\\
&=-\iint_{\mathbb{R}^3\times\mathbb{R}^3} \partial^\alpha\left(\frac{u^n}{\sqrt{M}}  \nabla_v\big(M+\sqrt{M} f^{n+1}\big)\right) \partial^\alpha f^{n+1} \mathrm{d} x \mathrm{d} v \nonumber\\
&\ \quad+\iint_{\mathbb{R}^3\times\mathbb{R}^3} \partial^\alpha\left(\frac{\theta^n}{\sqrt{M}} \Delta_v\big(M+\sqrt{M} f^{n+1}\big)\right) \partial^\alpha f^{n+1} \mathrm{d} x \mathrm{d} v \nonumber\\
 & \leq C\|(u^n, \theta^n)\|_{H^2}\|f^{n+1}\|_{H_{x,v}^2}+C\|u^n\|_{H^2}\|f^{n+1}\|_{H_{x,v}^2}\|\partial^\alpha f^{n+1}\|_\nu \nonumber\\
 &\quad\ +C\|\theta^n\|_{H^2}\sum_{\left|\alpha\right| \leq 2}\|\partial^{\alpha} f^{n+1}\|_\nu\|\partial^\alpha f^{n+1}\|_\nu .
\end{align}
Note that
\begin{align*}
\int_{\mathbb{R}^3}\left\langle-\mathbf{L} \partial^\alpha f^{n+1}, \partial^\alpha f^{n+1}\right\rangle \mathrm{d} x \geq \lambda_0\left\|\left\{\mathbf{I}-\mathbf{P}_0\right\} \partial^\alpha f^{n+1}\right\|_\nu^2 .
\end{align*}
By combining $\tilde\lambda \|\mathbf{P}_0 \partial^\alpha f^{n+1}\|_\nu^2$ 
for some $\tilde \lambda>0$ with \eqref{G2.42}, we obtain that, for $|\alpha| \leq 2$,  
\begin{align}
& \frac{1}{2} \frac{\mathrm{d}}{\mathrm{d} t} \sum_{|\alpha| \leq 2}\|\partial^\alpha f^{n+1}\|^2+\lambda_9 \sum_{|\alpha| \leq 2}\|\partial^\alpha f^{n+1}\|_\nu^2\nonumber \\
& \quad  \leq C\|f^{n+1}\|_{{H_{x,v}^2}}^2+C\left\|\left(u^n, \theta^n\right)\right\|_{H^2}\|f^{n+1}\|_{{H_{x,v}^2}}\nonumber  \\
&\qquad+C\|u^n\|_{H^2}\|f^{n+1}\|_{{H_{x,v}^2}}^2+C\|(u^n, \theta^n)\|_{H^2} \sum_{|\alpha| \leq 2}\|\partial^\alpha f^{n+1}\|_\nu^2, 
\end{align}
holds for some $\lambda_9>0$.
Furthermore, for all $0 \leq t \leq T \leq T^*$,
\begin{align}\label{G2.44}
& \frac{1}{2} \frac{\mathrm{d}}{\mathrm{d} t}\left\|f^{n+1}\right\|_{H_{x, v}^2}^2+\lambda_{10} \sum_{|\alpha|+|\beta| \leq 2}\|\partial_\beta^\alpha f^{n+1}\|_\nu^2 \nonumber \\
& \quad \leq C\|f^{n+1}\|_{H_{x, v}^2}^2+C\|(u^n, \theta^n)\|_{H^2}^2+C\|u^n\|_{H^2}\|f^{n+1}\|_{H^2}^2 \nonumber \\
& \qquad+C\|(u^n, \theta^n)\|_{H^2} \sum_{|\alpha|+|\beta| \leq 2}\|\partial_\beta^\alpha f^{n+1}\|_\nu^2, 
\end{align}
holds for some $\lambda_{10}>0$.
By the classical   theory developed in \cite{MN-jmku-1980}, we known  that the system 
\eqref{G2.39}--\eqref{G2.41} have a unique solution 
$(\rho^{n+1}, u^{n+1}, \theta^{n+1})$ 
with $\rho^{n+1} \geq \tilde \varrho>0$, and
\begin{gather*}
 \rho^{n+1} \in  C\left([0, T], H^2\right) \cap C^1\left([0, T], H^1\right), \\
 u^{n+1}, \theta^{n+1} \in  C\left([0, T], H^2\right) \cap C^1\left([0, T], L^2\right) .
\end{gather*}
In order to estimate $\frac{\mathrm{d}}{\mathrm{d} t}\left\|\left(\rho^{n+1}, u^{n+1}, \theta^{n+1}\right)\right\|_{H^2}^2$,
we apply $\partial^\alpha(|\alpha| \leq 2)$ to the equations \eqref{G2.38}-\eqref{G2.41} and multiply by $\partial^\alpha \rho^{n+1}, \partial^\alpha u^{n+1}, \partial^\alpha \theta^{n+1}$ respectively,
  and then take integration and summation to obtain
\begin{align}\label{G2.45}
&\frac{1}{2} \frac{\mathrm{d}}{\mathrm{d} t}\|(\rho^{n+1}, u^{n+1}, \theta^{n+1})\|_{H^2}^2+\lambda_{11} \sum_{|\alpha| \leq 2} \int_{\mathbb{R}^3}\big(|\nabla_x \partial^\alpha u^{n+1}|^2+|\nabla_x \partial^\alpha u^{n+1}|^2\big) \mathrm{d} x \nonumber \\
&\quad \leq  C(1+\|\rho^n\|_{H^2}^2+\|u^n\|_{H^2}^2)\|\rho^{n+1}\|_{H^2}^2+C(1+\|(\rho,u,\theta)\|_{H^2}^2)\|(u^{n+1}, \theta^{n+1})\|_{H^2}^2\nonumber  \\
&\qquad+C(1+\|\theta^n\|_{H^2}^2)\|\rho^n\|_{H^2}^2
+C(1+\|u^n\|_{H^2}^2)\|u^n\|_{H^2}^2+C(1+\|u^n\|_{H^2}^2)\|\theta^n\|_{H^2}^2\nonumber \\
&\qquad+C(1+\|\theta^n\|_{H^2}^2+\|u^n\|_{H^2}^2+\|u^n\|_{H^2}^4)\|f^n\|_{H^2}^2+C\|f^{n+1}\|_{{H_{x,v}^2}}^2,
\end{align}
for some $\lambda_{11}>0$.
Combining \eqref{G2.44} with \eqref{G2.45}, we obtain
\begin{align}\label{G2.46}
& \frac{1}{2} \frac{\mathrm{d}}{\mathrm{d} t}\left(\|f^{n+1}\|_{H_{x, v}^2}^2+\|\rho^{n+1}\|_{H^2}^2+\|u^{n+1}\|_{H^2}^2+\|\theta^{n+1}\|_{H^2}^2\right)\nonumber  \\
& \qquad+\lambda_{12} \sum_{|\alpha|+|\beta| \leq 2}\|\partial_\beta^\alpha f^{n+1}\|_\nu^2+\lambda \sum_{|\alpha| \leq 2}\left(\|\nabla_x \partial^\alpha u^{n+1}\|^2+\|\nabla_x\partial^\alpha \theta^{n+1}\|^2\right) \nonumber \\
 &\quad\leq C\left(1+\|u^n\|_{H^2}^2\right)\|f^{n+1}\|_{H_{x, v}^2}^2+C\left(1+\|\rho^n\|_{H^2}^2+\|u^n\|_{H^2}^2\right)\|\rho^{n+1}\|_{H^2}^2 \nonumber \\
& \qquad+C\left(1+\|(\rho,u,\theta)\|_{H^2}^2\right)\|(u^{n+1}, \theta^{n+1})\|_{H^2}^2+C\left(1+\|\theta^n\|_{H^2}^2+\|u^n\|_{H^2}^2+\|u^n\|_{H^2}^4\right)\|f^n\|_{H^2}^2\nonumber \\
&\qquad+C\left(1+\|\theta^n\|_{H^2}^2\right)\|\rho^n\|_{H^2}^2
+C\left(1+\|u^n\|_{H^2}^2\right)\|u^n\|_{H^2}^2+C\left(1+\|u^n\|_{H^2}^2\right)\|\theta^n\|_{H^2}^2\nonumber \\
&\qquad+C\left(1+\|\theta^n\|_{H^2}^2+\|u^n\|_{H^2}^2+\|u^n\|_{H^2}^4\right)\|f^n\|_{H^2}^2+C\|(u^n, \theta^n)\|_{H^2} \sum_{|\alpha|+|\beta| \leq 2}\|\partial_\beta^\alpha f^{n+1}\|_\nu^2,
\end{align}
for some $\lambda_{12}>0$.

Next, we define   $B_n(T)$ as
\begin{align*}
B_n(T):=\sup _{0 \leq t \leq T}\left\{\|\rho^n(t)\|_{H^2}^2+\|u^n(t)\|_{H^2}^2+\|\theta^n(t)\|_{H^2}^2+\|f^n(t)\|_{H_{x, v}^2}^2\right\} ,
\end{align*}
and suppose that  $B_n(T) \leq B_0$ and $B_n(0) \leq  {B_0}/{3}$ for some $B_0>0$. 
Then,  integrating \eqref{G2.46} over $[0, T]$ implies that  there exists a positive constant $\lambda_{13}>0$ such that
\begin{align}\label{G2.47}
&B_{n+1}(T)+\lambda_{13} \int_0^T\bigg(\sum_{|\alpha| \leq 2}\left(\left\|\nabla \partial^\alpha u^{n+1}\right\|^2+\|\nabla \partial^\alpha \theta^{n+1}\|^2\right)+\sum_{|\alpha|+|\beta| \leq 2}\|\partial_\beta^\alpha f^{n+1}\|_\nu^2\bigg) \mathrm{d} t \nonumber\\
\leq\,& B_{n+1}(0)+C\left(1+B_n(T)\right) B_{n+1}(T) T+C \left(1+B_n^2(T)\right)B_n(T) T\nonumber \\
&+C B_n^{\frac{1}{2}}(T) \int_0^T \sum_{|\alpha|+|\beta| \leq 2}\|\partial_\beta^\alpha f^{n+1}\|_\nu^2 \mathrm{d} t\nonumber \\
\leq\, & \frac{B_0}{3}+C\left(1+B_0\right) T B_{n+1}(T)+C\left(B_0+B_0^3\right) T +C B_0^{\frac{1}{2}} \sum_{|\alpha|+|\beta| \leq 2} \int_0^T\|\partial_\beta^\alpha f^{n+1}\|_\nu^2 \mathrm{d} t .
\end{align}
Notice that for $t \leq T^*$, it holds
\begin{align*}
\big(1-C\left(1+B_0\right) T\big) B_{n+1}(T) \leq \frac{B_0}{3}+C\left(B_0+B_0^3\right) T .
\end{align*}
By selecting $T^*$ with $CT^*(1+B_0^2) \leq {1}/{3}$, where $B_0$ is sufficient small, one has
$
B_{n+1} \leq B_0 .
$
With these preparations, we    prove that $\|f^{n+1}\|_{H_{x, v}^2}^2$ and 
 $\|(\rho^{n+1},u^{n+1},\theta^{n+1})\|_{H^2}$ are  continuous in $t$. Taking similar arguments to 
  that using in \eqref{G2.44}, we obtain that
\begin{align}\label{G2.48}
& \quad\left|\left\|f^{n+1}(t)\right\|_{H_{x, v}^2}^2-\left\|f^{n+1}(s)\right\|_{H_{x, v}^2}^2\right| \nonumber\\
&=\left|\int_s^t \frac{\mathrm{d}}{\mathrm{d} \eta}\left\|f^{n+1}(\eta)\right\|_{H_{x, v}^2}^2 \mathrm{d} \eta\right| \nonumber\\
& \leq C B_0^{\frac{1}{2}} \sum_{|\alpha|+|\beta| \leq 2} \int_s^t\|\partial_\beta^\alpha f^{n+1}\|_v^2 \mathrm{d} \eta+C(B_0+B_0^{\frac{3}{2}})|t-s| .
\end{align}
In addition, $\|\partial_\beta^\alpha f^{n+1}\|_\nu^2$ is integrable over $\left[0, T^*\right]$. Besides, we also have
\begin{align*}
& \quad\left|\left\|(\rho^{n+1},u^{n+1},\theta^{n+1}) (t)\right\|_{H^2}^2-\left\|(\rho^{n+1},u^{n+1},\theta^{n+1})(s)\right\|_{H^2}^2\right| \\
&=\left|\int_s^t \frac{\mathrm{d}}{\mathrm{d} \eta}\left\|(\rho^{n+1},u^{n+1},\theta^{n+1})(\eta)\right\|_{H^2}^2 \mathrm{d} \eta\right| \\
& \leq C(B_0+B_0^{2}+B_0^{3})|t-s| .
\end{align*}
Therefore, 
$$
\left(f^n, \rho^n, u^n, \theta^n\right) \in \mathcal{X}\left(0, T^* ; B_0\right),
$$ holds for all $n \geq 0$.

Now we prove $\{\left(f^n, \rho^n, u^n, \theta^n\right)\}_{n\geq 0}$ is a Cauchy sequence.
Through direct calculation,   we know that $(f^{n+1}-f^n, $ $\rho^{n+1}-\rho^n, u^{n+1}-u^n, \theta^{n+1}-\theta^n )$ satisfies:
\begin{align*}
&\partial_t(\left.f^{n+1}-f^n\right)+v \cdot \nabla_x\left(f^{n+1}-f^n\right)-\mathbf{L}\left(f^{n+1}-f^n\right) \\
&\quad =-u^n \cdot\left(\nabla_v-\frac{v}{2}\right)\left(f^{n+1}-f^n\right)-\left(u^n-u^{n-1}\right) \cdot\left(\nabla_v-\frac{v}{2}\right) f^n+\sqrt{M} v \cdot\left(u^n-u^{n-1}\right) \\
&\qquad+\frac{\theta^n}{\sqrt{M}}  \Delta_v\left(\sqrt{M}\left(f^{n+1}-f^n\right)\right)+\frac{\Delta_v M\left(\theta^n-\theta^{n-1}\right)}{\sqrt{M}}
+\frac{\left(\theta^n-\theta^{n-1}\right)\Delta_v(\sqrt{M}f^n)}{\sqrt{M}} , \\
&\partial_t(\left.\rho^{n+1}-\rho^n\right)+u^n \cdot \nabla_x\left(\rho^{n+1}-\rho^n\right) \\
&\quad =-\left(u^n-u^{n-1}\right) \cdot \nabla_x \rho^n-\left(1+\rho^n\right) \operatorname{div}\left(u^{n+1}-u^n\right)-\left(\rho^n-\rho^{n-1}\right) \operatorname{div} u^n, \\
&\partial_t\left(u^{n+1}-u^n\right)-\frac{\mu_1}{1+\rho^n} \Delta_x\left(u^{n+1}-u^n\right)-\frac{\mu_1+\mu_2}{1+\rho^n} \nabla_x\left({\rm div}(u^{n+1}-u^n)\right) \\
&\quad=\left(\frac{1}{1+\rho^n}-\frac{1}{1+\rho^{n-1}}\right)\left(\mu_1\Delta_x u^n+(\mu_1+\mu_2)\nabla_x({\rm div}_x  u^n)+b^{n}-u^n-u^n a^{n-1}\right) \\
&\quad\quad-\left(\left(u^n-u^{n-1}\right) \cdot \nabla_x u^n+u^{n} \cdot \nabla_x\left(u^{n+1}-u^{n}\right)\right)-\nabla_x\left(\theta^{n+1}-\theta^{n}\right) \\
&\quad\quad+\left(\frac{1+\theta^n}{1+\rho^n} \nabla_x\left(\rho^{n+1}-\rho^{n}\right)-\left(\frac{1+\theta^n}{1+\rho^n}-\frac{1+\theta^{n-1}}{1+\rho^{n-1}}\right) \nabla_x \rho^{n}\right) \\
&\quad\quad+\frac{1}{1+\rho^n}\left(\left(b^{n+1}-b^{n}\right)-\left(u^{n+1}-u^{n}\right)-u^n (a^n-a^{n-1})-a^n(u^{n+1}-u^n)\right), \\
&\partial_t\left(\theta^{n+1}-\theta^n\right)-\frac{\kappa}{1+\rho^n} \Delta_x\left(\theta^{n+1}-\theta^n\right) \\
&\quad=- \left(u^n-u^{n-1}\right) \cdot \nabla_x \theta^n-u^{n} \cdot \nabla_x\left(\theta^{n+1}-\theta^{n}\right) -\operatorname{div}\left(u^{n+1}-u^{n}\right) \\
&\quad\quad-\left(\theta^{n+1}-\theta^{n}\right) \operatorname{div} u^n-\operatorname{div}\left(u^n-u^{n-1}\right) \theta^n \\
&\quad\quad+\left(\frac{1}{1+\rho^n}-\frac{1}{1+\rho^{n-1}}\right)\left(\kappa\Delta_x \theta^n-\sqrt{6} \omega^n+3 \theta^n-3 \theta^n a^n\right) \\
&\quad\quad+\left(\frac{1}{1+\rho^n}-\frac{1}{1+\rho^{n-1}}\right)\left(\left|u^n\right|^2-2 u^n \cdot b^n+a^n\left|u^n\right|^2+2\mu_1|D(u^n)|^2+\mu_2|{\rm div}_x  u^n|^2\right) \\
&\quad\quad+\frac{1}{1+\rho^{n-1}}\left(\left|u^n\right|^2-\left|u^{n-1}\right|^2-2 u^n \cdot b^n+2 u^{n-1} \cdot b^{n-1}+a^n\left|u^n\right|^2-a^{n-1}\left|u^{n-1}\right|^2\right) \\
&\quad\quad-\frac{1}{1+\rho^{n-1}}\Big(3\left(a^n \theta^n-a^{n-1} \theta^{n-1}\right)-2\mu_1\big (|D(u^n)|^2- |D(u^{n-1})|^2\big)\Big) \\
&\quad\quad-\frac{1}{1+\rho^{n}}\Big(\sqrt{6}\left(\omega^{n+1}-\omega^{n}\right)-3\left(\theta^{n+1}-\theta^{n}\right)\Big)+\frac{\mu_2}{1+\rho^{n-1}}\big (|{\rm div}_x  u^n|^2- |{\rm div}_x  u^{n-1}|^2\big) .
\end{align*}
Taking similar arguments to that using in \eqref{G2.46}, we obtain that there exists a positive constant $\lambda_{14}>0$ such that
\begin{align*}
&\frac{1}{2} \frac{\mathrm{d}}{\mathrm{d} t} \Big(\|f^{n+1}-f^n\|_{H_{x, v}^1}^2+\|\rho^{n+1}-\rho^n\|_{H^1}^2+\|u^{n+1}-u^n\|_{H^1}^2+\|\theta^{n+1}-\theta^n\|_{H^1}^2\Big) \\
&\qquad  +\lambda_{14}\bigg(\sum_{|\alpha|+|\beta| \leq 1}\|\partial_\beta^\alpha(f^{n+1}-f^n)\|_\nu^2+\lambda\|\nabla_x(u^{n+1}-u^n), \nabla_x(\theta^{n+1}-\theta^n)\|_{H^1}^2\bigg) \\
\,&\quad \leq \tilde{C}\left(\|f^{n+1}-f^n\|_{H_{x, v}^1}^2+\|\rho^{n+1}-\rho^n\|_{H^1}^2+\|u^{n+1}-u^n\|_{H^1}^2+\|\theta^{n+1}-\theta^n\|_{H^1}^2\right) \\
& \qquad +\tilde{C}\left(\|f^n-f^{n-1}\|_{H_{x, v}^1}^2+\|\rho^n-\rho^{n-1}\|_{H^1}^2+\|u^n-u^{n-1}\|_{H^1}^2+\|\theta^n-\theta^{n-1}\|_{H^1}^2\right),
\end{align*}
where $\tilde{C}>0$ is a constant   depending on $\|(\rho^n, \rho^{n-1})\|_{H^2},\|(u^n, u^{n-1})\|_{H^2},\|(\theta^n, \theta^{n-1})\|_{H^2}$, $\|(f^n, f^{n-1})\|_{H_{x, v}^2}$ and $\sum_{|\alpha|+|\beta| \leq 2}\|\partial_\beta^\alpha f^n\|_\nu^2$. From \eqref{G2.47}, we    infer that
\begin{align*}
\int_0^{T^*} \sum_{|\alpha|+|\beta| \leq 2}\|\partial_\beta^\alpha f^n\|_\nu^2 \mathrm{d} s,
\end{align*}
is sufficient small  for all $n\geq 0$. Therefore, there exists a constant  $\eta_1\in(0,1)$ such that
\begin{align}\label{G2.49}
& \sup _{0<t \leq T^*}\big\{\|f^{n+1}-f^n\|_{H^1}
+\|\rho^{n+1}-\rho^n\|_{H^1}+\|u^{n+1}-u^n\|_{H^1}+\|\theta^{n+1}-\theta^n\|_{H^1}\big\}\nonumber\\
& \quad \leq \eta_1 \sup _{0<t \leq T^*}\big\{\|f^n-f^{n-1}\|_{H^1}+\|\rho^n-\rho^{n-1}\|_{H^1}+\|u^n-u^{n-1}\|_{H^1}+\|\theta^n-u^{n-1}\|_{H^1}\big\}.
\end{align}

According to \eqref{G2.49}, we deduce that $\left\{(f^n, \rho^n, u^n, \theta^n)\right\}_{n \geq 0}$ is a Cauchy sequence in the Banach space $C\left(\left[0, T^*\right], H_{x,v}^1\right) \times\left(C\left(\left[0, T^*\right], H^1\right)\right)\times\left(C\left(\left[0, T^*\right], H^1\right)\right)\times\left(C\left(\left[0, T^*\right], H^1\right)\right)$. Due to   $(f, \rho, u, \theta)$ is the limit of $(f^n, \rho^n, u^n, \theta^n)$, 
it is readily to verify that   $(f, \rho, u, \theta)$ satisfies the problem
 \eqref{G1}--\eqref{G5} through letting $n \rightarrow \infty$. 

Using the similar argument to \eqref{G2.48} implies that $f \in C\left(\left[0, T^*\right], H^1\right)$. 
Furthermore, we    infer that $(f, \rho, u, \theta) \in \mathcal{X}\left(0, T^* ; B_0\right)$.

To prove the uniqueness of strong solutions, we  assume that  $(\bar{f}, \bar{\rho}, \bar{u}, \bar{\theta}) \in \mathcal{X}\left(0, T^*, B_0\right)$ is an another solution to the Cauchy problem \eqref{G1}--\eqref{G5}. Similar to \eqref{G2.49}, it holds that
\begin{align*}
& \sup _{0<t \leq T^*}\big\{\|f-\bar{f}\|_{H^1}+\|\rho-\bar{\rho}\|_{H^1}+\|u-\bar{u}\|_{H^1}+\|\theta-\bar{\theta}\|_{H^1}\big\} \\
& \quad \leq \eta_2 \sup _{0<t \leq T^*}\big\{\|f-\bar{f}\|_{H^1}+\|\rho-\bar{\rho}\|_{H^1}+\|u-\bar{u}\|_{H^1}+\|\theta-\bar{\theta}\|_{H^1}\big\},
\end{align*}
for some  $0<\eta_2<1$, which implies that   $f \equiv \bar{f}, \rho \equiv \bar{\rho}, u \equiv \bar{u}, \theta \equiv \bar{\theta}$.
Thus, we complete the proof for the local existence and uniqueness of strong solutions
to the problem \eqref{G1}--\eqref{G5}.
\end{proof}

\subsection{Global  existence of strong solutions} 
By Theorem \ref{T2.7}, 
we know that the problem \eqref{G1}-\eqref{G5} has a unique  local strong solution.
According to the uniform a priori estimates \eqref{G2.37}, 
the local existence result, and continuity argument, 
we    get the global existence of strong solutions to 
the  problem  \eqref{G1}--\eqref{G5}. 
Finally, by the maximum principle, we    conclude that
$F=M+\sqrt{M}f\geq0$. Therefore, 
we finish the proof of the existence  part of Theorem \ref{T1.1}.

\section{Time-decay rates of the strong solutions and their gradients} 
In this section, we deal with the time decay rates of the global strong solution $(\rho, u, \theta, f)$
to the  problem \eqref{G1}--\eqref{G5} in $\mathbb{R}^{3}$. To archive this goal,
we first consider the linearized system corresponding to \eqref{G1}--\eqref{G4} as follows
\begin{equation}\label{C1}
\left\{\begin{aligned}
&\partial_t f+v \cdot \nabla_x f-u \cdot v\sqrt{M}-\theta\left(|v|^2-3\right) \sqrt{M}-\mathbf{L} f=\mathfrak{S}, \\
&\partial_t \rho+{\rm div}_x  u=0, \\
&\partial_t u-\mu_1\Delta_x u-(\mu_1+\mu_2)\nabla_x({\rm div}_x  u)+\nabla_x \theta+\nabla_x \rho+u-b=0, \\
&\partial_t \theta-\kappa\Delta_x \theta+{\rm div}_x  u+\sqrt{3}(\sqrt{3} \theta-\sqrt{2} \omega)=0 .\\
\end{aligned}
\right.\end{equation}
Here $ \mathfrak{S}$     takes the form of
\begin{equation*}
\mathfrak{S}=\mathrm{div}_v G-\frac{1}{2}v\cdot G+h,
\end{equation*}
with 
\begin{equation*}
\ \ G=\big(G_1(t,x,v), G_2(t,x,v), G_3(t,x,v)\big)\in\mathbb{R}^{3},\ \ h=h(t,x,v)\in\mathbb{R},
\end{equation*}
satisfying
\begin{equation*}
\mathbf{P}_{0}G=0,\quad\mathbf{P}_{1}G=0,\quad \mathbf{P}h=0.
\end{equation*}
We supplement the system \eqref{C1} with the initial data 
\begin{equation}\label{C2}
(\rho, u, \theta, f)|_{t=0}=\big(\rho(0,x) ,u(0,x), \theta(0,x), f(0,x,v)\big), \quad (x, v)\in\mathbb{R}^{3}_x\times \mathbb{R}^3_v.
\end{equation}

Now we denote  $U=(\rho ,u ,\theta ,f )$ and $U_0=(\rho_0,u_0,\theta_0,f_0)$ to simplify 
the presentation of the  problem \eqref{C1}--\eqref{C2}.
By Duhamel's principle, we have
\begin{equation}\label{NJK3.3}
U(t)=\mathbb{A}(t)U_{0}+\int_{0}^{t}\mathbb{A}(t-s)(\mathfrak{S}(s),0,0,0){\rm d}s,
\end{equation}
for all $t\geq 0$, where $\mathbb{A}(t)$ is  the solution operator to 
\eqref{C1}--\eqref{C2} with $\mathfrak{S}=0$. 
Below we establish the following time-decay property of $U(t)$.
\begin{thm}\label{C3}
Let $1\leq q\leq2$. For any $\alpha$, $\alpha^{\prime}$ with 
$\alpha^{\prime}\leq\alpha$ and $m=|\alpha-\alpha^{\prime}|$, we have 
\begin{align}
& \|\partial^{\alpha}\mathbb{A}(t)U_{0}\|_{\mathcal{Z}_{2}}
\leq C\left(1+t\right)^{-\frac{3}{2}\left(\frac{1}{q}-\frac{1}{2}\right)-\frac{m}{2}}
\left(\|\partial^{\alpha^{\prime}}U_{0}\|_{\mathcal{Z}_{q}}+\|\partial^{\alpha}U_{0}\|_{\mathcal{Z}_{2}}\right),
\,\, 
t\geq 0,
\label{T2}\\ 
&\left\|\partial^{\alpha}\int_{0}^{t}\mathbb{A}(t-s)
\big(\mathfrak{S}(s),0,0,0\big){\rm d}s\right\|_{\mathcal{Z}_{2}}^{2} \leq\,C\int_0^t(1+t-s)^{-3\left(\frac{1}{q}-\frac{1}{2}\right)-m}\nonumber\\
& \qquad \quad
\times \left(\left\|\partial^{\alpha^{\prime}}\big(G(s),\nu^{-1/2}h(s)\big)\right\|_{\mathcal{Z}_{q}}^{2}
+\left\|\partial^{\alpha}\big(G(s),\nu^{-1/2}h(s)\big)\right\|_{\mathcal{Z}_{2}}^{2}\right){\rm d}s, \,\, 
t\geq 0. \label{C4}
\end{align}
\end{thm}
\begin{proof}
Applying  Fourier transform in $x$ to \eqref{C1}, we obtain that 
\begin{equation}\label{C5}\left\{
\begin{aligned}
&\partial_t \hat{f}+i v \cdot k \hat{f}-\hat{u} \cdot v \sqrt{M}-\hat{\theta}\left(|v|^2-3\right) \sqrt{M}
=\mathbf{L} \hat{f}+{\rm div}_v  \hat{G}-\frac{1}{2} v \cdot \hat{G}+\hat{h}, \\
&\partial_t \hat{\rho}+i k \cdot \hat{u}=0, \\
&\partial_t \hat{u}+\mu_1|k|^2 \hat{u}+(\mu_1+\mu_2)k(k\cdot\hat{u})+i k \hat{\theta}+i k \hat{\rho}+\hat{u}-\hat{b}=0, \\
&\partial_t \hat{\theta}+\kappa|k|^2 \hat{\theta}+i k \hat{u}+\sqrt{3}(\sqrt{3} \hat{\theta}-\sqrt{2} \hat{\omega})=0 .
\end{aligned}  \right.
\end{equation}
By taking the inner product of the first equation
in \eqref{C5} with  $\Bar{\hat{f}}$ and integrating the result over $\mathbb{R}^3$ in $v$, we have
\begin{align*}
& \frac{1}{2} \partial_t\|\hat{f}\|_{L_v^2}^2+\mathfrak{Re} 
\int_{\mathbb{R}^3}\langle-\mathbf{L}\{\mathbf{I}-\mathbf{P}\} \hat{f} \mid\{\mathbf{I}-\mathbf{P}\} \hat{f}\rangle \mathrm{d} v
+|\hat{b}|^2+2|\hat{\omega}|^2-\mathfrak{Re}(\hat{u} \mid \hat{b})-\mathfrak{Re}(\hat{\theta} \mid \sqrt{6} \hat{\omega})\\
& \quad=\mathfrak{Re} \int_{\mathbb{R}^3}\left({\rm div}_v   \hat{G}-\frac{1}{2} v \cdot \hat{G} \rvert\{\mathbf{I}-\mathbf{P}\} \hat{f}\right) \mathrm{d} v+\mathfrak{Re} \int_{R^3}(\hat{h} \mid\{\mathbf{I}-\mathbf{P}\} \hat{f}) \mathrm{d} v,
\end{align*}
where we have used the following facts:
\begin{align*}
{\rm div}_v  G-\frac{1}{2} v \cdot G\perp {\rm Range}( \mathbf{P}), \quad \int_{R^3}(\hat{h} \mid \mathbf{P} \hat{f}) \mathrm{d} v=0 .
\end{align*}
Similarly, taking the inner product of remaining equations in \eqref{C1}, we obtain
\begin{equation*} \left\{
\begin{aligned}
&\frac{1}{2} \partial_t|\hat{\rho}|^2+\mathfrak{Re}( i k \hat{u} \mid \hat{\rho})=0, \\
 &\frac{1}{2} \partial_t|\hat{u}|^2+\mathfrak{Re}( i k \hat{\theta} \mid \hat{u})+\mathfrak{Re}( i k \hat{\rho} \mid \hat{u})+|\hat{u}|^2+\mu_1|k|^2|\hat{u}|^2+(\mu_1+\mu_2)|k\cdot\hat{u}|^2-\mathfrak{Re}(\hat{b} \mid \hat{u})=0,\\
& \frac{1}{2} \partial_t|\hat{\theta}|^2+\mathfrak{Re}( i k \hat{u} \mid \hat{\theta})+3|\hat{\theta}|^2+\kappa|k|^2|\hat{\theta}|^2-\mathfrak{Re}(\sqrt{6} \hat{\omega} \mid \hat{\theta})=0.
 \end{aligned}\right.
\end{equation*}

Using the property \eqref{G1.12} and the Cauchy's inequality, we have
\begin{align}\label{G3.6}
& \frac{1}{2} \partial_t\left(\|\hat{f}\|_{L_v^2}^2+|\hat{\rho}|^2
+|\hat{u}|^2+|\hat{\theta}|^2\right)+\lambda_0|\{\mathbf{I}-\mathbf{P}\} \hat{f}|_\nu^2 \nonumber\\
& \qquad +|\hat{u}-\hat{b}|^2+(\mu_1+\mu_2)|k\cdot\hat{u}|^2+|\sqrt{2} \hat{\omega}-\sqrt{3} \hat{\theta}|^2+\mu_1|k|^2|\hat{u}|^2+\kappa
|k|^2|\hat{\theta}|^2\nonumber\\
\,&\quad \leq  C\left(\|\hat{G}\|^2+\|\nu^{-\frac{1}{2}} \hat{h}\|^2\right) .
\end{align}
Now we   deal with the estimates of $a, b,$ and $ \omega$.
Similar to the processing of \eqref{G2.13}--\eqref{G2.17}, one has 
\begin{equation*} \left\{
\begin{aligned}
& \partial_t a+{\rm div}_x  b=0, \\
& \partial_t b_i+\partial_i a+\frac{2}{\sqrt{6}} 
\partial_i \omega+\sum_{j=1}^3 \partial_{x_j} \Gamma_{i, j}(\{\mathbf{I}-\mathbf{P}\} f)-u_i+b_i=0, \\
& \partial_t \omega+\sqrt{2}(\sqrt{2} \omega-\sqrt{3} \theta)
+\frac{2}{\sqrt{6}} {\rm div}_x  b+\sum_{i=1}^3 \partial_{x_i} \Upsilon_i(\{\mathbf{I}-\mathbf{P}\} f)=0, \\
& \partial_j b_i+\partial_i b_j-\frac{2}{\sqrt{6}} \delta_{i j}
\left(\frac{2}{\sqrt{6}} {\rm div}_x  b+\sum_{i=1}^3 \partial_{x_i} 
\Upsilon_i(\{\mathbf{I}-\mathbf{P}\} f)\right) \\
& \quad=-\partial_t \Gamma_{i, j}\{\mathbf{I}-\mathbf{P}\} f
+\Gamma_{i, j}\left(\mathfrak{l}+\mathfrak{S}\right), \\
& \frac{5}{3} \partial_i \omega-\frac{2}{\sqrt{6}} \sum_{j=1}^3 \partial_{x_j} 
\Gamma_{i, j}(\{\mathbf{I}-\mathbf{P}\} f)
=-\partial_t \Upsilon_i(\{\mathbf{I}-\mathbf{P}\} f)+\Upsilon_i\left(\mathfrak{l}+\mathfrak{S}\right),
\end{aligned}\right.
\end{equation*}
where $\mathfrak{l}$ is defined in \eqref{NJK2.181}.
Similarly, by taking  Fourier transform in $x$, we get
\begin{equation*} \left\{
\begin{aligned}
& \partial_t \hat{a}+i k \cdot \hat{b}=0, \\
& \partial_t \hat{b}_i+i k_i \hat{a}+\frac{2}{\sqrt{6}} i k_i \hat{\omega}+\sum_{j=1}^3 i k_j \Gamma_{i, j}(\{\mathbf{I}-\mathbf{P}\} \hat{f})-\hat{u}_i+\hat{b}_i=0, \\
& \partial_t \hat{\omega}+\sqrt{2}(\sqrt{2} \hat{\omega}-\sqrt{3} \hat{\theta})+\frac{2}{\sqrt{6}} i k \cdot \hat{b}+\sum_{i=1}^3 i k_i \Upsilon_i(\{\mathbf{I}-\mathbf{P}\} \hat{f})=0, \\
& i k_j \hat{b}_i+i k_i \hat{b}_j-\frac{2}{\sqrt{6}} \delta_{i j}\left(\frac{2}{\sqrt{6}} i k \cdot \hat{b}+\sum_{i=1}^3 i k_i \Upsilon_i(\{\mathbf{I}-\mathbf{P}\} \hat{f})\right) \\
& \quad=-\partial_t \Gamma_{i, j}\{\mathbf{I}-\mathbf{P}\} \hat{f}+\Gamma_{i, j}\big(\hat{\mathfrak{l}}+\hat{\mathfrak{S}}\big), \\
& \frac{5}{3} i k_i \hat{\omega}-\frac{2}{\sqrt{6}} \sum_{j=1}^3 i k_j 
\Gamma_{i, j}(\{\mathbf{I}-\mathbf{P}\} \hat{f})=-\partial_t \Upsilon_i(\{\mathbf{I}-\mathbf{P}\} \hat{f})
+\Upsilon_i\big(\hat{\mathfrak{l}}+\hat{\mathfrak{S}}\big) .
\end{aligned}\right.
\end{equation*}
By direct calculation, we   
 obtain there exists a positive constant $\lambda_{15}>0$, such that
\begin{align}
& \partial_t \mathfrak{Re} \sum_{i, j}\Big( i k_i \hat{b}_j+i k_j \hat{b_i} 
\mid \Gamma_{i, j}(\{\mathbf{I}-\mathbf{P}\} \hat{f})\Big)
+\lambda_{15}\left(|k|^2|\hat{b}|^2+|k \cdot \hat{b}|^2\right) \label{G3.71}\\
& \quad \leq \varepsilon|k|^2\left(|\hat{a}|^2+|\hat{\omega}|^2\right)
+C\left(1+|k|^2\right)|\{\mathbf{I}-\mathbf{P}\} \hat{f}|_{L_v^2}^2+C\left(|\hat{u}-\hat{b}|^2+\|\hat{G}\|^2+\|\nu^{-\frac{1}{2}} \hat{h}\|^2\right),\nonumber \\
& \partial_t \mathfrak{Re} \sum_i\left( i k_i \hat{\omega} 
\mid \Upsilon_i(\{\mathbf{I}-\mathbf{P}\} \hat{f})\right)+|k|^2|\hat{\omega}|^2 \label{G3.72}\\
& \quad \leq \varepsilon|k|^2|\hat{b}|^2+C\left(1+|k|^2\right)|\{\mathbf{I}-\mathbf{P}\} \hat{f}|_{L_v^2}^2+C\left(|\sqrt{2} \hat{\omega}-\sqrt{3} \hat{\theta}|^2+\|\hat{G}\|^2+\|\nu^{-\frac{1}{2}} \hat{h}\|^2\right), \nonumber\\
& \partial_t \mathfrak{Re}\bigg(\hat{a} \bigg\lvert\, i 
\frac{\sqrt{6}}{5} \sum_{j=1}^3 k_j \Upsilon_j(\{\mathbf{I}-\mathbf{P}\} 
\hat{f})-i k \cdot \hat{b}\bigg)+\frac{3}{4}|k|^2|\hat{a}|^2 \label{G3.73}\\
& \quad \leq \frac{5}{4}|k \cdot \hat{b}|^2+C\left(1+|k|^2\right)|\{\mathbf{I}-\mathbf{P}\} \hat{f}|_{L_v^2}^2+C\left(|\hat{u}-\hat{b}|^2+\|\hat{G}\|^2+\|\nu^{-\frac{1}{2}} \hat{h}\|^2\right), \nonumber\\
& \partial_t \mathfrak{Re}(\hat{u} \mid i k \hat{\rho})
+\frac{3}{4}|k|^2|\hat{\rho}|^2 \leq C|k|^2\left(|\hat{\theta}|^2
+|\hat{u}|^2\right)+C|\hat{u}-\hat{b}|^2+|k\cdot\hat{u}|^2 .\label{G3.7}
\end{align}
Define 
\begin{align}
\mathcal{E}_1(\hat{f}):= & \frac{1}{1+|k|^2}\bigg\{\sum_{i, j}
\left(i k_i \hat{b}_j+i k_j \hat{b}_i \mid \Gamma_{i, j}(\{\mathbf{I}-\mathbf{P}\} \hat{f})\right)
\nonumber\\
& +\sum_i\left( i k_i \hat{\omega} \mid \Upsilon_i(\{\mathbf{I}-\mathbf{P}\} \hat{f})\right)
+\kappa_1\bigg(\hat{a} \bigg\lvert\, i \frac{\sqrt{6}}{5} \sum_{j=1}^3 k_j \Upsilon_j(\{\mathbf{I}-\mathbf{P}\} 
\hat{f})-i k \cdot \hat{b} \bigg)\bigg\},
\end{align}
where $0<\kappa_1\ll 1$ is a constant.
We know that  there exists a positive constant $\lambda_{16}>0$, such that
\begin{align}\label{G3.9}
& \partial_t \mathfrak{Re} \mathcal{E}_1(\hat{f})+\frac{\lambda_{16}|k|^2}{1+|k|^2}\left(|\hat{a}|^2+|\hat{b}|^2+|\hat{\omega}|^2\right)\nonumber \\
& \quad \leq C\left(|\{\mathbf{I}-\mathbf{P}\} \hat{f}|_{L_v^2}^2+|\hat{u}-\hat{b}|^2+|\sqrt{2} \hat{\omega}-\sqrt{3} \hat{\theta}|^2\right)+C\left(\|\hat{G}\|^2+\|\nu^{-\frac{1}{2}} \hat{h}\|^2\right) .
\end{align}
Define 
\begin{align*}
{\mathcal E}_{M}\big(\hat{U}(t,k)\big)
:=\|\hat{f}\|_{L_v^2}^2+|\hat{\rho}|^2+|\hat{u}|^2
+|\hat{\theta}|^2+\kappa_2\mathfrak{Re}\mathcal{E}_1(\hat{f})+\kappa_3 \mathfrak{Re}\frac{1}{1+|k
|^2}(\hat{u} \mid i k\hat{\rho}),
\end{align*}
where $\kappa_2$ and $\kappa_3$ are constants that satisfy $0<\kappa_2, \kappa_3\ll 1$. 
According to  \eqref{G3.6}, \eqref{G3.7} and \eqref{G3.9}, we have
\begin{equation*}
\partial_{t}\mathcal{E}_{M}\big(\hat{U}(t,k)\big)+\frac{\lambda_{17}|k|^{2}}{1+|k|^{2}}
\mathcal{E}_{M}(\hat{U}(t,k))\leq C\left(\|\hat{G}\|_{L_{v}^{2}}^{2}
+\|\nu^{-\frac{1}{2}}\hat{h}\|_{L_{v}^{2}}^{2}\right),
\end{equation*}
for some $\lambda_{17}>0$.
Using Gronwall's inequality, we obtain 
\begin{equation*}
\mathcal{E}_{M}(\hat{U}(t,k))\leq Ce^{\frac{-\lambda_{17}|k|^{2}t}{1+|k|^{2}}}
\mathcal{E}_{M}(\hat{U}(0,k))+C\int_{0}^{t}e^{\frac{-\lambda_{17}|k|^{2}(t-s)}{1+|k|^{2}}}
\left(\|\hat{G}\|_{L_{v}^{2}}^{2}+\|\delta^{-\frac{1}{2}}\hat{h}\|_{L_{v}^{2}}^{2}\right)ds.
\end{equation*}
As mentioned in \cite{DFT-cmp-2010, KS-thesis-1984}, the expected time-decay properties 
\eqref{T2} and  \eqref{C4} can be obtained directly via the above estimate. 
For brevity, below we only prove \eqref{T2}.

Now we set $\mathfrak{S}(s)=0$. Therefore,
\begin{align*}
\|\partial^{\alpha}\mathbb{A}(t)U_{0}\|_{\mathcal{Z}_{2}}^2\leq &\,C \int_{\mathbb{R}^3} |k|^{2\alpha}e^{\frac{-\lambda_{17}|k|^{2}t}{1+|k|^{2}}}
\mathcal{E}_{M}(\hat{U}(0,k))\mathrm{d}k\\
\leq&\, C \int_{|k|\leq 1} |k|^{2\alpha}e^{\frac{-\lambda_{17}|k|^{2}t}{1+|k|^{2}}}
\mathcal{E}_{M}(\hat{U}(0,k))\mathrm{d}k+\int_{|k|\geq 1} |k|^{2\alpha}e^{\frac{-\lambda_{17}|k|^{2}t}{1+|k|^{2}}}
\mathcal{E}_{M}(\hat{U}(0,k))\mathrm{d}k\\
\leq&\, C \left(1+t\right)^{-3\left(\frac{1}{q}-\frac{1}{2}\right)-m}
\|\partial^{\alpha^{\prime}}U_{0}\|_{\mathcal{Z}_{q}}^2+Ce^{\frac{-\lambda_{17}t}{2}}\|\partial^{\alpha}U_{0}\|_{\mathcal{Z}_{2}}^2,
\end{align*}
since
\begin{align*}
\int_{|k|\leq 1} |k|^{2\alpha}e^{\frac{-\lambda_{17}|k|^{2}t}{1+|k|^{2}}}
\mathcal{E}_{M}(\hat{U}(0,k))\mathrm{d}k\leq&\, C \int_{|k|\leq 1} |k|^{2|\alpha-\alpha^{\prime}|}e^{\frac{-\lambda_{17}|k|^{2}t}{1+|k|^{2}}}|k|^{2\alpha^{\prime}}
\mathcal{E}_{M}(\hat{U}(0,k))\mathrm{d}k\\
\leq &\, C \Big\||k|^{m} 
e^{\frac{-\lambda_{17}|k|^{2}t}{2+2|k|^{2}}}
\Big\|_{L^{\frac{2q}{2-q}}}^2\|\partial^{\alpha}U_{0}\|_{\mathcal{Z}_{q}}^2.
\end{align*}
Here we have used Hausdorff-Young and H\"{o}lder's inequalities.
Thanks to change of variables and the property of Gamma function, we have
\begin{align*}
\Big\||k|^{m} e^{\frac{-\lambda_{17}|k|^{2}t}{2+2|k|^{2}}}\Big\|_{L^{\frac{2q}{2-q}}}^{\frac{2q}{2-q}}\leq&\, C\int_{0}^{\infty} |k|^{m(\frac{2q}{2-q})}e^{-\lambda_{17}|k|^{2}t(\frac{2q}{2-q})}\mathrm{d}k\\
\leq &\, C\int_{0}^{\infty} |x|^\frac{m(\frac{2q}{2-q})+1}{2}e^{-\lambda_{17}tx(\frac{2q}{2-q})}\mathrm{d}x\\
\leq & \,C(1+t)^{\frac{-m(\frac{2q}{2-q})-3}{2}}\int_{0}^{\infty} x^{\frac{m(\frac{2q}{2-q})+3}{2}-1}e^{-x}\mathrm{d}x\\
\leq&\, C \left(1+t\right)^{{\frac{-m(\frac{2q}{2-q})-3}{2}}},
\end{align*}
which implies
\begin{align*}
\Big\||k|^{m} e^{\frac{-\lambda_{17}|k|^{2}t}{2+2|k|^{2}}}\Big\|_{L^\frac{2q}{2-q}}^2
\leq&\, C \left(1+t\right)^{-m-3(\frac{1}{q}-\frac{1}{2})}.  
\end{align*}
\end{proof}

Now we consider the nonlinear problem \eqref{G1}--\eqref{G5}. 
 By Duhamel's principle, the solution $U=(\rho,u,\theta,f)$ to \eqref{G1}--\eqref{G5} can  be rewritten as
\begin{equation}\label{NJK3.11}
U(t)=\mathbb{A}(t)U_{0}+\int_{0}^{t}\mathbb{A}(t-s)(\mathfrak{T}(s),\mathfrak{H}_{1},
\mathfrak{H}_{2},\mathfrak{H}_{3}) {\rm d}s
=:\sum_{i=1}^{4}\mathfrak{K}_{i}(t),
\end{equation}
where   $\mathfrak{T}$  takes the form
\begin{align*}
\mathfrak{T}={\rm div}_v  G-\frac{1}{2} v \cdot G+h+a u \cdot v \sqrt{M}-u \cdot b\sqrt{M}+a \theta\left(|v|^2-3\right) \sqrt{M}+u\cdot v(b\cdot v\sqrt{M})
\end{align*}
with
\begin{align*}
G=-u\left\{\mathbf{I}-\mathbf{P}_{\mathbf{0}}-\mathbf{P}_{\mathbf{1}}\right\} f, 
\quad h=\theta \left\{\mathbf{I}-\mathbf{P}_2\right\} \frac{1}{\sqrt{M}} 
\nabla_v\cdot\left(\sqrt{M}\Big(\nabla_v f-\frac{v}{2} f\Big)\right) .
\end{align*}
And $\mathfrak{H}_i(i=1,2,3)$ reads 
\begin{align*}
\mathfrak{H}_1=&-{\rm div}_x (\rho u), \\
\mathfrak{H}_2=& -u \cdot \nabla_x u+\frac{\rho-\theta}{1+\rho} \nabla_x \rho+\frac{\rho}{1+\rho}(u-b)-\frac{1}{1+\rho} a u-\frac{\rho\mu_1}{1+\rho} \Delta_x u-\frac{(\mu_1+\mu_2)\rho}{1+\rho}\nabla_x({\rm div}_x  u), \\
\mathfrak{H}_3 =&-  u \cdot \nabla_x \theta-\theta {\rm div}_x  u+\frac{\sqrt{3} \rho}{1+\rho}(\sqrt{3} \theta-\sqrt{2} \omega) \\
& +\frac{1}{1+\rho}\left((1+a)|u|^2-3 a \theta-2 u \cdot b-\rho\kappa \Delta_x \theta+2\mu_1|D(u)|^2+\mu_2|{\rm div}u|^2\right).
\end{align*}
Therefore, $\mathfrak{K}_{i}(i=1,2,3,4)$ takes the form
\begin{align*}
&\mathfrak{K}_{1}(t) =\mathbb{A}(t)U_{0},  \\
&\mathfrak{K}_{2}(t) =\int_{0}^{t}\mathbb{A}(t-s)\left({\rm div}_v G-\frac{1}{2}v\cdot G+h,0,0,0\right) {\rm d}s ,  \\
&\mathfrak{K}_{3}(t) =\int_{0}^{t}\mathbb{A}(t-s)\big(au\cdot v\sqrt{M}-u\cdot b\sqrt{M}+u\cdot v(b\cdot v\sqrt{M})+a\theta(|v|^2-3)\sqrt{M},0,0,0\big) {\rm d}s ,  \\
&\mathfrak{K}_{4}(t) =\int_{0}^{t}\mathbb{A}(t-s)(0,\mathfrak{H}_{1},\mathfrak{H}_{2},\mathfrak{H}_{3}) {\rm d}s . 
\end{align*}

From \eqref{G2.36}, one has
\begin{equation}\label{G3.11}
\frac{\mathrm{d}}{\mathrm{d} t} \mathcal{E}(t)+\lambda_{18} \mathcal{D}(t) \leq 0,
\end{equation}
for some $\lambda_{18}>0$.
By the definitions of $\mathcal{E}(t)$ and $\mathcal{D}(t)$, it follows that
\begin{align}\label{G3.12}
\mathcal{E}(t) & \leq C\left(\|\{\mathbf{I}-\mathbf{P}\} f\|_{H^2}^2+\|(a, b, \omega)\|_{H^2}^2+\|(\rho, u, \theta)\|_{H^2}^2\right) \nonumber\\
& \leq C\left(\mathcal{D}(t)+\|f\|_{L^2}^2+\|(\rho, u, \theta)\|_{L^2}^2\right) .
\end{align}
Combining \eqref{G3.11} with \eqref{G3.12}, we get
\begin{equation}
\frac{\mathrm{d}}{\mathrm{d} t} \mathcal{E}(t)+\lambda_{19} \mathcal{E}(t) \leq C\| U\|_{\mathcal{Z}_{2}}^{2},
\end{equation}
for some $\lambda_{19}>0$.
Using the Gronwall's inequality, we obtain
\begin{equation}\label{G3.14}
\mathcal{E}(t) \leq e^{-\lambda_{19} t} \mathcal{E}(0)+C \int_0^t e^{-\lambda_{19}(t-s)}\| U\|_{\mathcal{Z}_{2}}^{2} \mathrm{d} s .
\end{equation}

Applying \eqref{T2} to $\mathfrak{K}_1$, we get
\begin{align*}
\left\|\mathfrak{K}_1(t)\right\|_{\mathcal{Z}_2} \leq C(1+t)^{-\frac{3}{4}}\|U_{0}\|_{{\mathcal{Z}_{2}}\bigcap {\mathcal{Z}_{1}}} .
\end{align*}
Using the H\"{o}lder’s and  Sobolev's inequalities and \eqref{C4} to $\mathfrak{K}_2$, we deduce that
\begin{align*}
\left\|\mathfrak{K}_2(t)\right\|_{\mathcal{Z}_2}^2 \leq \,& C \int_0^t(1+t-s)^{-\frac{3}{2}}\left(\left\|u \cdot\left\{\mathbf{I}-\mathbf{P}_{\mathbf{0}}-\mathbf{P}_{\mathbf{1}}\right\} f\right\|_{\mathcal{Z}_1 \cap \mathcal{Z}_2}^2\right. \\
& \left.+\|\nu^{-\frac{1}{2}} \theta \left\{\mathbf{I}-\mathbf{P}_2\right\} \frac{1}{\sqrt{M}} \nabla_v\cdot \left(\sqrt{M}(\nabla_v f-\frac{v}{2} f)\right)\|_{\mathcal{Z}_1 \cap \mathcal{Z}_2}^2\right) \mathrm{d} s \\
\leq\, & C \int_0^t(1+t-s)^{-\frac{3}{2}} \mathcal{E}^2(s) \mathrm{d} s+C \int_0^t(1+t-s)^{-\frac{3}{2}}\|\theta\|_{H^2}^2\|\{\mathbf{I}-\mathbf{P}\} f\|^2 \mathrm{d} s \\
\leq \,& C \int_0^t(1+t-s)^{-\frac{3}{2}} \mathcal{E}^2(s) \mathrm{d} s.
\end{align*}
Similarly, we    compute $\mathfrak{K}_3$ and $\mathfrak{K}_4$ directly as follows
\begin{align*}
\left\|\mathfrak{K}_3(t)\right\|_{\mathcal{Z}_2} & \leq C \int_0^t(1+t-s)^{-\frac{3}{4}}\|\big(a u \cdot v M^{\frac{1}{2}}, u \cdot b\sqrt{M}\big)\|_{\mathcal{Z}_1 \cap \mathcal{Z}_2}\mathrm{d} s \\
 & \leq C \int_0^t(1+t-s)^{-\frac{3}{4}}\|\big(a \theta\left(|v|^2-3\right) \sqrt{M},u\cdot v(b\cdot v\sqrt{M})\big)\|_{\mathcal{Z}_1 \cap \mathcal{Z}_2}\mathrm{d} s \\
& \leq C \int_0^t(1+t-s)^{-\frac{3}{4}} \mathcal{E}(s)\mathrm{d}s,\\
\left\|\mathfrak{K}_4(t)\right\|_{L^2} & \leq C \int_0^t(1+t-s)^{-\frac{3}{4}}\|(\mathfrak{H}_1,\mathfrak{H}_2,\mathfrak{H}_3)\|_{\mathcal{Z}_1 \cap \mathcal{Z}_2}\mathrm{d} s \\
& \leq C \int_0^t(1+t-s)^{-\frac{3}{4}} \mathcal{E}(s) \mathrm{d} s. 
\end{align*}
Therefore, we obtain  
\begin{align}
\| U\|_{\mathcal{Z}_{2}}^{2}\leq\,& C(1+t)^{-\frac{3}{4}}\|U_{0}\|_{{\mathcal{Z}_{2}}\bigcap {\mathcal{Z}_{1}}}^2+C \int_0^t(1+t-s)^{-\frac{3}{2}} \mathcal{E}^2(s) \mathrm{d} s\nonumber\\
&+C\bigg(\int_0^t(1+t-s)^{-\frac{3}{4}} \mathcal{E}(s) \mathrm{d} s\bigg)^2.
\end{align}
Define
\begin{equation}
\mathcal{E}_{\infty}(t):=\sup _{0 \leq s \leq t}(1+s)^{\frac{3}{2}} \mathcal{E}(s) .
\end{equation}
By a direct calculation,  it holds
\begin{equation}\label{G3.17}
\|U\|_{\mathcal{Z}_{2}}^{2}
\leq C(1+t)^{-\frac{3}{2}}\left(\|U_{0}\|_{{\mathcal{Z}_{2}}\bigcap\mathcal{Z}_{1}}^{2}
+\mathcal{E}_{\infty}^{2}(t)\right).
\end{equation}
Then, substituting \eqref{G3.17} into \eqref{G3.14}, we arrive at
\begin{equation}
\mathcal{E}_{\infty}(t)\leq C\Big(\|U_{0}\|_{{\mathcal{Z}_{2}}\bigcap\mathcal{Z}_{1}}^{2}
+\mathcal{E}_{\infty}^{2}(t)\Big),
\end{equation}
for all $t\geq 0$, since $\|U_0\|_{{\mathcal{H}^{2}}\bigcap\mathcal{Z}_{1}}$ is sufficient small. Obviously, 
\begin{equation*}
\mathcal{E}(t)\leq C(1+t)^{-\frac{3}{2}}\|U_0\|_{{\mathcal{H}^{2}}\bigcap\mathcal{Z}_{1}}^{2},
\end{equation*}
which implies that 
\begin{equation*}
\|f\|_{{{L_v^2(H_x^{2})}}}+\|(\rho, u, \theta)\|_{H^2} \leq C(1+t)^{-\frac{3}{4}}\|U_0\|_{{\mathcal{H}^{2}}\bigcap\mathcal{Z}_{1}},
\end{equation*}
for all $t\geq0$. Hence we  complete the proof of \eqref{G1.13}.

Next, we   prove the time decay rate \eqref{G1.14}. Define the 
high-order energy functional $\mathcal{H}(t)$ and high-order dissipation rate  $\mathcal{M}(t)$ by
\begin{align} 
\mathcal{H}(t) & :=\sum_{1 \leq|\alpha| \leq 2}\left\|\partial^\alpha (f, \rho, u, \theta)\right\|^2+\tau_4\mathcal{E}_0^{\prime}(t)+\tau_5\sum_{|\alpha|= 1} \int_{\mathbb{R}^3} \partial^\alpha u \cdot \partial^\alpha \nabla_x \rho \mathrm{d} x\nonumber\\
&\quad\ +\tau_6\sum_{k=1} C_k \sum_{\substack{|\beta|=k \\
|\alpha|+|\beta| \leq 2}}\left\|\partial_\beta^\alpha\{\mathbf{I}-\mathbf{P}\} f\right\|^2, \label{hh}\\
\mathcal{M}(t) & := \sum_{1\leq|\alpha| \leq 2}\left(\left\|\{\mathbf{I}-\mathbf{P}\} \partial^\alpha f\right\|_\nu^2+\left\|\partial^\alpha(b-u)\right\|^2+\|\partial^\alpha(\sqrt{2} \omega-\sqrt{3} \theta)\|^2+\left\|\nabla_x \partial^\alpha(u, \theta)\right\|^2\right)\nonumber\\
&\quad\ +\sum_{|\alpha|=1}\|\partial^{\alpha}\nabla_x(a,b,\rho,\omega)\|^2+\sum_{\substack{| \alpha|=1,|\beta|= 1}}\left\|\partial_\beta^\alpha\{\mathbf{I}-\mathbf{P}\} f\right\|_\nu^2,\label{mm}
\end{align}   
where $\tau_4, \tau_5$ and $\tau_6$ are sufficiently small constants and  
\begin{align}
\mathcal{E}_0^{\prime}(t):= & \sum_{|\alpha| = 1} \sum_{i, j} \int_{\mathbb{R}^3} \partial^\alpha\left(\partial_j b_i+\partial_i b_j\right) \partial_x^\alpha \Gamma_{i, j}(\{\mathbf{I}-\mathbf{P}\} f) \mathrm{d} x\nonumber \\
& +\sum_{|\alpha| =1} \sum_i \int_{\mathbb{R}^3} \partial^\alpha \partial_i \omega \partial^\alpha \Upsilon_i\{\mathbf{I}-\mathbf{P}\} f \mathrm{d} x \nonumber\\
& +\frac{2}{21} \sum_{|\alpha| = 1} \int_{\mathbb{R}^3} \partial_x^\alpha a \partial^\alpha\bigg(\frac{\sqrt{6}}{5} \sum_i \partial_i \Upsilon_i\{\mathbf{I}-\mathbf{P}\} f-{\rm div}_x  b\bigg) \mathrm{d} x.
\end{align}

By taking the similar arguments to those used in the proofs of  Lemmas \ref{L2.3} to   \ref{L2.6}, we obtain
\begin{equation}\label{G3.22}
\frac{\mathrm{d}}{\mathrm{d} t} \mathcal{H}(t)+\lambda_{20} \mathcal{M}(t) \leq C\left(\mathcal{H}^{\frac{1}{2}}(t)+\mathcal{H}(t)+\mathcal{H}^2(t)\right) \mathcal{M}(t),
\end{equation}
for some $\lambda_{20}>0$.
Adding $\tilde \lambda_2 \|\nabla U\|_{\mathcal{Z}_2}^2$ for some 
$\tilde \lambda_2>0$ to both sides of the inequality  \eqref{G3.22} yields
\begin{equation}
\frac{\mathrm{d}}{\mathrm{d} t} \mathcal{H}(t)+\lambda_{21} \mathcal{M}(t) \leq C\|\nabla U\|_{\mathcal{Z}_2}^2,
\end{equation}
for some $\lambda_{21}>0$.

Now we define
\begin{equation}\label{G3.24}
\mathcal{H}_{\infty}(t):=\sup _{0 \leq s \leq t}(1+s)^{\frac{5}{2}} \mathcal{H}(s) .
\end{equation}
By taking the similar arguments   to the estimates of  $\|U\|_{\mathcal{Z}_{2}}$ and 
the proof of Theorem \ref{C3}, we obtain
\begin{equation}\label{D15}
\|\nabla U\|_{_{\mathcal{Z}_{2}}}^{2}\leq C(1+t)^{-\frac{5}{2}}
\left(\|U_0\|_{{\mathcal{H}^{2}}\bigcap\mathcal{Z}_{1}}^{2}
+\mathcal{H}_{\infty}^{2}(t)\right).
\end{equation}
Thanks to  \eqref{G3.22}, \eqref{G3.24} and Gronwall's inequality, we arrive at 
\begin{equation}
\mathcal{H}(t) \leq e^{-\lambda_{21} t} \mathcal{H}(0)+ C(1+t)^{-\frac{5}{2}}
\left(\|U_0\|_{{\mathcal{H}^{2}}\bigcap\mathcal{Z}_{1}}^{2}
+\mathcal{H}_{\infty}^{2}(t)\right).
\end{equation}
Therefore,
\begin{equation*}
\mathcal{H}_{\infty}(t)\leq C\left(\|U_0\|_{{\mathcal{H}^{2}}\bigcap\mathcal{Z}_{1}}^{2}
+\mathcal{H}_{\infty}^{2}(t)\right).
\end{equation*}
Since $\|U_0\|_{{\mathcal{H}^{2}}\bigcap\mathcal{Z}_{1}}^{2}$ is sufficient small, the above inequality implies that
\begin{equation*}
\mathcal{H}_{\infty}(t)\leq C\|U_0\|_{{\mathcal{H}^{2}}\bigcap\mathcal{Z}_{1}}^{2},
\end{equation*}
for all $t\geq 0$. 
 Here and below we denote the norm $\|(f_0,\rho_0, u_0,\theta_0)\|_{\mathcal{H}^{m}}$   as
\begin{align*}
\|(f_0,\rho_0, u_0,\theta_0)\|_{\mathcal{H}^{m}}&=\|\rho_0\|_{H^{m}}+\|u_0\|_{H^{m}}+\|\theta_0\|_{H^{m}}+\|f_0\|_{H^{m}_{x,v}},    
\end{align*}
for presentation simplicity.

By the definition of $\mathcal{H}_{\infty}(t)$, we get
\begin{equation*}
\mathcal{H}(t)\leq C(1+t)^{-\frac{5}{2}}\|U_0\|_{{\mathcal{H}^{2}}\bigcap\mathcal{Z}_{1}}^{2},
\end{equation*}
which implies that 
\begin{equation*}
\|\nabla_x f\|_{L_v^2(H_{x}^{1})}+\|\nabla_x(\rho, u, \theta,f)\|_{H^1} \leq C(1+t)^{-\frac{5}{4}},
\end{equation*}
for all $t\geq 0$. Hence we get \eqref{G1.14}.
Thus, we finally finish the proofs on the decay rates of the 
strong solutions and their gradients.  To complete the whole  proof 
 of Theorem \ref{T1.1}, the remainder is to show the time-decay rates of the highest-order derivatives of solutions, 
 this is the task for next section. 

 {\section{Time-decay rates of the highest-order derivatives of strong solutions}}

In this section, we shall establish the optimal time-decay rates of the
highest-order derivatives of $(f,\rho,u,\theta)$
and complete the whole proof of Theorem \ref{T1.1}. 

For a  function $g(x)\in L^2(\mathbb{R}^3)$, we   define a frequency decomposition of it as follows
\begin{align}
g^{L}(x)=\phi_0(D_x)g(x),\quad g^{H}(x)=\phi_1(D_x)g(x),
\end{align}
where $D_x=\frac{1}{\sqrt{-1}}{(\partial_{x_1},\partial_{x_2},\partial_{x_3})}$, and  $\phi_0(D_x)$ 
and $\phi_1(D_x)$ are two pseudo-differential operators
with regard to the smooth cut-off functions $\phi_0(k)$ and $\phi_1(k)$ satisfying
$0\leq \phi_0(k)\leq 1$, $\phi_1(k)=1-\phi_0(k)$, and
\begin{align}\nonumber%
\phi_0(k)=\begin{cases}1, \quad|k|\leq \frac{r_0}{2},\\0, \quad|k|>{r_0},\end{cases}
\end{align}
for a fixed constant $r_0>0$.
Notice that $ g^{L}(x)$ and $ g^{H}(x) $ are low-frequency and high-frequency 
of $g$ satisfying the following identity 
\begin{align}
 g(x)=g^{L}(x)+g^{H}(x).   
\end{align}

By Plancherel's theorem, we have the following results.
\begin{lem}[\cite{WW-SCM-2022}]\label{NJKL4.1}
Let $g\in H^2(\mathbb{R}^3)$. Then, it hold 
\begin{gather*}
 \|g^H\|_{L^2}\leq C\|\nabla g\|_{L^2},\quad \|g^H\|_{L^2}\leq C\|\nabla^2 g\|_{L^2}, \\
 \|\nabla^2 g^L\|_{L^2}\leq C\|\nabla g^L\|_{L^2}, 
\end{gather*}
for some constant $C>0$.
\end{lem}

In order to obtain the estimate of second-order derivatives of strong solutions,
we need the following lemma.
\begin{lem}[\cite{LZ-Cpaa-2020}]\label{NJKL4.2}
Let $h$ and $g$ be two Schwarz functions. Then, for $r\geq 0$, one has
\begin{align}
\|\nabla^{k}(gh) \|_{L^r\left(\mathbb{R}^3\right)} & \leq C\|g\|_{L^{r_1}\left(\mathbb{R}^3\right)}\|\nabla^{k}h\|_{L^{r_2}\left(\mathbb{R}^3\right)}+C\|h\|_{L^{r_3}\left(\mathbb{R}^3\right)}\|\nabla^{k}g\|_{L^{r_4}\left(\mathbb{R}^3\right)},
\end{align}
with $1<r,r_2,r_4<\infty$ and $r_i(i=1,\dots,4)$ satisfy the following 
\begin{align*}
\frac{1}{r_1}+\frac{1}{r_2}=\frac{1}{r_3}+\frac{1}{r_4}=\frac{1}{r}.   
\end{align*}
\end{lem}

Now we give the following estimates on the  second-order derivatives of $(f,\rho,u,\theta)$. 
\begin{lem}\label{NJKL4.3}
For strong solutions to the problem \eqref{G1}--\eqref{G5},  there exists a positive constant $\lambda_{22}>0$, such that
\begin{align}\label{NJKG4.7}
&\frac{1}{2} \frac{\mathrm{d}}{\mathrm{d} t}  \|\nabla_x^2 (f, \rho, u,\theta )\|^2+ \lambda_{22} \big\{\|\{\mathbf{I}-\mathbf{P}\} \nabla_x^2 f\|_\nu^2 +\left\|\nabla_x^2(b-u)\right\|^2\nonumber \\
&\qquad\quad+\|\nabla_x^2(\sqrt{2} \omega-\sqrt{3} \theta)\|^2+\left\|\nabla_x^3(u, \theta)\right\|^2\big\} 
  \leq C\delta\|\nabla^2(a,b,w,\rho)\|_{L^2}^2,
\end{align}
for all $0 \leq t<T$. 
\end{lem}    
\begin{rem}
In fact, based on the previous results, we only  need to consider the case   $|\alpha|=2$ in \eqref{G2.12}. Using
  Lemmas \ref{L2.1}--\ref{NJKL4.2}, and H\"{o}lder's, Sobolev's and Young's inequalities, it's easy to get \eqref{NJKG4.7} by taking  similar arguments  to the process
in Lemma \ref{L2.3}. For brevity, we omit the proof here.
\end{rem}

Below, we begin to  estimate the dissipation of $(a^H, b^H, \omega^H,\rho^H)$.
Applying $\phi_1(D_x)$ to the system \eqref{G2.13}--\eqref{G2.17}, one has
\begin{align}\label{NJK4.8}
&\partial_t a^H+{\rm div}_x  b^H=0,\\
\label{NJK4.9}
& \partial_t b_i^H+\partial_i a^H+\frac{2}{\sqrt{6}} \partial_i \omega^H+\sum_{j=1}^3 \partial_{x_j} \Gamma_{i, j}(\{\mathbf{I}-\mathbf{P}\} f)^H=-b_i^H+u_i^H+(u_i a)^H, \\
\label{NJK4.10}
& \partial_t \omega^H+\sqrt{2}(\sqrt{2} \omega^H-\sqrt{3} \theta^H)-\sqrt{6} a \theta^H+\frac{2}{\sqrt{6}} {\rm div}_x  b^H-\frac{2}{\sqrt{6}} (u \cdot b)^H \nonumber\\
& \quad\ \ +\sum_{i=1}^3 \partial_{x_i} \Upsilon_i(\{\mathbf{I}-\mathbf{P}\} f)^H=0, \\
\label{NJK4.11}
& \partial_j b_i^H+\partial_i b_j^H-\left(u_i b_j+u_j b_i\right)^H\nonumber \\
& \quad\ \ -\frac{2}{\sqrt{6}} \delta_{i j}\Big(\frac{2}{\sqrt{6}} {\rm div}_x  b^H-\frac{2}{\sqrt{6}} (u \cdot b)^H+\sum_{i=1}^3 \partial_{x_i} \Upsilon_i(\{\mathbf{I}-\mathbf{P}\} f)^H\Big)\nonumber \\
&\ \ =-\partial_t \Gamma_{i, j}(\{\mathbf{I}-\mathbf{P}\} f)^H+\Gamma_{i, j}(\mathfrak{l}+\mathfrak{r}+\mathfrak{s})^H, \\
\label{NJK4.12}
& \frac{5}{3}\left(\partial_i \omega^H-\omega u_i^H-\sqrt{6} (\theta b_i)^H\right)-\frac{2}{\sqrt{6}} \sum_{j=1}^3 \partial_{x_j} \Gamma_{i, j}(\{\mathbf{I}-\mathbf{P}\} f)^H\nonumber \\
&\ \ =-\partial_t \Upsilon_i(\{\mathbf{I}-\mathbf{P}\} f)^H+\Upsilon_i(\mathfrak{l}+\mathfrak{r}+\mathfrak{s})^H,
\end{align}
where $\mathfrak{l}, \mathfrak{r}$ and  $\mathfrak{s}$ are defined by \eqref{NJK2.181}--\eqref{NJK2.18}.
Therefore, we    define the temporal functional $\mathcal{E}_0^H(t)$ as follows  
\begin{align}\label{NJK4.13}
\mathcal{E}_0^H(t) =\, & \sum_{i, j} \int_{\mathbb{R}^3} \nabla_x\left(\partial_j b_i^H+\partial_i b_j^H\right) \nabla_x \Gamma_{i, j}(\{\mathbf{I}-\mathbf{P}\} f)^H \mathrm{d} x \nonumber\\
& + \sum_i \int_{\mathbb{R}^3} \nabla_x \partial_i \omega^H \nabla_x \Upsilon_i(\{\mathbf{I}-\mathbf{P}\} f)^H \mathrm{d} x \nonumber\\
& +\frac{2}{21} \int_{\mathbb{R}^3} \nabla_x a^H \nabla_x\bigg(\frac{\sqrt{6}}{5} \sum_i \partial_i \Upsilon_i(\{\mathbf{I}-\mathbf{P}\} f)^H-{\rm div}_x  b^H\bigg) \mathrm{d} x .
\end{align} 

Applying  similar  arguments  to that   in the proof of Lemma \ref{L2.4}, 
we have the following result. 
\begin{lem}\label{NJKL4.4}
For strong solutions to the system \eqref{G1}-\eqref{G5}, there exists a positive constant $\lambda_{23}>0$, such that
\begin{align}\label{NJK4.14}
& \frac{\mathrm{d}}{\mathrm{d} t} \mathcal{E}_0^H(t)+\lambda_{23}\big(\|\nabla_x^2 a^H\|_{L^2}^2+\|\nabla_x^2 b^H\|_{L^2}^2+\|\nabla_x^2 \omega^H\|_{L^2}^2 \big)\nonumber\\
\leq\, & C\left(\|\{\mathbf{I}-\mathbf{P}\}\nabla^2_x f\|_{L_v^2\left(L_x^2\right)}^2+\|\nabla_x^2(u-b)\|_{L^2}^2+\|\nabla_x^2(\sqrt{2} \omega-\sqrt{3} \theta)\|_{L^2}^2\right) \nonumber\\
& +C\delta\|\nabla^2(a,b,w,u)\|_{L^2}^2,
\end{align}
for all $0 \leq t<T$. 
\end{lem}   
\begin{rem}
Thanks to the Lemma \ref{NJKL4.1}, we    deduce that
\begin{align*}
\|(u-b)^H\|_{H^1}^2\leq& \,C\|\nabla_x^2(u-b)\|_{L^2}^2,\\
\|\{\mathbf{I}-\mathbf{P}\}\nabla_x f\|_{L_v^2\left(H_x^1\right)}^2\leq\,&
\|\{\mathbf{I}-\mathbf{P}\}\nabla^2_x f^H\|_{L_v^2\left(L_x^2\right)}^2,\\
\|\nabla_x (\sqrt{2} \omega-\sqrt{3} \theta)^H\|_{L^2}^2\leq\,&
\|\nabla_x^2(\sqrt{2} \omega-\sqrt{3} \theta)\|_{L^2}^2.
\end{align*}
\end{rem}

Lemma \ref{NJKL4.1} and the above  properties of high-frequency show  their features  distinguishing  
  from the functions themselves obtained in 
Lemma \ref{L2.4}. Thus, we have obtained a  more refined estimate 
with regard to $\nabla^2(a^H,b^H,\omega^H)$. Finally, we consider the estimate for the remaining term $\nabla^2 \rho^H$  which  can  be obtained in \eqref{G2.29} with $|\alpha|=1$.

\begin{lem}\label{NJKL4.5}
For strong solutions to the problem \eqref{G1}--\eqref{G5}, there exists a positive constant $\lambda_{24}>0$, such that
\begin{align}\label{NJK4.15}
& \frac{\mathrm{d}}{\mathrm{d} t}  \int_{\mathbb{R}^3} \nabla_x u \cdot \nabla_x^2 \rho^H \mathrm{d} x+\lambda_{24}\|\nabla_x^2 \rho^H\|_{L^2}^2\nonumber \\
\, & \qquad \leq C\left(\|\nabla_x^2(u-b)\|_{L^2}^2+\|\nabla_x^2 \rho^L \|_{L^2}^2+\|\nabla_x^2 \theta\|_{H^1}^2+\|\nabla_x^2 u\|_{H^1}^2\right)\nonumber \\
&\qquad \qquad +C\delta\|\nabla^2(a,b,w,\rho,u)\|_{L^2}^2,
\end{align}
for all $0 \leq t<T$ with any $T>0$.
\end{lem}

\begin{rem}
By  applying  Lemmas \ref{L2.1}, \ref{NJKL4.1}, and \ref{NJKL4.2}, and the smallness  assumption \eqref{G2.1}, we  easily obtain \eqref{NJK4.15} through the same way used in Lemma
\ref{L2.5}. For simplicity, we omit the proof here.
\end{rem}

In the following, we define the energy functional ${\mathcal{E}}_1(t)$
and the corresponding dissipation $\mathcal{D}_1(t)$:
\begin{align*}
\mathcal{E}_1(t):=\,&\|\nabla_x^2 (f, \rho, u,\theta )\|^2+\tau_7\mathcal{E}_0^H(t)+ \tau_8 \int_{\mathbb{R}^3} \nabla_x u \cdot \nabla_x^2 \rho^H \mathrm{d} x ,\\
\mathcal{D}_1(t):=\,& \big(\|\{\mathbf{I}-\mathbf{P}\} \nabla_x^2 f\|_\nu^2+\left\|\nabla_x^2(b-u)\right\|^2+\|\nabla_x^2(\sqrt{2} \omega-\sqrt{3} \theta)\|^2+\left\|\nabla_x^3(u, \theta)\right\|^2\big)\\
&+\big(\|\nabla_x^2 a^H\|_{L^2}^2+\|\nabla_x^2 b^H\|_{L^2}^2+\|\nabla_x^2 \omega^H\|_{L^2}^2 \big)+\|\nabla_x^2 \rho^H\|_{L^2}^2,
\end{align*}
where $0<\tau_7, \tau_8\ll 1$ are small constants to be determined later.

Adding $\eqref{NJKG4.7}$ to $\tau_7\times \eqref{NJK4.14}$ and
 $\tau_8\times \eqref{NJK4.15}$, we    infer that 
\begin{align}\label{NJK4.16}
\frac{\mathrm{d}}{\mathrm{d} t} \mathcal{E}_1(t)+\lambda_{25} \mathcal{D}_1(t) \leq C \|\nabla^2_x (a^L, b^L, \omega^L, \rho^L, u^L, \theta^L)\|_{L^2}^2,    
\end{align} 
for some $\lambda_{25}>0$.
Here, we have utilized   Lemma \ref{NJKL4.1} and the following facts:
\begin{align}
\|\nabla^2_x(a, b, \omega, \rho)\|\leq\,& 
\|\nabla^2_x(a^L, b^L, \omega^L, \rho^L)\|+\|\nabla^2_x(a^H, b^H, \omega^H, \rho^H)\|,\\
\|\nabla^2_x(u, \theta)\|_{L^2}^2\leq\,&
2\big(\|\nabla^2_x(u^L, \theta^L)\|_{L^2}^2+\|\nabla^2_x(u^H, \theta^H)\|_{L^2}^2\big)\nonumber\\
\leq\,& C\big(\|\nabla^2_x(u^L, \theta^L)\|_{L^2}^2+\|\nabla^3_x(u, \theta)\|_{L^2}^2      \big).
\end{align}
Using Young's inequality, we know that  $\mathcal{E}_1(t)$ is equivalent 
to $\|\nabla^2_x f(t)\|_{L_v^2\left(L_x^2\right)}^2+\|\nabla^2_x(\rho,u,\theta)\|_{L^2}^2.$

Then, applying $\phi_0(D_x)$ to \eqref{NJK3.3}, one has
\begin{align}
U^{L}(t)=\mathbb{A}^{L}(t)U_{0}+\int_{0}^{t}\mathbb{A}^{L}(t-s)(\mathfrak{S}(s),0,0,0){\rm d}s,    
\end{align}
where $U^L(t)=(f^L(t),\rho^L(t),u^L(t),\theta^L(t))$.

\begin{thm}\label{NJKT4.6}
Let $1\leq q\leq2$. For any $m$, $n$ with 
$0\leq m\leq n$, we have 
\begin{align}
&\|\partial^{n}\mathbb{A}^L(t)U_{0}\|_{\mathcal{Z}_{2}}
\leq C\left(1+t\right)^{-\frac{3}{2}\left(\frac{1}{q}-\frac{1}{2}\right)-\frac{n-m}{2}}
\|\partial^{m}U_{0}\|_{\mathcal{Z}_{q}},\label{NJK4.20}\\
&\left\|\partial^{n}\int_{0}^{t}\mathbb{A}^L(t-s)(\mathfrak{S}(s),0,0,0){\rm d}s\right\|_{\mathcal{Z}_{2}}^{2}\nonumber\\
&\,\qquad \leq C\int_0^t(1+t-s)^{-3\left(\frac{1}{q}-\frac{1}{2}\right)-(n-m)}
\left\|\partial^{m}(G(s),\nu^{-1/2}h(s))\right\|_{\mathcal{Z}_{q}}^{2}
{\rm d}s,  \label{NJK4.21}
\end{align}
for all $t\geq 0$.    
\end{thm}
\begin{rem}
Through the same method used in Theorem \ref{C3}, we  can  directly obtain \eqref{NJK4.20}--\eqref{NJK4.21}. 
Compared with Theorem \ref{C3}, Theorem \ref{NJKT4.6} only contains
the low-frequency part, which is the difference between Theorem \ref{C3}
and Theorem \ref{NJKT4.6}. 
\end{rem}
At the end of this section, we shall establish the estimate of $\|\nabla^2_x f(t)\|_{L_v^2\left(L_x^2\right)}^2+\|\nabla^2_x(\rho,u,\theta)\|_{L^2}^2$.
Applying Gronwall's inequality to \eqref{NJK4.16}, we have
\begin{align}\label{NJKG4.22}
\mathcal{E}_1(t) \leq e^{-\lambda_{25} t} \mathcal{E}_1(0)+C \int_0^t e^{-\lambda_{25}(t-s)}\| \nabla^2_x U^L(s)\|_{\mathcal{Z}_{2}}^{2} \mathrm{d} s.    
\end{align}
Lemma \ref{NJKL4.1} implies
\begin{align}\label{NJK4.23}
\| \nabla^2_x U^L(t)\|_{\mathcal{Z}_{2}}^{2}\leq C\| \nabla_x U^L(t)\|_{\mathcal{Z}_{2}}^{2}.      
\end{align}
Therefore, it's essential to estimate $\| \nabla_x U^L(t)\|_{\mathcal{Z}_{2}}^{2}$. To do this, we
apply $\phi_0(D_x)$ to \eqref{NJK3.11}, then we get
\begin{align}
\nabla_x U^L(t)=&\nabla_x\mathfrak{K}^L_{1}(t)+\nabla_x\mathfrak{K}^L_{2}(t)
+\nabla_x\mathfrak{K}^L_{3}(t)+\nabla_x\mathfrak{K}^L_{4}(t).
\end{align}

Below we give the estimate of $\nabla_x\mathfrak{K}^L_{i}(t),  i=1,\dots,4$.
To begin with, we define
\begin{align}
\mathcal{E}_{2}(t):=\sup _{0 \leq s \leq t}\bigg\{&\sum_{n=0,1}(1+s)^{\frac{3}{4}+n} \Big(\|\nabla^n_x f(t)\|_{L_v^2\left(L_x^2\right)}^2+\|\nabla^n_x(\rho,u,\theta)\|_{L^2}^2\Big) \Big.\nonumber\\
\Big.&+(1+s)^{\frac{7}{4}} \Big(\|\nabla^2_x f(t)\|_{L_v^2\left(L_x^2\right)}^2+\|\nabla^2_x(\rho,u,\theta)\|_{L^2}^2\Big)     \bigg\} .    
\end{align}
Taking $q=\frac{4}{3},\  n=1$ in \eqref{NJK4.20} and using Lemma \ref{NJKL4.1}, we have
\begin{align*}
\|\nabla_x\mathfrak{K}^L_{1}(t)\|_{\mathcal{Z}_{2}}\leq  C(1+t)^{-\frac{7}{8}}
\|U_0\|_{\mathcal{Z}_{\frac{4}{3}}} 
\leq   C(1+t)^{-\frac{7}{8}}
\|U_0\|_{\mathcal{Z}_{1}}^{\frac{1}{2}}\|U_0\|_{\mathcal{Z}_{2}}^{\frac{1}{2}}.
\end{align*}
Taking $q=1,\  m=0$ and $q=\frac{3}{2},\  m=1$ in \eqref{NJK4.20} respectively, we get
\begin{align*}
\left\|\nabla_x\mathfrak{K}^L_3(t)\right\|_{\mathcal{Z}_2}  \leq\,& C \int_0^{\frac{t}{2}}(1+t-s)^{-\frac{5}{4}}\|\big(a u \cdot v M^{\frac{1}{2}}, u \cdot b\sqrt{M}\big)\|_{\mathcal{Z}_1}\mathrm{d} s \\
&+ C \int_{\frac{t}{2}}^{t}(1+t-s)^{-\frac{1}{4}}\|\nabla_x\big(a u \cdot v M^{\frac{1}{2}}, u \cdot b\sqrt{M}\big)\|_{\mathcal{Z}_{\frac{3}{2}}}\mathrm{d} s \\
\leq\,& C \int_0^{\frac{t}{2}}(1+t-s)^{-\frac{5}{4}}\big(\|au\|_{L^1}+\|bu\|_{L^1}   \big)\mathrm{d}s\\
&+ C\int_{\frac{t}{2}}^{t}(1+t-s)^{-\frac{1}{4}}\big(\|\nabla_x u\|_{L^2}\|a\|_{L^6}+\|\nabla_x a\|_{L^2}\|u\|_{L^6}   \big)\mathrm{d}s\\
&+ C\int_{\frac{t}{2}}^{t}(1+t-s)^{-\frac{1}{4}}\big(\|\nabla_x u\|_{L^2}\|b\|_{L^6}+\|\nabla_x b\|_{L^2}\|u\|_{L^6}   \big)\mathrm{d}s\\
\leq&\, C \int_0^{\frac{t}{2}}(1+t-s)^{-\frac{5}{4}}\big(\|a\|_{L^2}\|u\|_{L^2}+\|b\|_{L^2}\|u\|_{L^2}   \big)\mathrm{d}s\\
&+ C\int_{\frac{t}{2}}^{t}(1+t-s)^{-\frac{1}{4}}\big(\|\nabla_x u\|_{L^2}^2+\|\nabla_x a\|_{L^2}^2+\|\nabla_x b\|_{L^2}^2  \big)\mathrm{d}s\\
\leq&\, C\mathcal{E}_2(t) \int_0^{\frac{t}{2}}(1+t-s)^{-\frac{5}{4}}(1+s)^{-\frac{3}{4}}\mathrm{d}s\\
&+ C\int_{\frac{t}{2}}^{t}(1+t-s)^{-\frac{1}{4}}(1+s)^{-\frac{7}{4}}\mathrm{d}s\\
\leq& \,C(1+t)^{-1}\mathcal{E}_2(t),
\end{align*}
where we have used the definition of $\mathcal{E}_2(t)$. 
Similarly, taking $q=1,\  m=0$ and $q=\frac{3}{2},\  m=0$ in \eqref{NJK4.20} respectively, we have
\begin{align*}
\left\|\nabla_x\mathfrak{K}^L_4(t)\right\|_{\mathcal{Z}_2}  \leq\,& C \int_0^{\frac{t}{2}}(1+t-s)^{-\frac{5}{4}}\|(\mathfrak{H}_1,\mathfrak{H}_2,\mathfrak{H}_3)\|_{\mathcal{Z}_1}\mathrm{d} s \\
&+ C \int_{\frac{t}{2}}^{t}(1+t-s)^{-\frac{3}{4}}\|\nabla_x (\mathfrak{H}_1,\mathfrak{H}_2,\mathfrak{H}_3)\|_{\mathcal{Z}_{\frac{3}{2}}}\mathrm{d} s\\
\leq&\, C \mathcal{E}_2(t) \int_0^{\frac{t}{2}}(1+t-s)^{-\frac{5}{4}}(1+s)^{-\frac{3}{4}}\mathrm{d}s\\
&+ C\mathcal{E}_2(t)\int_{\frac{t}{2}}^{t}(1+t-s)^{-\frac{3}{4}}(1+s)^{-\frac{5}{4}}\mathrm{d}s\\
\leq&\, C(1+t)^{-1}\mathcal{E}_2(t).
\end{align*}

Taking $q=1,\  m=0$ and $q=\frac{3}{2},\  m=1$ in \eqref{NJK4.21} respectively
and using the same method above, we    obtain
\begin{align*}
&\left\|\nabla_x\mathfrak{K}^L_2(t)\right\|_{\mathcal{Z}_2}^2 \\
\leq\,& C \int_0^{\frac{t}{2}}(1+t-s)^{-\frac{5}{2}}\left\|u \cdot\left\{\mathbf{I}-\mathbf{P}_{\mathbf{0}}-\mathbf{P}_{\mathbf{1}}\right\} f\right\|_{\mathcal{Z}_1}^2\mathrm{d} s \\
&+ C \int_{\frac{t}{2}}^{t}(1+t-s)^{-\frac{1}{2}}\|\nabla_x\big(u \cdot\left\{\mathbf{I}-\mathbf{P}_{\mathbf{0}}-\mathbf{P}_{\mathbf{1}}\right\} f\big)\|_{\mathcal{Z}_{\frac{3}{2}}}^2\mathrm{d} s \\
&+ C \int_0^{\frac{t}{2}}(1+t-s)^{-\frac{5}{2}}\Big\|\nu^{-\frac{1}{2}} \theta \cdot\left\{\mathbf{I}-\mathbf{P}_2\right\} \frac{1}{\sqrt{M}} {\rm div}_v \Big(\sqrt{M}\Big(\nabla_v f-\frac{v}{2} f\Big)\Big)\Big\|_{\mathcal{Z}_1}^2\mathrm{d} s \\
&+ C \int_{\frac{t}{2}}^{t}(1+t-s)^{-\frac{1}{2}}\Big\|\nabla_x\nu^{-\frac{1}{2}} \theta \cdot\left\{\mathbf{I}-\mathbf{P}_2\right\} \frac{1}{\sqrt{M}} {\rm div}_v \Big(\sqrt{M}\Big(\nabla_v f-\frac{v}{2} f\Big)\Big)\Big\|_{\mathcal{Z}_{\frac{3}{2}}}^2\mathrm{d} s \\
\leq&\, C \mathcal{E}_2^2(t) \int_0^{\frac{t}{2}}(1+t-s)^{-\frac{5}{2}}(1+s)^{-\frac{3}{2}}\mathrm{d}s\\
&+ C\mathcal{E}_2^2(t)\int_{\frac{t}{2}}^{t}(1+t-s)^{-\frac{1}{2}}(1+s)^{-\frac{7}{2}}\mathrm{d}s\\
\leq&\, C(1+t)^{-5/2}\mathcal{E}_2^2(t).
\end{align*}
Thanks to the above estimates, we    deduce that
\begin{align}\label{NJK4.25}
\|\nabla_x U^L(t)\|_{\mathcal{Z}_{2}}\leq\,& 
C(1+t)^{-\frac{7}{8}}
\Big(\mathcal{E}_2(t)+\|U_0\|_{\mathcal{Z}_{1}}^{\frac{1}{2}}\|U_0\|_{\mathcal{Z}_{2}}^{\frac{1}{2}}\Big).
\end{align}
Similarly,
\begin{align}\label{NJK4.26}
\| U^L(t)\|_{\mathcal{Z}_{2}}\leq\,& C(1+t)^{-\frac{3}{8}}\Big(\mathcal{E}_2(t)
+\|U_0\|_{\mathcal{Z}_{1}}^{\frac{1}{2}}\|U_0\|_{\mathcal{Z}_{2}}^{\frac{1}{2}}\Big).
\end{align}

Plugging \eqref{NJK4.23} and \eqref{NJK4.25} into \eqref{NJKG4.22} yields
\begin{align}\label{NJK4.28}
\|\nabla_x^2 U(t)\|_{\mathcal{Z}_{2}}^2\leq\,& C(1+t)^{-\frac{7}{4}}\big(\mathcal{E}_2^2(t)+\|U_0\|_{\mathcal{Z}_{1}}^2+\|U_0\|_{\mathcal{H}^2}^2\big).
\end{align}
In addition, for $n=0,1$, Lemma \ref{NJKL4.1} gives
\begin{align*}
\|\nabla_x^n U(t)\|_{\mathcal{Z}_{2}}^2\leq\,& C\|\nabla_x^n U^L(t)\|_{\mathcal{Z}_{2}}^2+C\|\nabla_x^n U^H(t)\|_{\mathcal{Z}_{2}}^2\\
\leq\,& C\|\nabla_x^n U^L(t)\|_{\mathcal{Z}_{2}}^2+C\|\nabla_x^2 U(t)\|_{\mathcal{Z}_{2}}^2\\
\leq\,& C(1+t)^{-\frac{3}{4}-n}\big(\mathcal{E}_2^2(t)+\|U_0\|_{\mathcal{Z}_{1}}^2+\|U_0\|_{\mathcal{H}^2}^2\big),
\end{align*}
which implies that
\begin{align}
\mathcal{E}_2(t) \leq\,& C\big(\mathcal{E}_2^2(t)+\|U_0\|_{\mathcal{Z}_{1}}^2+\|U_0\|_{\mathcal{H}^2}^2\big).   
\end{align}
Since $\|U_0\|_{\mathcal{Z}_{1}}^2+\|U_0\|_{\mathcal{H}^2}^2$ is small enough,
it follows that
\begin{align}
\mathcal{E}_2(t) \leq\,& C\big(\|U_0\|_{\mathcal{Z}_{1}}^2+\|U_0\|_{\mathcal{H}^2}^2\big),   
\end{align}
for all $t\geq 0$. Through the definition of $\mathcal{E}_2(t)$, we get
\begin{align}\label{NJKG4.31}
\| f(t)\|_{L_v^2(L_{x}^{2})}+\|(\rho, u, \theta)(t)\|_{L^2(\mathbb{R}^3)} \leq C(1+t)^{-\frac{3}{8}}\big(\|U_0\|_{\mathcal{Z}_{1}}+\|U_0\|_{\mathcal{H}^2}\big),\nonumber\\
\|\nabla_x f(t)\|_{L_v^2(L_{x}^{2})}+\|\nabla_x (\rho, u, \theta)(t)\|_{L^2(\mathbb{R}^3)} \leq C(1+t)^{-\frac{7}{8}}\big(\|U_0\|_{\mathcal{Z}_{1}}+\|U_0\|_{\mathcal{H}^2}\big),\nonumber\\
\|\nabla_x^2 f(t)\|_{L_v^2(L_{x}^{2})}+\|\nabla_x^2(\rho, u, \theta)(t)\|_{L^2(\mathbb{R}^3)} \leq C(1+t)^{-\frac{7}{8}}\big(\|U_0\|_{\mathcal{Z}_{1}}+\|U_0\|_{\mathcal{H}^2}\big).
\end{align}

In the following we give the optimal time-decay rates of $\|\nabla^2_x f(t)\|_{L_v^2\left(L_x^2\right)}^2+\|\nabla^2_x(\rho,u,\theta)\|_{L^2}^2$.
Define
\begin{align}
\mathcal{T}_{\infty}(t):=\sup_{0\leq s \leq t}\sum_{n=0}^2(1+s)^{\frac{3}{2}+n}\big(\|\nabla^n_x f(t)\|_{L_v^2\left(L_x^2\right)}^2+\|\nabla^n_x(\rho,u,\theta)\|_{L^2}^2\big).    
\end{align}
Taking $m=0,\  q=1$ in Theorem \ref{NJKT4.6}, we have
\begin{align*}
\|\nabla^n_x\mathfrak{K}^L_{1}(t)\|_{\mathcal{Z}_{2}}\leq\,& C(1+t)^{-\frac{3}{4}-\frac{n}{2}}
\|U_0\|_{\mathcal{Z}_1}, \quad n=0,1,2.
\end{align*}
Similarly, using \eqref{NJKG4.31} and taking $m=0,\  q=1$ and $m=1, \ q=-\frac{3}{2} $ in \eqref{NJK4.20} respectively, yields
\begin{align*}
\left\|\nabla_x^n\mathfrak{K}^L_3(t)\right\|_{\mathcal{Z}_2}  \leq\,& C \int_0^{\frac{t}{2}}(1+t-s)^{-\frac{3}{4}-\frac{n}{2}}\|\big(a u \cdot v M^{\frac{1}{2}}, u \cdot b\sqrt{M}\big)\|_{\mathcal{Z}_1}\mathrm{d} s \\
&+ C \int_{\frac{t}{2}}^{t}(1+t-s)^{-\frac{1}{4}}\|\nabla_x^n\big(a u \cdot v M^{\frac{1}{2}}, u \cdot b\sqrt{M}\big)\|_{\mathcal{Z}_{\frac{3}{2}}}\mathrm{d} s \\
\leq\,& C \int_0^{\frac{t}{2}}(1+t-s)^{-\frac{3}{4}-\frac{n}{2}}\big(\|a\|_{L^2}\|u\|_{L^2}+\|b\|_{L^2}\|u\|_{L^2}   \big)\mathrm{d}s\\
&+ C\int_{\frac{t}{2}}^{t}(1+t-s)^{-\frac{1}{4}}\big(\|\nabla_x^n u\|_{L^2}\|\nabla a\|_{L^2}+\|\nabla_x^n a\|_{L^2}\|\nabla u\|_{L^2}   \big)\mathrm{d}s\\
&+ C\int_{\frac{t}{2}}^{t}(1+t-s)^{-\frac{1}{4}}\big(\|\nabla_x^n u\|_{L^2}\|\nabla b\|_{L^2}+\|\nabla_x^n b\|_{L^2}\|\nabla u\|_{L^2}   \big)\mathrm{d}s\\
\leq&\, C\big(\|U_0\|_{\mathcal{Z}_{1}}+\|U_0\|_{\mathcal{H}^2}\big)\mathcal{T}_{\infty}^{\frac{1}{2}} (t) \int_0^{\frac{t}{2}}(1+t-s)^{-\frac{3}{4}-\frac{n}{2}}(1+s)^{-\frac{9}{8}}\mathrm{d}s\\
&+ C\big(\|U_0\|_{\mathcal{Z}_{1}}+\|U_0\|_{\mathcal{H}^2}\big)\mathcal{T}_{\infty}^{\frac{1}{2}} (t)\int_{\frac{t}{2}}^{t}(1+t-s)^{-\frac{1}{4}}(1+s)^{-\frac{13}{8}-\frac{n}{2}}\mathrm{d}s\\
\leq& \,C(1+t)^{-\frac{3}{4}-\frac{n}{2}}\big(\|U_0\|_{\mathcal{Z}_{1}}+\|U_0\|_{\mathcal{H}^2}\big)\mathcal{T}_{\infty}^{\frac{1}{2}} (t), \quad n=0,1,2.
\end{align*}
For the term $\nabla_x^n\mathfrak{K}^L_2(t)$, using \eqref{NJKG4.31} and taking $q=1, \ q=\frac{3}{2}$ in \eqref{NJK4.21}, we obtain
\begin{align*}
&\left\|\nabla_x^n\mathfrak{K}^L_2(t)\right\|_{\mathcal{Z}_2}^2 \\
\leq\,& C \int_0^{\frac{t}{2}}(1+t-s)^{-\frac{3}{2}-n}\left\|u \cdot\left\{\mathbf{I}-\mathbf{P}_{\mathbf{0}}-\mathbf{P}_{\mathbf{1}}\right\} f\right\|_{\mathcal{Z}_1}^2\mathrm{d} s \\
&+ C \int_{\frac{t}{2}}^{t}(1+t-s)^{-\frac{1}{2}}\|\nabla_x^n\big(u \cdot\left\{\mathbf{I}-\mathbf{P}_{\mathbf{0}}-\mathbf{P}_{\mathbf{1}}\right\} f\big)\|_{\mathcal{Z}_{\frac{3}{2}}}^2\mathrm{d} s \\
&+ C \int_0^{\frac{t}{2}}(1+t-s)^{-\frac{3}{2}-n}\Big\|\nu^{-\frac{1}{2}} \theta \cdot\left\{\mathbf{I}-\mathbf{P}_2\right\} \frac{1}{\sqrt{M}} \nabla_v\cdot \left(\sqrt{M}\Big(\nabla_v f-\frac{v}{2} f\Big)\right)\Big\|_{\mathcal{Z}_1}^2\mathrm{d} s \\
&+ C \int_{\frac{t}{2}}^{t}(1+t-s)^{-\frac{1}{2}}\Big\|\nabla^n_x\nu^{-\frac{1}{2}} \theta \cdot\left\{\mathbf{I}-\mathbf{P}_2\right\} \frac{1}{\sqrt{M}} \nabla_v\cdot \left(\sqrt{M}\Big(\nabla_v f-\frac{v}{2} f\Big)\right)\Big\|_{\mathcal{Z}_{\frac{3}{2}}}^2\mathrm{d} s \\
\leq&\, C(1+t)^{-\frac{3}{2}-n}\big(\|U_0\|_{\mathcal{Z}_{1}}^2+\|U_0\|_{\mathcal{H}^2}^2\big)\mathcal{T}_{\infty} (t), \quad n=0,1,2.
\end{align*}
For the remaining term $\nabla_x^n\mathfrak{K}^L_4(t)$, it follows from
\eqref{NJK4.20} and \eqref{NJKG4.31} that
\begin{align*}
\left\|\nabla_x^n\mathfrak{K}^L_4(t)\right\|_{\mathcal{Z}_2}  \leq\,& C \int_0^{\frac{t}{2}}(1+t-s)^{-\frac{3}{4}-\frac{n}{2}}\|(\mathfrak{H}_1,\mathfrak{H}_2,\mathfrak{H}_3)\|_{\mathcal{Z}_1}\mathrm{d} s \\
&+ C \int_{\frac{t}{2}}^{t}(1+t-s)^{-\frac{1}{4}-\frac{n}{2}}\Big(\| \mathfrak{H}_1\|_{\mathcal{Z}_{\frac{3}{2}}}+\Big\|\mathfrak{H}_2-\frac{\rho}{1+\rho}(u-b)+\frac{1}{1+\rho} a u\Big\|_{\mathcal{Z}_{\frac{3}{2}}}\big)\mathrm{d} s\\
&+ C \int_{\frac{t}{2}}^{t}(1+t-s)^{-\frac{1}{4}-\frac{n}{2}}\Big(\Big\|\mathfrak{H}_3-\frac{\sqrt{3} \rho}{1+\rho}(\sqrt{3} \theta-\sqrt{2} \omega)\Big.\\
\Big.&\qquad\qquad-\frac{1}{1+\rho}\big((1+a)|u|^2-3 a \theta-2 u \cdot b\big)\Big\|_{\mathcal{Z}_{\frac{3}{2}}}\Big)\mathrm{d} s\\
&+ C \int_{\frac{t}{2}}^{t}(1+t-s)^{-\frac{1}{4}}\Big\|\nabla_x^n\Big(\frac{\rho}{1+\rho}(u-b)-\frac{1}{1+\rho} a u\Big)\Big\|_{\mathcal{Z}_{\frac{3}{2}}}\mathrm{d} s\\
&+ C \int_{\frac{t}{2}}^{t}(1+t-s)^{-\frac{1}{4}}\Big\|\nabla_x^n\Big(\frac{\sqrt{3} \rho}{1+\rho}(\sqrt{3} \theta-\sqrt{2} \omega)\Big)\Big\|_{\mathcal{Z}_{\frac{3}{2}}}\mathrm{d} s\\
&+ C \int_{\frac{t}{2}}^{t}(1+t-s)^{-\frac{1}{4}}\Big\|\nabla_x^n\Big(\frac{1}{1+\rho}\big((1+a)|u|^2-3 a \theta-2 u \cdot b\big)\Big)\Big\|_{\mathcal{Z}_{\frac{3}{2}}}\mathrm{d} s\\
\leq&\, C\big(\|U_0\|_{\mathcal{Z}_{1}}+\|U_0\|_{\mathcal{H}^2}\big)\mathcal{T}_{\infty}^{\frac{1}{2}} (t) \int_0^{\frac{t}{2}}(1+t-s)^{-\frac{3}{4}-\frac{n}{2}}(1+s)^{-\frac{9}{8}}\mathrm{d}s\\
&+ C\big(\|U_0\|_{\mathcal{Z}_{1}}+\|U_0\|_{\mathcal{H}^2}\big)\mathcal{T}_{\infty}^{\frac{1}{2}} (t)\int_{\frac{t}{2}}^{t}(1+t-s)^{-\frac{1}{4}-\frac{n}{2}}(1+s)^{-\frac{17}{8}}\mathrm{d}s\\
&+ C\big(\|U_0\|_{\mathcal{Z}_{1}}+\|U_0\|_{\mathcal{H}^2}\big)\mathcal{T}_{\infty}^{\frac{1}{2}} (t)\int_{\frac{t}{2}}^{t}(1+t-s)^{-\frac{1}{4}}(1+s)^{-\frac{13}{8}-\frac{n}{2}}\mathrm{d}s\\
\leq& \,C(1+t)^{-\frac{3}{4}-\frac{n}{2}}\big(\|U_0\|_{\mathcal{Z}_{1}}+\|U_0\|_{\mathcal{H}^2}\big)\mathcal{T}_{\infty}^{\frac{1}{2}} (t),\quad n=0, 1.
\end{align*}
Similar to the above estimate, for $n=2$, it holds
\begin{align*}
\left\|\nabla_x^2\mathfrak{K}^L_4(t)\right\|_{\mathcal{Z}_2}  \leq\,& C \int_0^{\frac{t}{2}}(1+t-s)^{-\frac{7}{4}}\|(\mathfrak{H}_1,\mathfrak{H}_2,\mathfrak{H}_3)\|_{\mathcal{Z}_1}\mathrm{d} s \\
&+ C \int_{\frac{t}{2}}^{t}(1+t-s)^{-\frac{5}{4}}\Big(\Big\| \mathfrak{H}_1\|_{\mathcal{Z}_{\frac{3}{2}}}+\Big\|\mathfrak{H}_2-\frac{\rho}{1+\rho}(u-b)+\frac{1}{1+\rho} a u\Big\|_{\mathcal{Z}_{\frac{3}{2}}}\Big)\mathrm{d} s\\
&+ C \int_{\frac{t}{2}}^{t}(1+t-s)^{-\frac{5}{4}}\Big(\Big\|\mathfrak{H}_3-\frac{\sqrt{3} \rho}{1+\rho}(\sqrt{3} \theta-\sqrt{2} \omega)\Big.\\
\Big.&\qquad\qquad-\frac{1}{1+\rho}\big((1+a)|u|^2-3 a \theta-2 u \cdot b\big)\Big\|_{\mathcal{Z}_{\frac{3}{2}}}\Big)\mathrm{d} s\\
&+ C \int_{\frac{t}{2}}^{t}(1+t-s)^{-\frac{3}{4}}\Big\|\nabla_x\Big(\frac{\rho}{1+\rho}(u-b)-\frac{1}{1+\rho} a u\Big)\Big\|_{\mathcal{Z}_{\frac{3}{2}}}\mathrm{d} s\\
&+ C \int_{\frac{t}{2}}^{t}(1+t-s)^{-\frac{3}{4}}\Big\|\nabla_x\Big(\frac{\sqrt{3} \rho}{1+\rho}(\sqrt{3} \theta-\sqrt{2} \omega)\Big)\Big\|_{\mathcal{Z}_{\frac{3}{2}}}\mathrm{d} s\\
&+ C \int_{\frac{t}{2}}^{t}(1+t-s)^{-\frac{3}{4}}\Big\|\nabla_x\Big(\frac{1}{1+\rho}\big((1+a)|u|^2-3 a \theta-2 u \cdot b\big)\Big)\Big\|_{\mathcal{Z}_{\frac{3}{2}}}\mathrm{d} s\\
\leq&\, C\big(\|U_0\|_{\mathcal{Z}_{1}}+\|U_0\|_{\mathcal{H}^2}\big)\mathcal{T}_{\infty}^{\frac{1}{2}} (t) \int_0^{\frac{t}{2}}(1+t-s)^{-\frac{7}{4}}(1+s)^{-\frac{9}{8}}\mathrm{d}s\\
&+ C\big(\|U_0\|_{\mathcal{Z}_{1}}+\|U_0\|_{\mathcal{H}^2}\big)\mathcal{T}_{\infty}^{\frac{1}{2}} (t)\int_{\frac{t}{2}}^{t}(1+t-s)^{-\frac{5}{4}}(1+s)^{-\frac{17}{8}}\mathrm{d}s\\
&+ C\big(\|U_0\|_{\mathcal{Z}_{1}}+\|U_0\|_{\mathcal{H}^2}\big)\mathcal{T}_{\infty}^{\frac{1}{2}} (t)\int_{\frac{t}{2}}^{t}(1+t-s)^{-\frac{3}{4}}(1+s)^{-\frac{17}{8}}\mathrm{d}s\\
\leq& \,C(1+t)^{-\frac{7}{4}}\big(\|U_0\|_{\mathcal{Z}_{1}}+\|U_0\|_{\mathcal{H}^2}\big)\mathcal{T}_{\infty}^{\frac{1}{2}} (t).
\end{align*}
Collecting  all the above estimates together, we finally get
\begin{align}\label{NJKG4.33}
\|\nabla_x^n U^L(t)\|_{\mathcal{Z}_{2}}\leq\,& 
C(1+t)^{-\frac{3}{4}-\frac{n}{2}}\big(\|U_0\|_{\mathcal{Z}_{1}}+\|U_0\|_{\mathcal{H}^2}\big)\mathcal{T}_{\infty}^{\frac{1}{2}} (t),\quad n=0,1,2.
\end{align}

Combining \eqref{NJKG4.22} with \eqref{NJKG4.33}, we deduce that
\begin{align}\label{NJKG4.34}
\mathcal{E}_1(t) \leq e^{-\lambda_{25} t} \mathcal{E}_1(0)+C(1+t)^{-\frac{7}{2}}\Big(\big(\|U_0\|_{\mathcal{Z}_{1}}^2+\|U_0\|_{\mathcal{H}^2}^2\big)\mathcal{T}_{\infty} (t)+\|U_0\|_{\mathcal{Z}_{1}}^2\Big).    
\end{align}
By the fact that $\mathcal{E}_1(t)$ is equivalent 
to $\|\nabla_x^2 U(t)\|_{\mathcal{Z}_{2}}^2$,
we have
\begin{align}\label{NJKG4.35}
\|\nabla_x^2 U(t)\|_{\mathcal{Z}_{2}}^2 \leq C(1+t)^{-\frac{7}{2}}\Big(\big(\|U_0\|_{\mathcal{Z}_{1}}^2+\|U_0\|_{\mathcal{H}^2}^2\big)\mathcal{T}_{\infty} (t)+\|U_0\|_{{\mathcal{H}^{2}}\bigcap\mathcal{Z}_{1}}^{2}\Big).    
\end{align}
Using \eqref{NJKG4.33}, \eqref{NJKG4.35} and Lemma \ref{NJKL4.1}, we arrive at
\begin{align}
\|\nabla_x^n U(t)\|_{\mathcal{Z}_{2}}^2\leq\,& C\|\nabla_x^n U^L(t)\|_{\mathcal{Z}_{2}}^2+C\|\nabla_x^n U^H(t)\|_{\mathcal{Z}_{2}}^2\nonumber\\
\leq\,& C\|\nabla_x^n U^L(t)\|_{\mathcal{Z}_{2}}^2+C\|\nabla_x^2 U(t)\|_{\mathcal{Z}_{2}}^2\nonumber\\
\leq\,& C(1+t)^{-\frac{3}{2}-n}\Big(\big(\|U_0\|_{\mathcal{Z}_{1}}^2+\|U_0\|_{\mathcal{H}^2}^2\big)\mathcal{T}_{\infty} (t)+\|U_0\|_{\mathcal{Z}_{1}}^2\Big),  \quad n=0,1.
\end{align}
Thanks to the definition of $\mathcal{T}_{\infty} (t)$, we obtain
\begin{align}
\mathcal{T}_{\infty} (t)\leq C\Big(\big(\|U_0\|_{\mathcal{Z}_{1}}^2+\|U_0\|_{\mathcal{H}^2}^2\big)\mathcal{T}_{\infty} (t)+\|U_0\|_{{\mathcal{H}^{2}}\bigcap\mathcal{Z}_{1}}^{2}\Big).    
\end{align}
Since $\|U_0\|_{\mathcal{Z}_{1}}^2+\|U_0\|_{\mathcal{H}^2}^2$ is sufficient small, one has
\begin{align}
\mathcal{T}_{\infty} (t)\leq C\|U_0\|_{{\mathcal{H}^{2}}\bigcap\mathcal{Z}_{1}}^{2},    
\end{align}
for all $t\geq 0$, which implies that 
\begin{align*}
 {\|\nabla_x^2 f(t)\|_{L_v^2(L_{x}^{2})}+\|\nabla_x^2(\rho, u, \theta)(t)\|_{L^2(\mathbb{R}^3)} }&\leq C(1+t)^{-\frac{7}{4}}, 
\end{align*}
for all $t\geq 0$.
Thus, we finally complete the proof of Theorem \ref{T1.1}.

\smallskip 

\begin{proof}[Proof of Corollary \ref{Cor1}]
Based on the results obtained in Theorem \ref{T1.1}, we can obtain 
the decay rates of solutions in general $L^p$ space with  $p\in [2,+\infty]$. 
In fact, by \eqref{G1.13} and \eqref{G1.14}, we know that  
\begin{align}
 \|(\rho,u,\theta)\|_{L^6}\leq \,& C\|\nabla_x(\rho,u,\theta)\|_{L^2}\leq C(1+t)^{-\frac{5}{4}},\nonumber \\ \label{A1}
 \| f\|_{L_v^2(L_x^{6})}\leq \,& C \|\nabla_x f\|_{L_v^2(L_x^{2})}\leq C(1+t)^{-\frac{5}{4}},     \\
 \|(\rho,u,\theta)\|_{L^2}\leq \,& C(1+t)^{-\frac{3}{4}},\nonumber\\ \label{A2}
 \|f\|_{L_v^2(L_x^{2})}\leq \,& C(1+t)^{-\frac{3}{4}}.
\end{align}

For $p\in [2,6]$, using the interpolation inequality, it follows from \eqref{A1}--\eqref{A2} that
\begin{align}
 \|(\rho,u,\theta)\|_{L^p}\leq \,& C\|(\rho,u,\theta)\|_{L^2}^{\zeta}\|(\rho,u,\theta)\|_{L^6}^{{1-\zeta}}\leq C(1+t)^{-\frac{3}{2}(1-\frac{1}{p})} , \\
 \|f\|_{L_v^2(L_x^{p})}\leq \,& C\|f\|_{L_v^2(L_x^{2})}^{\zeta}\|f\|_{L_v^2(L_x^{6})}^{{1-\zeta}}\leq C(1+t)^{-\frac{3}{2}(1-\frac{1}{p})},
\end{align}
where $\zeta= ({6-p})/ {2p} \in [0,1]$.

By Sobolev's embedding theorem, we get
\begin{align}
&\|(\rho,u,\theta)\|_{L^{\infty}}\leq C\|(\rho,u,\theta)\|_{H^{2}}\leq C(1+t)^{-\frac{3}{4}} , \\
&\|f\|_{L_v^2(L_x^{\infty})}\leq C\|f\|_{L_v^2(H_x^{2})}\leq C(1+t)^{-\frac{3}{4}}.
\end{align}

Below we   want to  achieve better results 
through the Gagliardo-Nierenberg's inequality  and the estimate \eqref{NJK1.17}. In fact, we  have 
\begin{align}
\|(\rho,u,\theta)\|_{L^{\infty}}\leq\,& C\|\nabla_x(\rho,u,\theta)\|_{L^{2}}^{\frac{1}{2}}\|\Delta_x(\rho,u,\theta)\|_{L^{2}}^{\frac{1}{2}}\nonumber\\\label{A3}
\leq\,& C(1+t)^{-\frac{3}{2}}, \\
\|f\|_{L_v^2(L_x^{\infty})}\leq\,& C\|\nabla_x f\|_{L_v^2(L_x^{2})}^{\frac{1}{2}}\|\Delta_x f\|_{L_v^2(L_x^{2})}^{\frac{1}{2}}\nonumber\\\label{A4}
\leq\,& C(1+t)^{-\frac{3}{2}}.
\end{align}

For $p\in [6,\infty]$, using the interpolation inequality again, it follows from \eqref{A1} and
\eqref{A3}--\eqref{A4} that
\begin{align}
\|(\rho,u,\theta)\|_{L^p}&\,\leq C\|(\rho,u,\theta)\|_{L^6}^{\zeta^{\prime}}\|(\rho,u,\theta)\|_{L^{\infty}}^{{1-\zeta^{\prime}}}\leq C(1+t)^{-\frac{3}{2}(1-\frac{1}{p})} , \\
\|f\|_{L_v^2(L_x^{p})}&\,\leq C\|f\|_{L_v^2(L_x^{6})}^{\zeta^{\prime}}\|f\|_{L_v^2(L_x^{\infty})}^{{1-\zeta^{\prime}}}\leq C(1+t)^{-\frac{3}{2}(1-\frac{1}{p})},
\end{align}
where $\zeta^{\prime}= {6}/{p} \in [0,1]$.

Now we estimate the $L^p$ norm of gradients to the solution with $p\in [2,6]$.
Using the interpolation inequality and the estimate \eqref{NJK1.17}, it follows that
\begin{align}
\|\nabla_x(\rho,u,\theta)\|_{L^p}&\,\leq C\|\nabla_x(\rho,u,\theta)\|_{L^2}^{\eta}\|\Delta_x(\rho,u,\theta)\|_{L^{2}}^{{1-\eta}}
\leq C(1+t)^{-\frac{3}{2}(\frac{4}{3}-\frac{1}{p})} , \\
\|\nabla_x f\|_{L_v^2(L_x^{p})}&\,\leq C\|\nabla_x f\|_{L_v^2(L_x^{2})}^{\eta}\|\Delta_x f\|_{L_v^2(L_x^{2})}^{1-\eta}\leq C(1+t)^{-\frac{3}{2}(\frac{4}{3}-\frac{1}{p})},
\end{align}
where $\eta= {6-p}/{2p} \in [0,1]$.
 Hence the whole proof of Corollary \ref{Cor1} is completed.   
 \end{proof}

\section{The periodic case} 

In this section, we show study the Cauchy problem \eqref{G1}--\eqref{G5} in  spatial periodic 
domain $\Omega:=\mathbb{T}^3$.  From the original system \eqref{I1},  we    obtain 
the following conservation laws by direct calculation
\begin{align*}
&\frac{\mathrm{d}}{\mathrm{d} t} \iint_{\mathbb{R}^3\times\mathbb{T}^3} F \mathrm{d} x \mathrm{d} v=0,\\
& \frac{\mathrm{d}}{\mathrm{d} t}\left(\int_{\mathbb{T}^3} n u \mathrm{d} x+\iint_{\mathbb{R}^3\times\mathbb{T}^3} v F \mathrm{d} x \mathrm{d} v\right)=0, \\
&\frac{\mathrm{d}}{\mathrm{d} t} \int_{\mathbb{T}^3} n \mathrm{d} x=0,\\
& \frac{\mathrm{d}}{\mathrm{d} t}\left(\int_{\mathbb{T}^3} n\Big(C_\textsl{v} \tilde{\theta}+\frac{1}{2}|u|^2\Big)
 \mathrm{d} x
+\iint_{\mathbb{R}^3\times\mathbb{T}^3} \frac{|v|^2}{2} F \mathrm{d} x \mathrm{d} v\right)=0.
\end{align*}
Thus, by the assumptions in  Theorem \ref{T1.2}, we obtain, for all $t \geq 0$, that  
\begin{align}\label{G4.1}
& \int_{\mathbb{T}^3} a \mathrm{d} x=0, \quad 
 \int_{\mathbb{T}^3} \rho \mathrm{d} x=0, \nonumber\\
&\int_{\mathbb{T}^3}(b+(1+\rho) u) \mathrm{d} x=0,\nonumber \\
& \int_{\mathbb{T}^3}(1+\rho)\left(\theta+\frac{1}{2}|u|^2\right)+\frac{\sqrt{6}}{2} \omega \mathrm{d} x=0.
\end{align}

 Since the proofs of local existence and uniqueness of strong solutions 
are   similar to   that in the whole space case, we shall omit them here for brevity. 
Thus, to complete the proof of Theorem \ref{T1.2}, 
the key step is to obtain the uniform a priori estimates of strong solutions with 
an exponential decay.   
According to the Poincaré's inequality and  \eqref{G4.1}, we obtain
\begin{align}\label{G4.2}
\|a\|_{L^2} & \leq C\|\nabla_x a\|_{L^2},\\
\label{G4.3}
\|\rho\|_{L^2}&\leq C\|\nabla_x \rho\|_{L^2}, \\
\label{G4.4}\|u+b\|_{L^2} & \leq\|b+u+\rho u\|_{L^2}+\|\rho u\|_{L^2} \nonumber\\
& \leq C\|\nabla_x(b+u+\rho u)\|_{L^2}+\|u\|_{L^{\infty}}\|\rho\|_{L^2} \nonumber\\
& \leq C\|\nabla_x(b, u)\|_{L^2}+C\|u\|_{H^2}\|\nabla_x \rho\|_{L^2}+C\|\rho\|_{H^2}\|\nabla_x u\|_{L^2}, \\
\label{G4.5}\|\sqrt{6} / 2 \omega+\theta\|_{L^2} & \leq\Big\|(1+\rho)\Big(\theta+\frac{1}{2}|u|^2\Big)+\frac{\sqrt{6}}{2} \omega\Big\|_{L^2}+\|\rho \theta\|_{L^2}+\left\|\rho|u|^2\right\|_{L^2}\nonumber \\
& \leq C\left(\|\theta\|_{H^2}+\|u\|_{H^2}^2\right)\|\nabla_x \rho\|_{L^2}+C\|\nabla_x \omega\|_{L^2}\nonumber \\
&\quad+C\left(1+\|\rho\|_{H^2}\right)\left(\|\nabla_x \theta\|_{L^2}+\|u\|_{H^2}\|\nabla_x u\|_{H^2}\right), \\
\label{G4.6}
\|\omega\|_{L^2}+\|\theta\|_{L^2} & \leq C\big(\|\sqrt{6} / 2 \omega+\theta\|_{L^2}+\|\sqrt{2} \omega-\sqrt{3} \theta\|_{L^2}\big) .
\end{align}

The energy functionals $\mathcal{E}(t)$ and the corresponding dissipation rate functional $\mathcal{D}(t)$  
are defined  as same as the whole space case.  Taking the similar arguments to that given  before, we obtain 
 that
\begin{equation}\label{G4.7}
\frac{\mathrm{d}}{\mathrm{d} t} \mathcal{E}(t)
+\lambda_{26} \mathcal{D}(t) \leq C\left(\mathcal{E}^{\frac{1}{2}}(t)+\mathcal{E}(t)
+\mathcal{E}^2(t)\right) \mathcal{D}(t), 
\end{equation}
for some $\lambda_{26}>0$.

Nowe we define
\begin{align*}
\mathcal{D}_{\mathbb{T}}(t):=\mathcal{D}(t)
+\tau_9\left(\|a\|_{L^2}^2+\|\rho\|_{L^2}^2\right)+\tau_{10}\|b+u\|_{L^2}^2+\tau_{11}\|\sqrt{6} / 2 
\omega+\theta\|_{L^2},
\end{align*}
where $\tau_9, \tau_{10}$ and $\tau_{11}$ are sufficient small. Note that
$$
\mathcal{D}_{\mathbb{T}}(t) \sim \sum_{|\alpha|+|\beta| \leq 2}\left\| \partial_x^\alpha\partial_v^\beta\{\mathbf{I}-\mathbf{P}\} f\right\|_\nu^2+\|(a, b, \rho, \omega)\|_{H^2}^2+\|(u, \theta)\|_{H^3}^2 .
$$
Combining \eqref{G4.2}--\eqref{G4.7} together, we deduce that
\begin{equation}
\frac{\mathrm{d}}{\mathrm{d} t} \mathcal{E}(t)+\lambda_{27} \mathcal{D}_{\mathbb{T}}(t) \leq C\left(\mathcal{E}^{\frac{1}{2}}(t)+\mathcal{E}(t)+\mathcal{E}^2(t)\right) \mathcal{D}_{\mathbb{T}}(t),
\end{equation}
for some $\lambda_{27}>0$.
Thus,
\begin{equation}
\frac{\mathrm{d}}{\mathrm{d} t} \mathcal{E}(t)+\lambda_{28} \mathcal{D}_{\mathbb{T}}(t) \leq 0,
\end{equation}
for some $\lambda_{28}>0$.
Since $\mathcal{E}(t)\leq C\mathcal{D}_{\mathbb{T}}(t) $, the above inequality implies that  
\begin{equation*}
\frac{\mathrm{d}}{\mathrm{d} t} \mathcal{E}(t)+\lambda_{29} \mathcal{E}(t) \leq 0 ,
\end{equation*}
for some $\lambda_{29}>0$.
Using Gronwall's inequality, we obtain
\begin{equation*}
\mathcal{E}(t) \leq e^{-\lambda_{29} t} \mathcal{E}(0), 
\end{equation*}
for all $t\geq 0$.
Therefore, we  complete the proof of Theorem \ref{T1.2}.

\medskip
{\bf Acknowledgements:} 
  Li and   Ni  are supported by NSFC (Grant Nos. 12331007, 12071212).  
And Li is also supported by the ``333 Project" of Jiangsu Province.
Wu is supported by Jiangsu Province ordinary university natural sciences
research project (Grant No. 24KJD110004).


\end{document}